\documentclass{amsart}
\usepackage{amsmath,amsthm,amssymb}
\usepackage[dvipsnames]{xcolor}
\usepackage{graphicx} 
\usepackage{color}
\usepackage{tikz}
\usepackage{tikz-cd}
\usepackage[all]{xy}
\usetikzlibrary{cd,positioning}

\usepackage{ytableau}

\usepackage[mathscr]{euscript}

\usepackage[colorinlistoftodos]{todonotes}

\newcommand{\A}{{\mathbb{A}}}
\newcommand{\Ad}{\mathrm{Ad}}
\newcommand{\ad}{\mathrm{ad}}

\newcommand{\af}{{\mathrm{af}}}

\newcommand{\al}{\alpha}

\newcommand{\Aut}{\mathrm{Aut}}
\newcommand{\Ax}{{\tilde{\mathbb{A}}_\af}}

\newcommand{\Bv}{B^\vee} 
\newcommand{\C}{\mathbb{C}}
\newcommand{\Cinf}{{C_\infty}}
\newcommand{\cl}{\mathrm{cl}}

\newcommand{\cC}{\mathcal{C}}
\newcommand{\cF}{\mathcal{F}}
\newcommand{\cFd}{\cF^\dagger}

\newcommand{\cZ}{\mathscr{Z}}
\newcommand{\cZGv}{\cZ_{G^\vee}}
\newcommand{\eps}{\varepsilon}

\newcommand{\fbv}{\mathfrak{b}^\vee}
\newcommand{\fg}{\mathfrak{g}}
\newcommand{\fgv}{\mathfrak{g}^\vee}

\newcommand{\fh}{\mathfrak{h}}
\newcommand{\fhv}{\fh^\vee}

\newcommand{\fso}{\mathfrak{so}}
\newcommand{\fsp}{\mathfrak{sp}}

\newcommand{\Gr}{\mathrm{Gr}}
\newcommand{\Gv}{G^\vee} 

\newcommand{\hq}{\hat{q}}

\newcommand{\hLa}{\hat{\Lambda}}
\newcommand{\hGn}{\hat{\Gamma}_{(n)}}
\newcommand{\Hom}{\mathrm{Hom}}
\newcommand{\hP}{\hat{P}}
 
\newcommand{\hz}{\hat{z}}  

\newcommand{\II}{I_Z}

\newcommand{\la}{\lambda}
\newcommand{\La}{\Lambda}

\newcommand{\nilA}{\A_\af}  
\newcommand{\om}{\varpi}

\newcommand{\pair}[2]{\langle #1\,,\,#2\rangle}
\newcommand{\Pet}{Z_{\nilA}(S)} 
\newcommand{\Petx}{{Z}_{\Ax}(S)}
\newcommand{\Ps}{\mathscr{P}}
\newcommand{\pt}{\mathrm{pt}}

\newcommand{\Quo}{\mathscr{A}_n}
\newcommand{\Quom}{\Quo^{\mathfrak{m}}}
\newcommand{\reg}{\Delta}

\newcommand{\Schub}{\mathfrak{S}}
\newcommand{\sh}[1]{D(#1)}
\newcommand{\SL}{\mathrm{SL}}
\newcommand{\SO}{\mathrm{SO}}
\newcommand{\SP}{\mathscr{SP}}
\newcommand{\Sp}{\mathrm{Sp}}
\newcommand{\Sym}{\mathrm{Sym}}

\newcommand{\trace}{\mathrm{tr}}
\newcommand{\trsp}[1]{{}^t\!{#1}} 

\newcommand{\Wafgr}{W_\af^0} 
 
\newcommand{\Wxgr}{{\Wx}^0} 
\newcommand{\Wx}{\tilde{W}_\af}
\newcommand{\Z}{\mathbb{Z}}

\newcommand{\jj}{\overline{j}}
\newcommand{\ii}{\overline{i}}
\newcommand{\pp}{\overline{p}}

\newcommand{\nn}{\overline{n}}
\newcommand{\kk}{\overline{k}}

\def\ket#1{\vert\mathord{ #1 }\rangle}
\def\bra#1{\langle\mathord{ #1 }\vert}

\newtheorem{thm}{Theorem}[section]
\newtheorem{lem}[thm]{Lemma}
\newtheorem{prop}[thm]{Proposition}
\newtheorem{cor}[thm]{Corollary}

\theoremstyle{remark}
\newtheorem{rem}[thm]{Remark}
\newtheorem{defn}[thm]{Definition}
\newtheorem{ex}[thm]{Example}

\usepackage{chngcntr}

\numberwithin{equation}{section}


\title[Equivariant homology of affine Grassmannian]{Equivariant homology of the symplectic affine Grassmannian and dual affine Schur $P$-functions}

\author[T.~Ikeda]{Takeshi Ikeda}
\address{Department of Pure and Applied Mathematics, School of Fundamental Science and
Engineering, Waseda University, 3-4-1 Okubo, Shinjuku-ku, Tokyo 169-8555, Japan}
\email{gakuikeda@waseda.jp}

\author[S.~Iwao]{Shinsuke Iwao}
\address{Faculty of Business and Commerce, Keio University, 4–1–1 Hiyosi, Kohoku-ku, Yokohama-si, Kanagawa
223-8521, Japan}
\email{iwao-s@keio.jp}

\author[M.~Shimozono]{Mark Shimozono}
\address{Department of Mathematics, 460 McBryde Hall, Virginia Tech, 225 Stanger St., Blacksburg, VA, 24061, USA}
\email{mshimo@math.vt.edu}

\begin{document}
\begin{abstract} We study the torus-equivariant homology 
$
H_*^T(\Gr_G)$ of the affine Grassmannian $
\Gr_G$, where $
G=\Sp_{2n}(\C)$ is the symplectic group. This homology admits a natural ring structure and a Schubert basis, giving rise to a well-defined  Schubert calculus. We  realize  
$H_*^T(\Gr_
G)$ in terms of symmetric functions. Our first main result introduces a new family of symmetric functions, called the \emph{dual affine Schur $P$-functions}, which represent the Schubert classes. These functions are defined through the action of the affine nil-Hecke algebra, and specialize, in the stable limit as $n\to \infty$, to the dual factorial 
$P$-functions of Nakagawa and Naruse.
Our second main result gives a precise comparison between this symmetric function model and the geometric construction of 
$H_*^T
 (\Gr_G
 )$ due to Ginzburg and Peterson, which identifies it with a coordinate ring of a centralizer family in the Langlands dual group. 
\end{abstract}

\maketitle

\tableofcontents

\section{Introduction}
Let \( G \) be a simple, simply connected complex linear algebraic group, and let \( T \subset G\) be a maximal torus of \( G \). We are interested in the torus-equivariant homology group \( H_*^T(\Gr_G) \) of the affine Grassmannian \( \Gr_G \). 
$S$-Hopf algebra and admits a Schubert basis over the equivariant cohomology ring 
$S = H_T^*(\mathrm{pt})$, the equivariant cohomology ring of a point. Although the product structure differs from the classical intersection product, ``Schubert calculus'' remains well-defined in this setting:  the multiplication of Schubert classes is governed by structure constants that are $S$-linear and retain a rich combinatorial interpretation. 
A fundamental result of Peterson~\cite{Pet} (see also~\cite{LS:Acta}) establishes a deep connection: the Schubert calculus on $H_*^T(\Gr_G)$ is essentially equivalent to that in the torus-equivariant quantum cohomology ring $QH_T^*(G/B)$ of the flag variety. This identification plays a central role in linking the geometry of affine Grassmannians and the quantum geometry of flag varieties.


This paper focuses on type C, where $G=\Sp_{2n}(\C)$ is the symplectic group; see \cite{Lam} and \cite{LLMSSZ} for known results in type A, and \cite{Pon} for types B and D. In the non-equivariant case, Lam, Schilling, and Shimozono \cite{LSS:C} realized the homology ring using symmetric functions, introducing the type C  analogue of \( k \)-\emph{Schur functions}, which correspond to Schubert classes. 

The first main result of this paper is the realization of the torus-equivariant homology ring \( H_*^T(\Gr_G) \) in terms of symmetric functions. In this setting, the special symmetric functions corresponding to Schubert classes can be defined via the action of the affine nil-Hecke algebra. We introduce a new class of functions, which we call the \emph{dual affine Schur \( P \)-functions}. These functions serve as an equivariant version of the type  C  analogue of \( k \)-\emph{Schur functions} for the symplectic group introduced in \cite{LSS:C}.
When \( n \), the rank of \( G \), is sufficiently large, this function coincides with the \emph{dual factorial \( P \)-function} introduced by Nakagawa and Naruse (\cite{NN}, \cite{NN:Uni}). 

The second main result is a detailed realization of the construction of $ H_*^T(\Gr_G)$ due to Ginzburg \cite{Gi} and Peterson \cite{Pet}. They showed that \( H_*^T(\Gr_G) \) is isomorphic to the coordinate ring of a centralizer family for the Langlands dual group \( G^\lor \)(see also Yun and Zhu \cite{YZ} for the case of integral coefficients). 
This viewpoint can be interpreted in connection with the type B Toda lattice and plays a crucial role in establishing correspondences with the quantum cohomology ring.

\subsection{Dual affine Schur $P$-functions}
\label{ssec:affine dual double Schur P}
The coefficient ring $S=H_T^*(\pt)$ of the equivariant homology $H_*^T(\Gr_G)$ is identified with 
the polynomial ring $\C[a_1,\ldots,a_n]$, where 
$\{a_i\}_{i=1}^n$ is a $\Z$-basis of the weight lattice $P$ of $\fg=\mathfrak{sp}_{2n}(\C),$ the Lie algebra of $G=\Sp_{2n}(\C).$
Let $\alpha_i = a_i - a_{i+1}$ for $1\le i\le n-1$ and
$\alpha_n = 2a_n$ be the simple roots.
We also set $\alpha_0=-\theta$, where $\theta=2a_1$ is the highest root.
This is the classical projection of the affine simple root, whose value is $\delta-\theta$ where $\delta$ is the affine null root.

Let $W$ and $W_\af$ be the Weyl group and the affine Weyl group of $\fg.$ 
$W$ is a subgroup of $W_\af$ and $W_\af$
is identified with the semidirect product $W\ltimes Q^\lor$, where 
$Q^\lor $ is the coroot lattice of $\fg.$ We have $Q^\lor=\bigoplus_{i=1}^n \Z\eps_i\subset \Hom_\Z(P,\Z)$ where $\{\eps_i\}_{i=1}^n$
is the dual basis of $\{a_i\}_{i=1}^n$. 
$W_\af$ acts on $S$ via
the projection 
$W_\af\rightarrow W$ given by $t_\la u\mapsto u$ for $\la\in Q^\vee$ and $u\in W$; this is called the \emph{level zero action}. 
Let $\A_\af$ be the \emph{level zero affine nil-Hecke algebra} of type $\mathrm{C}_n^{(1)}$ (see \S \ref{SS:level zero extended affine nil-Hecke} for its definition). The algebra $\A_\af$ contains a Demazure type element
$$
A_i=\alpha_i^{-1}(1-s_i)\quad \text{for $i\in I_\af$},
$$
where $s_i$ is the standard generator of $W_\af.$
Since the $A_i$ satisfy the same braid relations as the $s_i$, for $w\in W_\af$, we can define $A_w\in \A_\af$. 
The algebra 
$\A_\af$ is a free left $S$-module with basis $\{A_w\mid w\in W_\af\}.$

Define $\hLa_S=S\otimes  \hLa$, where $\hLa$ is the graded completion of the ring $\Lambda$ of symmetric functions in $y=(y_1,y_2,\ldots).$
The left $S$-module structure of $\hLa_S$ can be  
extended to an $\A_\af$-module structure such that 
$W$ acts on $\hLa_S$ via the natural action of $W$ on $S$, and the translation element
$t_{\eps_i}\in W_\af$
acts on $\hLa_S$ by 
multiplication by 
\begin{equation}
\Omega(a_i|y):=\prod_{j=1}^\infty\frac{1+a_iy_j}{1-a_i y_j}.
\end{equation}

Let $W_\af^0$ denote the set of minimal length representatives for the coset space $W_\af/W$.
It is a natural index set for the Schubert basis $\{\sigma_w\}_{w}$ of $H_*^T(\Gr_{\Sp_{2n}}).$ 
For $w\in \Wafgr$, define
the \emph{dual affine Schur $P$-function} by 
\begin{equation}
\hP_w^{(n)}(y|a)=A_w(1).
\end{equation}
Let $\hGn$ be the $S$-span of $\hP_{w}^{(n)}(y|a)\;(w\in \Wafgr)$. 
The first main result is
the following.
\begin{thm}\label{thm:main isom gamma}
There is an isomorphism of $S$-algebras
$$
\gamma: H_*^T(\Gr_{\Sp_{2n}})
\rightarrow \hGn
$$
such that the Schubert 
class $\sigma_w$ of 
$
H_*^T(\Gr_{\Sp_{2n}})$ corresponds to $\hP_w^{(n)}(y|a).$
\end{thm}
The specialization of $\hP_w^{(n)}(y|a)$ at $a=0$, recovers
the type C analogue of $k$-Schur function $\hP_w^{(n)}(y)$ introduced in \cite{LSS:C}.

\subsection{Relation to infinite rank case} 
\label{SS:relation to infinite rank}
In the infinite rank limit, the dual affine  Schur $P$-function is the dual factorial Schur $P$-function $\hP_\la(y|a)$ defined by Nakagawa and Naruse; this correspondence will be made precise below.


Nakagawa and Naruse \cite{NN} introduced a family of functions $\hP_\la(y|a)$ indexed by the set $\SP$ of strict partitions, where $a$ is an infinite sequence $(a_1,a_2,\ldots).$
Let $\C[a]$ be the polynomial ring in the variables $a_i$, and $\hat{\Gamma}_{\C[a]}:=\C[a]\otimes\hat{\Gamma}$, where $\hat{\Gamma}$ is the completion of the ring $\Gamma$ of Schur's $Q$-functions with respect to the natural grading.
The functions $\hP_\la(y|a)$, called the \emph{dual factorial Schur $P$-functions}, form a $\C[a]$-basis of $\hat{\Gamma}_{\C[a]}$, and are identified with the Schubert basis of in the equivariant homology of the infinite dimensional Langrangian Grassmannian $LG(\infty).$

There is a bijection $\Wafgr\to \Ps_C^n$ denoted $w\mapsto \la_w$ \cite[Lemma 24]{EE} \cite[Lemma 5.3]{LSS:C} where $\Ps_C^n$ is the set of partitions $\lambda=(\lambda_1,\lambda_2,\ldots)$ whose part sizes are at most $2n$, which have at most one part of each size $k$ for $1\le k\le n$.

For $i\ge 1$, define $a_i^{(n)}$ by 
\begin{align}\label{E:an}
(a_1^{(n)},\ldots,a_{2n}^{(n)})&=
(a_1,\ldots,a_n,-{a}_n,\ldots,-{a}_1) \\
\label{E:an periodic}
a_{i+2nk}^{(n)}&=a_i^{(n)}\qquad\text{for $k\ge 1$.}
\end{align}

\begin{thm}\label{thm:special}
If $w\in W_\af^0$ satisfies  $\ell(w)\le 2n$ then
$\la_w$ is a strict partition and 
$$
\hP_{w}^{(n)}(y|a)=
\hP_{\la_w}(y|a^{(n)}).
$$

\end{thm}

Similar results are known for the nonequivariant type A case 
by Lapointe and Morse in \cite[Property 39]{LM}.
See \cite[Proposition A.4]{ISY} for the corresponding result 
for the equivariant $K$-homology.



\subsection{Ginzburg–Peterson construction}
The isomorphism constructed in Theorem \ref{thm:main isom gamma} is most naturally understood in the context of the family of centralizers 
$\cZGv$
introduced by Peterson (see \S\ref{sec:cent} for the definition).
This is an affine variety over \( \mathfrak{h} \),  
defined as a closed subvariety of \( \Bv \times \mathfrak{h} \), where \( \Bv \) is a Borel subgroup of the dual group \( \Gv = \SO_{2n+1}(\mathbb{C}) \) of \( G \).  
Let \( z_{ij} \) be the coordinate functions on the Borel subgroup \( B^\lor \).  

Note that the translation element $t_\lambda\;(\lambda\in P)$ can be considered as an element in $H_*^T(\Gr_{\Sp_{2n}(\mathbb{C})})$ (see \S \ref{ssec:affnil}). 

\begin{thm}\label{thm:centralizer}
The following diagram commutes, where all maps \( \phi, \beta, \gamma \) are isomorphisms of $S$-algebras:
\[
\begin{tikzcd}
H_*^T(\Gr_{\Sp_{2n}(\mathbb{C})})\arrow[d,"\gamma", swap]&   \\
\hGn    & \mathbb{C}[\cZGv]\arrow[l,"\beta"]\arrow[lu, "\phi", swap]
\end{tikzcd}
\]
These maps satisfy  
\[
\phi(z_{ii}^{\pm 1}) = t_{\pm \eps_i},\quad
\gamma(t_{\pm \eps_i})
= \Omega(\pm a_i|y)
= \beta(z_{ii}^{\pm 1}).
\]
Each of \( H_*^T(\Gr_{\Sp_{2n}(\mathbb{C})}) \), \( \hGn \), and \( \mathbb{C}[\cZGv] \) is equipped with the structure of a left \( \nilA \)-module, and each of the maps \( \phi, \beta, \gamma \) is an isomorphism of \( \nilA \)-modules.
\end{thm}

Our arguments are based on Peterson's theory (\cite{Pet}) on the affine nil-Hecke algebra.
Let $\Pet$ be the centralizer algebra of $S$ in $\nilA.$ It is known that $\Pet$ is a commutative $S$-algebra. 
A fundamental result due to Peterson is that $\Pet$ is isomorphic to $H_*^T(\Gr_G)$ as $S$-algebras.

In order to construct the map $\phi$, it is necessary to prove that $\C[\cZGv]$ is flat over $S.$ 
In fact, Peterson \cite{Pet} proved the flatness by using sheaf cohomology of a line bundle on closely related varieties, the so-called Peterson varieties. 
Our proof of the flatness is based on a concrete form of the generators of the defining ideal of $\cZGv.$ 

The construction of $\beta$ is most naturally understood in the context of the Toda lattice. In fact, the map arises in Kostant's solution of the type $\mathrm{B}_n$ Toda lattice. 
More details on the Toda lattice and equivariant quantum cohomology will be discussed in forthcoming papers.

\subsection{Pieri classes}
For \( 1 \leq i \leq 2n \), define \( \rho_i \in \Wafgr \) as follows:  
\begin{equation}
\rho_i =
\begin{cases}
    s_{i-1} \cdots s_1 s_{0} & \text{for } 1 \leq i \leq n, \\
    s_{2n+1-i} \cdots s_{n-1} s_{n}
    s_{n-1} \cdots s_1 s_0 & \text{for } n+1 \leq i \leq 2n.
\end{cases}
\end{equation}
It was shown in \cite{LSS:C} that  
\[
H_*(\Gr_{\Sp_{2n}}) \cong \mathbb{Z}[P_1(y), P_2(y), \ldots, P_{2n}(y)],
\]
where \( P_i(y) \) is the Schur $P$ function corresponding to single row partition $(i)$, which corresponds to the non-equivariant Schubert basis associated with \( \rho_i \).  
We obtain the equivariant analogue of this result.  

\begin{thm}
\label{thm:Pieri classes}
\( H_*^T(\Gr_{\Sp_{2n}}) \) is generated by \( \sigma_{\rho_i} \) for \( 1 \leq i \leq 2n \) as an \( S \)-algebra.
\end{thm}
We derive this result as a corollary to Theorem \ref{thm:centralizer}.
Via the isomorphism \( \gamma \), the special Schubert class  
\( \sigma_{\rho_i} \) is identified with \( \hP_{\rho_i}^{(n)}(y|a) \), which is equal to  $\hP_{(i)}(y|a^{(n)})$ by Theorem \ref{thm:special}.
This yields a concrete algebraic model
\[
S[\hP_{1}(y|a^{(n)}), \ldots, \hP_{2n}(y|a^{(n)})]
\]
for \( H_*^T(\Gr_{\Sp_{2n}}) \).
The following fundamental  questions remain open:
\begin{itemize}
    \item Pieri rule:
 Find an explicit formula for the expansion of $\hP_{\rho_i}^{(n)}(y|a)\hP_w^{(n)}(y|a)$
   in terms of the Schubert basis.
    \item Giambelli formula: 
    Express $\hP_\la^{(n)}(y|a)$ explicitly as a polynomial in  $\hP_{\rho_i}^{(n)}(y|a)$ for $1\le i\le 2n.$
\end{itemize}
These problems lie at the heart of understanding the Schubert calculus in equivariant $K$-homology of $\Gr_{\Sp_{2n}}$.
As for the equivariant Pieri rule, there is a conjecture by Lam and Shimozono \cite[Conjecture 14]{LS:pieri}.


\subsection*{Related works}
A type A analogue of the current article was partly 
studied in \cite{LS:double Kostant}.
In \cite{ISY}, the
equivariant $K$-homology of type A affine Grassmannian was studied.
Based on this, the connection with the integrable system called the relativistic Toda lattice was studied in \cite{IINY}.

The current paper should be extended to equivariant $K$-homology of $\Gr_{\Sp_{2n}}$ and will be discussed in a forthcoming paper \cite{ISY2}. 
$H^T_*(\Gr_{\Sp_{2n}}(\C))$ is known to have an $S$-Hopf algebra structure. 
The $T$-equivariant cohomology $H_T^*(\Gr_{\Sp_{2n}}(\C))$ is naturally the dual to $H_*^T(\Gr_{\Sp_{2n}}(\C))$ as an $S$-Hopf algebra. More detailed discussions about the duality  with $H_T^*(\Gr_{\Sp_{2n}}(\C))$ will be pursued elsewhere.

\subsection*{Acknowledgements}
Special thanks to Kazuma Shimomoto for helpful discussions on commutative ring theory. We also thank Noriyuki Abe, Koushik Brahma, Satoshi Naito, Thomas Lam, Hiroshi Naruse, and Kohei Yamaguchi for fruitful discussions.
We thank Vladimir N.~Ivanov for sharing with us his unpublished result on $\hat{P}$-functions.
This work was supported by JSPS KAKENHI Grant Numbers 23K03056 and 22K03239.

\subsection*{Organization}

In \S \ref{sec:Phat}, we review the dual factorial $P$-functions introduced by Nakagawa and Naruse, and introduce the $\hat{q}$-functions, which will be used later in \S \ref{sec:cent}.  
In \S \ref{sec:Ginzburg-Peterson}, we review basic results on the level-zero affine nil-Hecke algebra and prove a general factorization formula for anti-dominant translation elements (\S \ref{ssec:factor}).  
In \S \ref{sec:C}, we fix notation and conventions for type~C.  In \S \ref{sec:Sym}, we prove our first main result, Theorem~\ref{thm:main isom gamma} (Theorem~\ref{thm:gamma}), and in \S \ref{ssec:extended_sym} we provide a realization of the equivariant $K$-homology of $\Gr_{\mathrm{PSp}_{2n}}$.  
Section~\ref{ssec:small} is devoted to the proof of Theorem~\ref{thm:special}, and in \S \ref{ssec:factorization} we establish the factorization formula for $\hat{P}_w^{(n)}(y|a)$ (Theorem~\ref{thm:hat P factor}).  
In \S \ref{sec:cent}, we prove Theorem~\ref{thm:centralizer}, and in \S \ref{sec:flat} we show that the coordinate ring of $\cZGv$ is flat over~$S$.  
Finally, in \S \ref{sec:Fermion}, we develop a free fermion formalism for $\hat{P}_\lambda(y|a)$ and, as an application, prove Proposition~\ref{prop:hatq_to_hatP} concerning $\hat{q}$-functions, which is used in the proof of Theorem~\ref{thm:Pieri classes}.

\section{Dual factorial $P$-functions}
\label{sec:Phat}
In this section, we recall the definition of Ivanov's factorial $Q$-functions, which represent the Schubert basis of the equivariant cohomology of the infinite Lagrangian Grassmannian $LG(\infty)$.
We also review the construction of $\hat{P}$-functions due to Nakagawa and Naruse \cite{NN:Uni}; these are the Schubert basis for the equivariant homology of $LG(\infty)$.

\subsection{Root system $\mathop{\mathrm{C}}_\infty$ and strict partitions}
\label{ssec:Cinf}
Let $\mathrm{C}_\infty$ be the following Dynkin diagram of infinite rank:
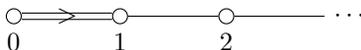
\begin{center}
\begin{tikzpicture}[node distance=1.5cm and 1.2cm, auto]
    \node[draw, circle, inner sep=2pt, label=below:{$0$}] (0) {};
    \node[draw, circle, inner sep=2pt, label=below:{$1$}, right=of 0] (1) {};
    \node[draw, circle, inner sep=2pt, label=below:{$2$}, right=of 1] (2) {};
    \node[right=of 2] (dots) {\(\cdots\)};
    
    \draw[double distance=2pt] (0) -- (1) node[midway, yshift=-1.5ex] {\(>\)};
    
    \draw (1) -- (2);
    \draw (2) -- (dots);
\end{tikzpicture}
\end{center}
Its weight lattice is by definition $P_\infty = \bigoplus_{k\in \Z_{>0}} \Z a_k.$
Let $I_\infty=\{0,1,2,\ldots\}$ be the index set of $\mathrm{C}_\infty.$ For $i\in I_\infty$, define simple roots 
\begin{equation}
    \alpha_0=-2a_1, \quad \alpha_i=a_i-a_{i+1}\quad\text{ for $i>0$}.
\end{equation} 
Let $P_\infty^*=\bigoplus_{k\in \Z_{>0}}\Z
\eps_k$ with pairing $P_\infty^*\times P_\infty\rightarrow \Z$ defined by  $\langle \eps_i,a_i\rangle=\delta_{ij}.$
The simple coroots are defined by 
\begin{equation}
\alpha_0^\lor=-\eps_1,\quad\alpha_i^\lor=\eps_i-\eps_{i+1}\quad \text{for $i>0$}.
\end{equation}
Let $W_\infty$ be the group generated by the following elements $s_i\in \Aut(P_\infty)$:
\begin{align}
  s_i(v) = v - \langle\alpha_i^\vee,v\rangle \alpha_i\qquad\text{for $i\in I_\infty, \quad v\in P_\infty$.}
\end{align}
We have
\begin{align*}
 s_0(a_1) &= -a_1, \quad s_0(a_k) = a_k \; \text{($k\ge 1$), and for $i>0,$} \\
s_i(a_i) &= a_{i+1}, \quad s_i(a_{i+1}) =a_i ,\quad 
s_i(a_k) = a_k \;\text{($k\ge 1,\;k\notin\{i,i+1\}$)}.
\end{align*}

A \emph{strict partition} $\la\in\SP$ is a strictly decreasing sequence of positive integers $\lambda=(\lambda_1>\ldots>\lambda_l>0)$; its length $\ell$ is denoted $\ell(\la)$.
Given $\la\in\SP$ let $\sh{\la}$ denote its \emph{shifted diagram}: $\sh{\la}=\{(i,j) \mid 1\le  i\le j < i + \la_i\}$. 
$\SP$ is a graded poset under containment of shifted diagrams. We denote the covering relation by $\mu\gtrdot\la$, that is, $\mu>\la$ and there is no $\nu$ such that $\mu>\nu>\la.$
Let $s_i \la = \mu$ if either $\mu\gtrdot \la$ or $\mu \lessdot \la$ and in either case the box giving the difference of the shifted shapes, is on diagonal $i$, where the diagonal index of a box $(r,c)$ in the $r$-th row and $c$-th column is $r-c$. 
If there is no $\la$-addable nor $\la$-removable box on diagonal $i$ we set $s_i \la=\la$.

This defines a transitive $W_\infty$-action on $\SP$.
The empty partition has a stabilizer generated by $s_i$ for $i>0$. This induces a bijection $\SP\to W_\infty^0$ denoted $\la\mapsto w_{\la}$, where $W_\infty^0=\{w\in W_\infty\mid ws_i>w\text{ for all $i>0$ } \}$. 
For $\la\le\mu$ with $\la,\mu\in \SP$, define $w_{\mu/\la} = w_{\mu} w_{\la}^{-1}$.

\begin{ex} Let $\mu=(5,3)$ and $\la=(3,1)$.
Then $w_{\mu}=s_2s_1s_0\cdot s_4s_3s_2s_1s_0 = s_2s_1s_4s_3\cdot s_0s_2s_1s_0$, $w_{\la}=s_0s_2s_1s_0$, and $w_{\mu/\la} = s_2s_1s_4s_3$:
\[
\begin{ytableau}
*(green) s_0 & *(green)s_1 &*(green) s_2 & s_3 & s_4 \\
\none & *(green) s_0 & s_1 & s_2\\
\end{ytableau}
\]
Observe that if $\ell(\la)=\ell(\mu)$ then $w_{\mu/\la}$ is a permutation, that is, an element of $W=\langle s_1,s_2,\ldots\rangle$.
\end{ex}

\subsection{Schur $Q$-functions}
\label{SS:Schur Q}
We refer \cite[Chap. III, \S 8]{Mac} for the basic properties of the Schur $Q$-functions.
 
Let $\Gamma(x):=\C[Q_1(x),Q_2(x),\ldots]$ be the ring of Schur $Q$-functions, where $Q_i=Q_i(x)$ is the symmetric function in $x=(x_1,x_2,\ldots)$ defined by 
\begin{equation}
\Omega(u|x) := \prod_{i=1}^\infty
    \frac{1+x_i u}{1-x_i u}=
\sum_{k=0}^\infty Q_k(x)u^k.
\end{equation}
For $\la\in\SP$ the Schur $Q$-function $Q_\lambda(x)\in \Gamma(x)$ can defined as $Q_\lambda(x;-1)$, the Hall-Littlewood symmetric function with $t=-1$ (see \S~\ref{sec: factorial Q} below for more direct construction of $Q_\la(x)$). Note that $Q_i(x)=Q_{(i)}(x)$, where $(i)$ is the strict partition of the single part $i$. 
We have $\Gamma(x)=\bigoplus_{\la\in\SP} \C Q_\lambda(x)$ and $\Gamma(x) = \C[p_1,p_3,p_5,\dotsc]$, where $p_{i}\;(i\ge 1)$ is the power sum of degree $i$.

Let $\Gamma(y)\times \Gamma(x)\rightarrow \C$ be the bilinear form defined by $\langle P_\lambda(y),Q_\mu(x)\rangle=
\delta_{\lambda\mu}$, where $P_\lambda(y)=2^{-\ell(\la)}Q_\lambda(y)$.
Equivalently, we have
\begin{equation}\label{E:PQ form}
\sum_{\lambda\in \SP}P_\lambda(y)Q_\la(x)=
\prod_{i,j}\frac{1+x_iy_j}{1-x_iy_j}.
\end{equation}

\subsection{Factorial $Q$-functions}
\label{sec: factorial Q}
For $\lambda\in\SP$ and a positive integer $n$, define
\begin{equation}
    Q_\lambda(x_1,\ldots,x_n|a)
    :=\frac{1}{(n-\ell(\lambda))!}
    \sum_{w\in S_n}
    w\left(
    \prod_{i=1}^{\ell(\lambda)}(\!(x_i|a)\!)^{\lambda_i}
    \prod_{i=1}^{\ell(\lambda)}\prod_{j=i+1}^n\frac{x_i+x_j}{x_i-x_j}
    \right),
\end{equation}
where the symmetric group $S_n$ acts by permuting the variables $x_1,\ldots,x_n$, and 
$$(\!(x|a)\!)^k:=\begin{cases}
2x & \text{for $k=1$,}\\
    2x(x-a_1)\cdots (x-a_{k-1}) & \text{for $k\ge 2$.}
\end{cases}$$ 
$Q_\la(x|a)$ is a symmetric polynomial in $x$.
The symmetric function, called the \emph{factorial Schur $Q$-function}, was introduced by Ivanov in \cite[Theorem 9.1]{Iv}. 
\begin{rem}\label{rem:Q}[Indexing and sign conventions for the equivariant parameters]
We follow the convention in \cite{IMN} for the index of $a_i$ variables in $Q_\lambda(x|a)$. 
In \cite{Iv}, it was assumed $a_1=0$ for the most part. In contrast, we do not use $a_1$ of \cite{Iv} and instead relabel $a_{i+1}$ of \cite{Iv} by $a_i$ for $i\ge 1$.
Recall also that the simple roots in \cite{IMN} are given by  $\alpha_0=2t_1,\; \alpha_i=t_{i+1}-t_i$ for $i\ge 1$. Accordingly, the variable $t_i$ in \cite{IMN} corresponds to $-a_i$ here. 
 \end{rem}
A formula expressing the factorial $Q$-function as a Schur-type Pfaffian is also available (\cite{Iv}). 
It is easy to see that $Q_\lambda(x_1,\ldots,x_n|a)$ vanishes unless $\ell(\lambda)\le n$.
We have the property 
\begin{align}\label{E:Ivanov stable}
Q_\lambda(x_1,\ldots,x_n,0|a)
=Q_\lambda(x_1,\ldots,x_n|a)
\end{align}
(\cite[Proposition 1.1]{Ivanov-dim}, see also \cite[Remark 3.1]{GPGQ}). The factorial $Q$-function, denoted by  $Q_\lambda(x|a),$ is by definition the stable limit of the family of polynomials $\{Q_\lambda(x_1,\ldots,x_n|a)\}_{n=1}^\infty.$
The limit $Q_\lambda(x|a)$ is an element of $\Gamma(x|a):=\C[a]\otimes_\C\Gamma(x)$, where $\C[a]=\C[a_1,a_2,\ldots].$
Then $Q_\la(x)$ is obtained from $Q_\la(x|a)$ by setting all $a_i$ to be zero.
The factorial $Q$-functions $Q_\lambda(x|a)$ for $\la\in \SP$ with $\la_1\le n$, are known to represent the Schubert classes of the torus equivariant cohomology of the Lagrangian Grassmannian \cite{Ike}.

\begin{rem} The value $(\!(a_i|a)\!)^k$ coincides with the localization of the Schubert class $\rho_k=s_{k-1}\dotsm s_1s_0$ at the torus-fixed point $\rho_i$ for $1\le k \le i$ of the Lagrangian Grassmannian. We note that the same value also appears in the affine Grassmannian case 
(see Lemma~\ref{L:loc value}). This observation suggests a more general construction of the ``factorial Schur'' functions for other types as well.
\end{rem}



\subsection{Dual factorial $P$-functions $\hP_\la(y|a)$}
\label{SS:dual Ivanov}
Let $\hat{\Gamma}(y)$ denote the set of formal sums $f(y)=\sum_{k\ge 0}f_k(y)$ with $f_k(y)\in \Gamma$ in $y=(y_1,y_2,\ldots)$ of degree $k.$
Define $\hat{\Gamma}(y|a):=\C[a]\otimes_\C\hat{\Gamma}(y)$.
The dual Ivanov functions $\hat{P}_\lambda(y|a)\in \hat{\Gamma}(y|a)$ are defined by
\begin{equation}\label{eq:hP and Q}
\langle \hat{P}_\lambda(y|a),{Q}_\mu(x|a)\rangle
=\delta_{\lambda\mu},
\end{equation}
where $\langle \cdot,\cdot\rangle: 
\hat{\Gamma}(y|a)
\times
\Gamma(x|a)\rightarrow \C[a]
$ is the $\C[a]$-bilinear form obtained by tensoring the pairing of \S \ref{SS:Schur Q} with $\C[a]$. Equivalently we have   \begin{equation}
     \prod_{i,j}\frac{1+x_iy_j}{1-x_iy_j}=\sum_{\lambda}\hP_\lambda(y|a)Q_\lambda(x|a).\label{eq:Cauchy}
    \end{equation}
The functions $\hP_\la(y|a)$ were introduced by Nakagawa and Naruse (\cite{NN}, \cite{NN:Uni}).  

\begin{thm}[Ivanov \cite{Iv:pri}]
\label{prop: dual Ivanov}
Let $\la\in \SP.$ We have
\begin{align}
\label{E:dual Ivanov into P}
\hat P_\lambda(y|a) &= \sum_{\substack{\mu\in \SP \\
\la\subset \mu \\ \ell(\mu)=\ell(\lambda)}} 
\det\left(h_{\mu_i-\lambda_j}(a_{1},\dots,a_{\lambda_{j}})\right)_{ij}  P_\mu (y). 
\end{align}
\end{thm}
\begin{proof}
This expansion is obtained as the dual to the formula in \cite[Theorem 10.2]{Iv}, which expresses $Q_\la(x|a)$ in terms of $Q_\mu(x).$
\end{proof}

\begin{rem}\label{rem:On_ivanov_function}

We have
\begin{align}
    \det \left(h_{\mu_i-\la_j}(a_1,\dotsc,a_{\la_{j}})\right)_{ij} = \Schub_{w_{\mu/\la}}(a)
\end{align}
where $\mathfrak{S}_w(a)$ denotes the Schubert polynomial \cite{LaSc} associated with a permutation $w$. 
See \cite[Proposition 8.2]{LLS:backstable} for a similar result in type A.
\end{rem}

\begin{cor}\label{cor:hat_P}   
For $r\ge 1$, we have 
\begin{equation} \label{cor-conj}
\hat P_{r}(y|a)=\sum_{i=0}^\infty
h_i(a_1,\ldots,a_r)
P_{r+i}(y).    
\end{equation}
\end{cor}

\begin{prop}\label{prop: Omega in hP} Let $b$ be a formal parameter. Then we have
    \begin{equation}
     \Omega(b|y)= 1+\sum_{k=1}^\infty (\!(b|a)\!)^k\,\hP_k(y|a).\label{eq:Omega(b)}
    \end{equation}
    In  particular, we have for $1\le i\le n,$
    \begin{align}
    \Omega(a_i|y) &= 1 + \sum_{k=1}^i (\!(a_i|a)\!)^k\,\hP_k(y|a).
    \label{eq:Omega in special classes} 
    \end{align}
\end{prop}
\begin{proof}
We specialize \eqref{eq:Cauchy} to $(x_1,x_2,\ldots)=(b,0,\ldots)$, and note that 
    \begin{align*}
    Q_{(i)}(b,0,\ldots,0)
    =Q_{(i)}(b|a) = (\!(b|a)\!)^i,
\end{align*}
so that 
\begin{equation}
 Q_\lambda(b,0,\ldots|a)=\begin{cases}
    (\!(b|a)\!)^i
    & \text{for $\la=(i)$ with $i\ge1$} \\
    0 & \text{otherwise.}
\end{cases}
\end{equation}
The equation \eqref{eq:Omega(b)} follows. 
Hence \eqref{eq:Omega in special classes} follows, since $(\!(b|a)\!)^i$ evaluates to zero when $b=a_i$ if $i>k$.
\end{proof}


\begin{rem} Plugging in \eqref{eq:Omega(b)} for $b$ and $-b$ into $\Omega(b|y)\Omega(-b|y)=1$, one obtains quadratic relations among the $\hP_i(y|a)$. 
In particular $\hP_{2i}(y|a)$ is expressed as a polynomial of $\hP_1(y|a)$, $\hP_3(y|a)$, $\ldots$, $\hP_{2i-1}(y|a)$ with coefficients in $\C[a_1,\dotsc,a_i]$.
\end{rem}

It is known by \cite{NN:Uni} that $\hat{\Gamma}(y|a)$ is an $\C[a]$-algebra and $\{\hat{P}_\la(y|a)\mid \la \in \SP\}$ is a formal $\C[a]$-basis. 
See \cite[\S 8.1]{LLS:backstable} for analogous results in type $A.$

\begin{ex} 
The direct sum $\bigoplus_{\la\in \SP}\C[a]\hP_\la(y|a)$ is not a ring. 
For example, we have
\begin{align*}
\hP_1(y|a)^2 = \hP_2(y|a)+\sum_{i=1}^\infty\left(
\left(\prod_{j=2}^{i+1} (a_1-a_j)\right)
\hP_{i+1}(y|a)
+(\!(a_1|a)\!)^{i}\hP_{i+1,1}(y|a)\right).
\end{align*}
\end{ex}
\begin{rem}   
As far as we know, an explicit Pieri-type expansion of $ {\hP}_{(i)}(y|a)\hP_\la(y|a)$ in terms of the ${\hP}_\mu(y|a)$ is not yet known, but it remains a fundamental problem.
See \cite[\S 8]{Molev}, \cite[\S 8.1]{LLS:backstable} for analogous results in type $A$. 
\end{rem}


\subsection{Left divided difference operators}
\label{SS:left divided differences}
This subsection presents another characterization of $\hP_\la(y|a)$ in terms of the left divided difference operators.



We define an action of \( W_\infty \) on the ring \( \Gamma(x|a) \) by
\begin{align}
s_0 f(x|a) &:= f(-a_1, x \mid -a_1, a_2, a_3, \ldots), \\
s_i  f(x|a) &:= f(x \mid \ldots, a_{i+1}, a_i, \ldots) \quad \text{for } i > 0,
\end{align}
for \( f \in \Gamma(x|a) \). 
The action of $s_0$ is well-defined on $\Gamma(x|a)$ and satidfies $s_0^2=1$, because $f(x|a)\in \Gamma(x|a)$ satisfies the $Q$-\emph{cancellation property}: $$f(t,-t,x|a)=f(x|a)$$ for arbitrary $t.$

We define an action of \( W_\infty \) on the ring \( \hat{\Gamma}(y|a) \) by setting 
\begin{align}
s_0  f(y|a) &:= \Omega(a_1 \mid y)\, f(y \mid -a_1, a_2, a_3, \ldots), \\
s_i  f(y|a) &:= f(y \mid \ldots, a_{i+1}, a_i, \ldots) \quad \text{for } i > 0,
\end{align}
for \( f \in \hat{\Gamma}(y|a) \).
The action of $s_0$ is well-defined since $\Omega(a_1|y)$ is an element of $\hat{\Gamma}(y|a)$ by Proposition \ref{prop: Omega in hP}.


Consider the diagonal action of  $W_\infty$ on $\hat{\Gamma}(y|a)\hat{\otimes}_{\C[a]} \Gamma(x|a)$, the completion 

\begin{lem}\label{L:Omega invariant}
Then the element $\Delta(x,y)=\prod_{i=1}^\infty\Omega(x_i|y)$ is $W_\infty$-invariant.
\end{lem}
\begin{proof} For $i>0$,  $\Delta(x,y)$ is obviously $s_i$-invariant because $\Delta(x,y)$ does not involve the parameters $a_i$. The action of $s_0$ on the tensor factor $\hat{\Gamma}(y|a)$ produces a factor of $\Omega(a_1|y)$, and the action of $s_0$ on the left factor multiplies by $\Omega(-a_1|y)$, and the factors cancel. 
\end{proof}

Define a linear operator $A_i$ on $\hat{\Gamma}(y|a)\otimes_{\C[a]} \Gamma(x|a)$ by
\begin{align}\label{E:delta def}
A_i = \frac{1}{\alpha_i}(1-s_i)\qquad\text{for $i\in I_\infty$.}
\end{align}
Note that we have 
\begin{equation}\label{eq:delta_i s_i}
s_i A_i=-A_is_i
=A_i.
\end{equation}
This yields $A_i^2=0.$ One can also verify the type $\Cinf$ braid relations 
$$
A_0A_1A_0A_1=A_1A_0A_1A_0, \quad A_i A_{i+1}A_i=A_{i+1}A_iA_{i+1}
\quad \text{for $i\ge 1$}
$$

Recall the action of $W_\infty$ on $\SP$ defined in \S \ref{ssec:Cinf}.
Let $\chi(\text{True})=1$ and $\chi(\text{False})=0$.

\begin{prop}\label{prop:delta P-hat} For all $\mu\in\SP$ and $i\in \Z_{\ge0}$
\begin{equation}\label{E:delta P-hat}
 A_i \hP_\mu(y|a) = \chi(s_i\mu>\mu) \hP_{s_i\mu}(y|a).
    \end{equation}
\end{prop}
\begin{proof}
From \cite{IMN} (see also Remark \ref{rem:Q}) we have
\begin{align}\label{E:A on Q}
A_i Q_\la(x|a) = - \chi(s_i\la<\la)  Q_{s_i\la}(x|a)
\qquad\text{for all $\la\in\SP$ and $i\in \Z_{\ge0}.$}
\end{align}
Fix $\mu\in\SP$ and $i\in \Z_{\ge0}$. From Lemma \ref{L:Omega invariant}, we have
\begin{align*}
0&=A_i \Delta(x,y)
=\sum_{\la\in \SP}A_i (\hP_\la(y|a)Q_\la(x|a))\nonumber\\&=
\sum_{\la\in \SP}A_i \hP_\la(y|a) Q_\la(x|a)+
\sum_{\la\in \SP}s_i \hP_\la(y|a)A_iQ_\la(x|a) \nonumber\\
&=\sum_{\la\in \SP}A_i \hP_\la(y|a) Q_\la(x|a)-\sum_{\la\in \SP}\chi(s_i\la<\la) s_i \hP_\la(y|a)Q_{s_i\la}(x|a).
\end{align*}
Taking the coefficient of $Q_\mu(x|a)$ we have
\begin{align*}
  0 &= A_i \hP_\mu(y|a) - \chi(\mu<s_i\mu) s_i \hP_{s_i\mu}(y|a)
\end{align*}
as required.
\end{proof}
\begin{rem}
This result follows from \cite[Proposition 5.1]{NN:Uni} if we use Lemma \ref{L:Omega invariant}. 
\end{rem}



\subsection{$Q$-functions in time variables}
We discuss $Q$-functions in terms of the power sum symmetric functions, which are relevant in the context of the Toda lattice.

Let $t_1,t_3,\dotsc,t_{2n-1}$ be indeterminates.
We define ${q}_k(t)=q_k(t_1,t_3,\ldots,t_{2n-1})\in \C[t_1,t_3,\dotsc,t_{2n-1}]$ by the following generating function:
\begin{equation}\label{eq:generating_q}
\exp\left(2\sum_{i=1}^n t_{2i-1} u^{2i-1} \right)
=\sum_{k=0}^\infty
{q}_k(t)u^k.    
\end{equation}
If we define the degree of $t_{2i-1}$ to be $2i-1$, the polynomial ${q}_k(t)$ is homogeneous of degree $k.$
Note in the context of the Toda lattice equation of type $B_n$, $t_1,t_3,\ldots,t_{2n-1}$ play the role of the time variables.

Let $p_{2i-1}$ denote the power sum of degree $(2i-1)$. 
It is well-known (see \cite[Chap. III, \S 8]{Mac}) that
\begin{equation}\label{eq:Pi}
    \C[P_1,P_2,\ldots,P_{2n}]
=\C[P_1,P_3,\ldots,P_{2n-1}]
=\C[p_1,p_3,\ldots,p_{2n-1}]\subset \Gamma,
\end{equation} 
and the functions $P_1,P_3,\ldots,P_{2n-1}$ are algebraically independent over $\C$.

The following result is well known. 
However, we include it here to clarify certain details, such as the rescaling of variables.
\begin{lem}\label{lem:q and p} Define a {graded} $\C$-algebra homomorphism
\begin{align*}
\C[t_1,t_3,\ldots,t_{2n-1}]&\overset{\iota}\longrightarrow
\C[P_1,P_2,\dotsc,P_{2n}], \\
    \iota(t_{2i-1})&=(2i-1)^{-1}p_{2i-1}\qquad\text{for $1\le i\le n$.}
\end{align*}
Then $\iota$ is an isomorphism, and it satisfies
\begin{equation}
\iota(q_k(t))= 2P_k\qquad\text{for $k\ge1$.}
\label{eq:qto2P}
\end{equation}
\end{lem}
 \begin{proof}
This may be viewed as a finite-variable analogue of \cite[III, \S 8 (8.5)]{Mac}.
\end{proof}



\subsection{$\hq$-functions}\label{sec:qhat}
To set the stage for our later developments in \S\ref{sec:cent}, we introduce a family of functions $\hq_{i}(y|c_1,\ldots,c_m)$.

Let $c_1,c_2,\ldots$ be an arbitrary sequence of parameters.
For $i\ge 0,$ and $m\ge 1$, define
\begin{equation}\label{eq:def_of_hatq}
\hq_{i}(y|c_1,\ldots,c_m)=\sum_{k=0}^\infty
h_{k}(c_1,\ldots,c_m)
{Q}_{i+k}(y).    
\end{equation}
Note that $i$ and $m$ are arbitrary here, but for $i>0$  we have
$$
\hq_{i}(y|a_1,\ldots,a_i)
=2\hat P_{i}(y|a),
$$
by Corollary \ref{cor:hat_P} and $Q_k(y)=2P_k(y).$
Note also that 
\begin{equation}
\hat{q}_0(y|c_1)=\sum_{k=0}^\infty
c_1^k
{Q}_{k}(y)=\Omega(c_1|y).\label{eq:inverseOmega}
\end{equation}
Define \( h_k(x_1,\ldots,x_n; z_1,\ldots,z_m) \in \C[x_1,\ldots,x_n, z_1,\ldots,z_m] \) by 
\begin{equation}
\sum_{k=0}^\infty
h_k(x_1,\ldots,x_n; z_1,\ldots,z_m) u^k
= \frac{\prod_{j=1}^m (1 - z_j u)}{\prod_{i=1}^n (1 - x_i u)}.
\label{eq:h super}
\end{equation}
Recall \( a^{(n)}_i \) from \eqref{E:an} and \eqref{E:an periodic}.

\begin{lem}\label{lem:hatQ in hatP}
Let \( b_i \) for \( 1 \le i \le 2n+1 \) be defined by  
\begin{equation}\label{E:b def}
(b_1,\ldots,b_{2n+1}) = (a_1,\ldots,a_n,0,-a_n,\ldots,-a_1).
\end{equation}
Then, for \( (i,j) \) such that \( 1 \le i \le j \le 2n+1 \), we have 
\begin{equation}
\hq_{j-i}(y | b_i,\ldots,b_j)
= \delta_{j-i,0} + 2 \sum_{k=0}^{2n+i-j}
h_{k}(b_i,\ldots,b_j; a^{(n)}_1,\ldots,a^{(n)}_{j-i+k-1})
\hP_{j-i+k}(y | a^{(n)}).
\label{eq:hat q affine gp}
\end{equation} 
\end{lem}

\begin{proof}
From the general formula \eqref{eq:hat q in hat P}, proved in Section~\ref{ssec:proof of hat q expansion}, we obtain
\begin{equation}
\hq_{j-i}(y | b_i,\ldots,b_j)
= \delta_{j-i,0} + 2 \sum_{k=0}^{\infty}
h_{k}(b_i,\ldots,b_j; a_1,\ldots,a_{j-i+k-1})
\hP_{j-i+k}(y | a).
\label{eq:hat q cite}
\end{equation}
We must show that if \( j - i + k > 2n \), then
\begin{equation}
h_{k}(b_i,\ldots,b_j; a^{(n)}_1,\ldots,a^{(n)}_{j-i+k-1}) = 0,
\label{eq:specialized h}
\end{equation}
so that the right-hand side of \eqref{eq:hat q cite} becomes a finite sum as in \eqref{eq:hat q affine gp} when the parameters \( a \) are specialized to \( a^{(n)} \).
The set \(\{a^{(n)}_1,\ldots,a^{(n)}_{j-i+k-1}\}\) of variables contains all elements \(\pm a_r\) for \(1 \le r \le n\).
Hence, by \eqref{eq:h super}, the parameters \( b_i,\ldots,b_j \) 
(except $b_{n+1}=0$ if it is present) cancel out in the left-hand side of \eqref{eq:specialized h}, which then reduces to the elementary symmetric function \( e_k \) in at most \( k-1 \) variables, which vanishes.
\end{proof}

\begin{ex}
For $n=2,\, i=4,\,j=5$, we have
\begin{equation}
\hat{q}_{1}(b_4,b_5)=
2\hP_1
-2(2a_1+a_2)\hP_2
+4(a_1+a_2)^2\hP_3
-4a_1^2(a_1+a_2)\hP_4,
\end{equation}
where $\hP_i=\hP_i(y|a^{(2)}).$
\end{ex}

\section{Affine nil-Hecke algebra and equivariant homology of the affine Grassmannian} 
\label{sec:Ginzburg-Peterson}

We review the definition and basic results on the affine nil-Hecke algebra and how it can be used to give an explicit algebraic construction of the equivariant homology ring of the affine Grassmannian.

\subsection{Root system and Weyl group}


Let $G$ be a simply connected complex simple linear algebraic group, and fix a maximal torus $T$ of $G$. Denote by $\fg$ and $\fh$ the Lie algebras of
$G$ and $T$ respectively. Let $\Phi\subset \fh^*$ be the root system associated with the pair $(G, T)$.
Choose a set of simple roots $\{\alpha_i\mid i\in I\}\subset \Phi$, where the index set $I$ is referred to as the Dynkin node set. 
Let $\theta\in \Phi$ denote the highest root. 
Let $\Phi^\lor\subset \fh$ denote the dual root system of $\Phi$; for each $\alpha\in \Phi$, let $\alpha^\lor\in \Phi^\lor$ be the corresponding coroot.
The fundamental weights
$\{\varpi_i\mid i\in I\}$ are defined as the dual basis of $\{\alpha_i^\lor\mid i\in I\}$. 
Similarly, the fundamental coweights
$\{\varpi_i^\lor\mid i\in I\}$ are defined as the dual basis of $\{\alpha_i\mid i\in I\}$. 
Denote by $P$ and $Q$ the weight and root lattices, respectively, and by $P^\lor$ and $Q^\lor$ the coweight and coroot lattices, respectively. 
Let $\langle \cdot,\cdot\rangle$ be the natural paring between $\fh$ and $\fh^*$. 
The Weyl group $W$ is the subgroup of $GL(\fh^*)$ generated $s_i\;(i\in I)$ defined by the reflections 
$$
s_i(\la)=\la-\langle \alpha_i^\lor,\lambda\rangle \alpha_i\quad\text{for $\la \in \fh^*$}.
$$ 
Moreover, $W$ is a Coxeter group equipped with the Bruhat order $\le $ and the length function $\ell:W\rightarrow \Z_{\ge 0}.$ The action of $W$ preserves $Q$ and $P$.

\subsection{Extended affine Weyl group}
\label{SS:extended affine Weyl}
The classical Weyl group $W$ also acts naturally on $\fh$ by 
$$
s_i(h)=h-\langle h,\alpha_i\rangle \alpha_i^\lor\quad \text{for $h\in \fh$}.
$$
The coweight and coroot lattices $P^\vee\supset Q^\vee$ are preserved by the action.  
The affine Weyl group is the semidirect product $W_\af \cong Q^\vee \rtimes W$ and the extended affine Weyl group $\Wx$ is by definition $\Wx\cong P^\vee\rtimes W$.

Denote by $X_\af$ the affine Dynkin diagram, $I_\af:=I\cup\{0\}$ its node set and $\Aut(X_\af)\le \mathfrak{S}_{I_\af}$ the group of automorphisms of $X_\af$ realized as a subgroup of permutations $\mathfrak{S}_{I_\af}$ of the set $I_\af$. Let $I_\af^s\subset I_\af$ be the subset of special (also called cominuscule) affine Dynkin nodes, those of the form $\sigma(0)$ where $0\in I_\af$ is the affine node and $\sigma\in \Aut(X_\af)$. There is a bijection $I_\af^s \cong P^\vee/Q^\vee$ such that $i\mapsto \varpi_i^\vee+Q^\vee$ where $\varpi_i^\vee$ is the fundamental coweight, with $\varpi_0^\vee=0$ by convention.
For each $i\in I_\af^s$, the map $x\mapsto x+\varpi_i^\vee+Q^\vee$ defines a permutation of the set $P^\vee/Q^\vee$ which induces a permutation on $I_\af^s$ that extends uniquely to a Dynkin automorphism we denote by $\pi_i\in \Aut(X_\af)$.\footnote{This convention differs from that in \cite{LS:double Kostant}.} 
%
%
Thus the elements $\{\pi_i\mid i\in I_\af^s\}$ form a subgroup $\Sigma$ of $\Aut(X_\af)$ that is isomorphic to $P^\vee/Q^\vee$.

The group $\Sigma$ embeds into $\Wx$ as follows.
For $i\in I_\af^s$,
let $u_i\in W$ be the shortest element such that $u_i^{-1}(\om_i^\lor)$ is anti-dominant. 
Note that 
\begin{equation}
\label{eq:left des}
    \text{$u_i$ has the unique left descent $s_i.$}
\end{equation}
One sees this holds because if $j\ne i$ then $\om_i^\lor$ is fixed by $s_j.$
\begin{ex}
In type $A_3$, we have
$$
\om_1^\lor
\overset{s_1}{\longmapsto}
\om_1^\lor-\alpha_1^\lor
\overset{s_2}{\longmapsto}
\om_1^\lor-\alpha_1^\lor-\alpha_2^\lor\overset{s_3}{\longmapsto}
\om_1^\lor-\alpha_1^\lor-\alpha_2^\lor-\alpha_3^\lor=-\om^\lor_3,
$$
and $u_1=s_1s_2s_3.$
\end{ex}
There is an injective group homomorphism
\begin{align}\label{E:F group}
  \Sigma &\to \Wx\\
  \label{E:length zero}
  \pi_i &\mapsto t_{\varpi_i^\vee} u_i \qquad\text{for $i\in I_\af^s$.}
\end{align}
We will consider $\Sigma$ as a subgroup of $\Wx$ via this map.
The group $\Sigma$ acts on $W_\af$ by conjugation: $\pi s_i \pi^{-1} = s_{\pi(i)}$ for all $\pi\in \Sigma$ and $i\in I_\af$ where the notation $\pi(i)$ means that $\pi$ is viewed as a permutation of $I_\af$. We have $\Wx\cong \Sigma \ltimes W_\af$.

\subsection{Level zero extended affine nil-Hecke algebra}
\label{SS:level zero extended affine nil-Hecke}
The following extended affine Weyl group version of Peterson's level zero nil-Hecke algebra appears in \cite{CP}, \cite{LS:double Kostant}.
The Weyl group $W$ acts naturally on the weight lattice $P$ and therefore on $S=H_T^*(\pt)\cong \Sym(P)$.
Let $\cl:\Wx\to W$ be the group homomorphism defined by $\cl( t_\la v)= v$ for $\la\in P^\vee$ and $v\in W$.
The \emph{level zero action} of $\Wx$ on $P$ is defined by
$w \mu = \cl(w)\mu$. There is an induced action of $\Wx$ on $S$. We shall write $\alpha_0=-\theta$; this is the classical projection of the $0$-th affine simple root, where $\theta$ is the highest classical root.

Let us consider a localization of $S$ defined by $$S^\reg =S[\alpha^{-1}\;(\alpha\in \Phi)],$$ on which $\tilde{W}_\af$ acts naturally. 
Let $S^\reg[\tilde{W}_\af]$ be the \emph{twisted group algebra} of $\tilde{W}_\af$ over $S^\reg$.
Linearly it is the tensor product $S^\reg \otimes_\Z \Z[\Wx]$ with product  
$$(f \otimes v)(g\otimes w) = f v(g) \otimes vw$$ for $f,g\in S^\reg$ and $v,w\in \Wx$.
We regard $S^\reg[W_\af]:=S^\reg \otimes_\Z \Z[W_\af]$ as a subalgebra of $S^\reg [\tilde{W}_\af].$ 
Note that the ring $S^\reg[\tilde{W}_\af]$ acts naturally on $S^\reg$.

For $i\in I_\af$, define $A_i\in S^\reg[\Wx]$ by
\begin{equation}
A_i = \alpha_i^{-1}(1-s_i).
\end{equation}
Then we have $A_i^2=0$, and the elements $A_i$ satisfy the braid relations. For $w\in W_\af$, one can define $A_w=A_{i_1}\cdots A_{i_l}$ where $w=s_{i_1}\cdots s_{i_l}$ is a reduced expression. 
We further define
\begin{align}\label{E:A extended}
A_{\pi v} = \pi A_v \qquad\text{for $\pi\in\Sigma$ and $v\in W_\af$.}
\end{align}

Define the left $S$-submodules 
$$\Ax=\sum_{w\in \Wx}S A_w,\quad 
\A_\af=\sum_{w\in W_\af}S A_w$$  of $S^\reg[\Wx].$  
The action of $S^\reg[\tilde{W}_\af]$ on $S^\reg$ restricts to an action of $\Ax$ on $S$.
Similarly $\nilA$ acts on $S$.
One may show that $\A_\af$ and $\Ax$ are rings. 
We can regard $\A_\af$ as a subring of $\Ax$. 
We call $\A_\af$ (resp. $\Ax$) the \emph{level $0$ (resp. extended) affine nil-Hecke algebra}. 
Note that $S$ is not central in $\nilA$ nor in $\Ax$.

\subsection{Peterson subalgebra}
\label{ssec:affnil}

Let $\Pet$ be the \emph{Peterson subalgebra}, the centralizer of $S$ in $\nilA$. 
Similarly let $\Petx$ be the centralizer of $S$ in $\Ax$.

Note that the translation element $t_\la\;(\la\in P^\lor)$ is an element of $\Petx$, because we consider the level zero action of $W_\af$ on $S$.
We have \cite[Lemma 4.1]{Lam} 
\begin{align}\label{E:rat Peterson}
\Petx^\reg&:=S^\reg\otimes_S\Petx=\bigoplus_{\la\in P^\vee} S^\reg t_\la,\\
\Pet^\reg&:=
S^\reg\otimes_S\Pet = \bigoplus_{\la\in Q^\vee} S^\reg t_\la.
\end{align}
Following Peterson \cite{Pet}, define the left action $*$ of $\tilde{W}_\af$ on $\Petx^\reg$ by 
\begin{align*}
  w * (f t_\mu) &
  = w f t_\mu w^{-1} 
  = w(f) t_{w(\mu)} &\qquad&\text{for $w\in W$, $f\in S^\reg$, $\mu\in P^\vee$} \\
  t_\la * (f t_\mu) &
  = f t_{\la+\mu}&&\text{for $\la,\mu\in P^\vee$.}
\end{align*}
This action is well-defined and restricted to an action of $W_\af$ on $\Pet^\reg$. Combined with the usual action of $S^\reg$ on $\Petx^\reg$, we obtain an action of $S^\reg[\tilde{W}_\af]$ on $\Petx^\reg$. One may show that this restricts to an action of $\Ax$ on $\Petx$, and further to an action of $\A_\af$ on $\Pet.$ 

\begin{prop}\label{prop:Pet *}
Let \( b, b' \in \Petx \). Then we have  
\begin{equation}
b \, b' = b * b'.
\end{equation}    
\end{prop}
\begin{proof}
It suffices to show the result for  $b,b'\in \Petx^\reg$.
For \( f, f' \in S^\reg \) and \( \mu, \mu' \in P^\lor \), let  
\[
b = f t_\mu, \quad b' = f' t_{\mu'}.
\]
Then, we compute  
$$b * b' = (f t_\mu) * (f' t_{\mu'}) 
       = f f' t_{\mu + \mu'} 
       = b b'.$$
Thus, the proposition follows.  
\end{proof}

Although $\Wx$ is not a Coxeter group, the length function on $\Wx$ is defined by $\ell(\sigma w)=\ell(w)$ for $\sigma \in \Sigma,\,w\in W_\af$.
Let $\Wxgr$ (resp. $\Wafgr$) be the affine Grassmannian elements of $\Wx$ (resp. $W_\af$), the elements $w$ that are of minimum length in their coset $w W$.
We have $\Wxgr = \Sigma \Wafgr$. 
The following is a straightforward extension to $\Petx$ of Peterson's characterization of the Schubert basis in $\Pet$.

\begin{thm} [\cite{Pet}, cf. \cite{ISY}]\label{thm:Pet} The Peterson algebra $\Petx$ is a commutative $S$-algebra. For each $w\in  \Wxgr$, there exists a unique element $j_w\in\Petx$ of the form  
\begin{equation}
    j_w = A_w + \sum_{v \in \Wx\setminus \Wxgr}j_w^v A_v,
\end{equation}
with coefficients  $j_w^v\in S$. 
Moreover, the subset $\{j_w\}_{w\in\Wxgr}$ forms an $S$-basis of $\Petx$. 
The analogous statement holds when one replaces \(\Petx,\Wxgr,\Wx\) by \(\Pet,W_\af^0,W_\af\), respectively.
\end{thm}

\begin{prop}[\cite{Pet}, cf. \cite{ISY}]\label{P:j basis by left action}
For $w\in \tilde{W}_\af$ and $v\in \Wxgr$ we have
\begin{align}\label{E:A on j}
  A_w * j_v = \begin{cases}
      j_{wv} & \text{if $wv\in \Wxgr$ and $\ell(w)+\ell(v)=\ell(wv)$} \\
      0 & \text{otherwise.}
  \end{cases}
\end{align}
In particular, for $w\in \Wxgr$ we have 
\begin{align}\label{E:j basis by left action}
    j_w = A_w * 1.
\end{align}
\end{prop}



\begin{ex}\label{ex:js0}
We have
$$
j_{s_0}=A_0*1={\alpha_0}^{-1}
(1-s_0*1)={\alpha_0}^{-1}
(1-t_{\theta^\lor}s_\theta*1 )
={\alpha_0}^{-1}
(1-t_{\theta^\lor}).
$$
\end{ex}
Recall that we identify $H_{T_\ad}^*(\mathrm{pt})$ with $H_T^*(\mathrm{pt})=S$, since we are working over $\C$, and the Lie algebras of $T$ and $T_\ad$ are naturally identified.
\begin{thm}[\cite{Pet}, \cite{Lam}, \cite{LS:Acta}]\label{thm:Pet and homology}
There are isomorphisms of $S$-algebras
\begin{align}
H_*^T(\Gr_G) &\cong \Pet, \quad
H_*^{T_\ad}(\Gr_{G_\ad}) \cong \Petx.
\end{align}
Moreover, for $w\in W_\af^0$ (resp. $w\in \Wxgr$) the $w$-th Schubert class in $H^T_*(\Gr_G)$ (resp. $H^{T_\ad}_*(\Gr_{G_\ad})$) is identified with $j_w\in \Pet$ (resp. $j_w \in \Petx$).
\end{thm}



\subsection{Map $j$}
Let $j: \Ax\rightarrow \Petx$ be the left $S$-module map defined by
\begin{equation}
j(A_w)=\begin{cases}
    j_w & \text{if $w\in \Wafgr$}\\
    0 & \text{otherwise.}
\end{cases}
\end{equation}
For $a\in \Petx$, we have
\begin{equation}
j(a)=a.
\end{equation}

\begin{prop}
For $a\in \Ax$ and $i\in I$, we have
\begin{align}
j(aA_i)&=0,\label{eq:j(aAi)}\\
j(as_i)&=j(a),\label{eq:jasi}\\
j(s_i*a)&=j(s_i a),\label{eq:jsi*}\\
j(A_i*a)&=j(A_i a).
\label{eq:jAi*a}
\end{align}
\end{prop}
\begin{proof}
For \eqref{eq:j(aAi)}, it suffices to show for $a=A_w.$ We have $$A_wA_i=\begin{cases}
A_{ws_i} & ws_i>w\\
0 & \text{otherwise}.
\end{cases}$$
If $w_si>w$, then $ws_i\notin \Wxgr.$ Equation \eqref{eq:jasi} holds since $$j(a)-j(as_i)=j(a(1-s_i))=j(a\alpha_i A_i)=0$$ by \eqref{eq:j(aAi)}.
We have $
j(s_i* a)=j(s_i as_i)=j(s_i a)
$ by \eqref{eq:j(aAi)}, so we obtain \eqref{eq:jsi*}.
For \eqref{eq:jAi*a}, we compute
\begin{align*}
\alpha_i j(A_i*a)&=
j(\alpha_i A_i*a)=
j((1-s_i)*a)\\
&=j((1-s_i)a)\quad \text{(by \eqref{eq:jsi*})}\\
&=
j(\alpha_i A_i a)=\alpha_ij( A_i a).
\end{align*}
By dividing by $\alpha_i$, we obtain \eqref{eq:jAi*a}.
\end{proof}
\subsection{Anti-dominant translation elements}
\begin{lem}[\cite{ISY}]\label{lem:j of anti-dom W invariant}
For any anti-dominant
$\la\in P^\lor,$ we have
$t_\la\in \tilde{W}_\af^0$, and 
$j_{t_\la}$ is $W$-invariant under the $*$ action. 
\end{lem}
\begin{proof}
The proof is the same as \cite[Lemma 2.20]{ISY}.
\end{proof}

\subsection{A criterion for membership in $\Wafgr\subset W_\af$}
Let $\Phi_\af=(\Phi\oplus \Z\delta)\setminus\{0\}$ denote the affine root system of the untwisted affine Lie algebra $\hat{\fg}$ associated with $\fg$. The set $\Phi_\af^{+}$ of affine positive roots consisting of $\gamma=\alpha+k\delta\in \Phi_\af$ such that either $k>0$ or both $\alpha\in \Phi^{+}$ and $k=0$. 
The set $\Phi_\af^{-}$ of affine negative roots is given by $-\Phi_\af^{+}.$

The \emph{level-zero action} of $\tilde {W}_\af$ on $P\oplus \Z\delta$ is given by 
\begin{equation}
wt_\la (\gamma+k\delta)=
w(\gamma)+(k-\langle \lambda,\gamma\rangle)\delta,
\end{equation}
for $w\in W,\,\la\in P^\lor,\,\gamma\in P,$ and $k\in \Z.$


\begin{lem} \label{L:inversion} For $x\in W_\af$ and $i\in I_\af$,
$x s_i > x$ if and only if $x(\alpha_i)\in \Phi_\af^+$.
\end{lem}
\begin{proof} Define the inversion set of $x$ by $S(x):=\{\alpha\in \Phi_\af^+\mid 
x(\alpha)\in \Phi_\af^-\}$.
It is known \cite[(2.2.4)]{Mac:aff} that for $x,y\in W_\af$, 
$\ell(xy)=\ell(x)+\ell(y)\Longleftrightarrow S(y)\subset S(xy)$.
Letting $y=s_i$, the following are equivalent: $xs_i>x$; $\{\alpha_i\} = S(s_i) \subset S(xs_i)$; $xs_i(\alpha_i)\in -\Phi_\af^+$; $x(\alpha_i)\in \Phi_\af^+$.
\end{proof}

\begin{lem}\label{lem:0Grass}
Let $x\in W_\af$. Then
\begin{equation}
    x\in W_\af^0 
    \Longleftrightarrow 
x(\alpha_i)\in \Phi_\af^+\qquad\text{for all $i\in I$.}
\end{equation}  
\end{lem}
\begin{proof} This follows by applying Lemma \ref{L:inversion} for $i\in I$.
\end{proof}






{
\subsection{Star action of the fundamental group $\Sigma$}

For $\sigma\in \Sigma$, we write $\sigma=t_\sigma u_\sigma$ with $t_\sigma$ a translation element corresponding to an element of $P^\lor$ and $u_\sigma\in W.$ In this notation, note that we have 
\begin{equation}
j_\sigma=(t_\sigma u_\sigma)*1=t_\sigma.
\end{equation}
Moreover, for $a\in \Pet,$ we have 
\begin{equation}
\sigma* a=t_\sigma\cdot  u_\sigma(a)
=j_\sigma\cdot  u_\sigma(a),
\end{equation}
where $u_\sigma(a)=u_\sigma*a=u_\sigma au_\sigma^{-1}.$
In particular, for $w\in \Wxgr$, we have 
\begin{equation}\label{eq:pi acts on j}
    \sigma*j_w=j_{\sigma w}=j_\sigma u_\sigma(j_{w}).
\end{equation}
}
\subsection{Factorization property }\label{ssec:factor}
Let $\la\in P^\lor$ be an anti-dominant coweight.
The corresponding translation element $t_\la$ is an affine Grassmannian element of $\Wx$. 
In particular, for $i\in I$, $t_{-\varpi_i^\lor}\in \Wxgr$. 
Denote by $i\mapsto i^s$ the map $I\to I^s$ defined by $\om_i^\vee + Q^\vee = \om_{i^s}^\vee+Q^\vee$.
In particular for $i\in I^s$, $i^s=i$.
There is a unique element $\kappa_i$ of $\Wafgr$ such that 
\begin{equation}
\label{eq:translation of negative cofundamental}
t_{-\varpi_i^\lor}
=
\pi_{i^s}^{-1} \kappa_i
\end{equation}
where, for special $j$,
$\pi_j\in \Sigma$ is defined in \S~\ref{SS:extended affine Weyl}.
For $i\in I^s$, $$\kappa_i=\pi_iu_i\pi_i^{-1},$$ with $u_i\in W$ as defined by \eqref{E:length zero}.



\begin{prop}\label{prop:kappa}
For $i\in I$ we have
$$j_{\kappa_i}=j_{\pi_{i^s}}j_{t_{-\varpi_i^\lor}}.$$
\end{prop}
\begin{proof} We have $\kappa_i = \pi_{i^s} t_{-\om_i^\vee}$.
By \eqref{eq:translation of negative cofundamental}, and \eqref{eq:pi acts on j}, we obtain
$$
j_{\kappa_i}
=j_{\pi_{i^s} t_{-\varpi_i^\lor}}
=j_{\pi_{i^s}} u_{i^s}(j_{t_{-\varpi_i^\lor}}).
$$
Since $j_{t_{-\varpi_i^\lor}}$ is $W$-invariant (Lemma \ref{lem:j of anti-dom W invariant}), the desired result follows.
\end{proof}

\begin{lem} \label{L:trans grass}
Let $i\in I$ and $v\in W_\af$ be such that $vt_{-\om_i^\vee}\in\Wxgr$
and $\ell(vt_{-\om_i^\vee})=\ell(v)+\ell(t_{-\om_i^\vee})$. Then $v\in \Wafgr$.
\end{lem}
\begin{proof} Let $v$ and $i$ be as in the hypotheses and suppose $v\not\in\Wafgr$. Then there is a $j\in I$ such that $v(\alpha_j)\in-\Phi_\af^+$.
We have $vt_{-\om_i^\vee}(\alpha_j) = v (\alpha_j + \delta_{ij} \delta) = v(\alpha_j) + \delta_{ij}\delta\in\Phi_\af^+$ due to the assumption that $v t_{-\om_i^\vee}\in \Wxgr$. This implies $i=j$ and $vs_i<v$.
But $s_i t_{-\om_i^\vee}<t_{-\om_i^\vee}$ since $t_{\om_i^\vee}(\alpha_i) = \alpha_i - \delta$.
This contradicts the length-additivity.
\end{proof}

\begin{lem}\label{lem:v bar}
Let $i\in I$ and $v\in W_\af$ such that $v\kappa_i\in W_\af^0$ and $\ell(v\kappa_i)=\ell(v)+\ell(\kappa_i)$.
Then $\bar{v}:=\pi_{i^s}^{-1}v\pi_{i^s}\in W_\af^0.$
\end{lem}
\begin{proof} We have $v \kappa_i = v \pi_{i^s} t_{-\om_i^\vee} = \pi_{i^s} \bar{v} t_{-\om_i^\vee}$.
Since the group $\Sigma$ preserves $\Phi_\af^+$, $\bar{v} t_{-\om_i^\vee}\in \Wxgr$.
Note that conjugation by $\pi_{i^s}^{-1}\in\Sigma$ is a length-preserving group automorphism of $\Wx$ that stabilizes $W_\af$.
Using also the fact that left or right multiplication by a length-zero element does not change the length, we deduce that
$\ell(\bar{v})+\ell(t_{-\om_i^\vee})= \ell(\bar{v} t_{-\om_i^\vee})$.
By Lemma \ref{L:trans grass} we deduce that $\bar{v}\in \Wafgr$.
\end{proof}

\begin{thm}\label{thm:factor kappa i}
Let $i\in I$, and  $v\in W_\af$ such that $v\kappa_i\in \Wafgr$ and $\ell(v\kappa_i)=\ell(v)+\ell(\kappa_i)$.
Then $\bar{v}:=\pi_{i^s}^{-1}v\pi_{i^s}\in\Wafgr$ and
$$
j_{v\kappa_i}=u_{i^s} (j_{\bar{v}}) j_{\kappa_i}.
$$
\end{thm}
\begin{proof}
We have $\bar{v}\in W_\af^0$ by Lemma \ref{lem:v bar}. 
We compute
\begin{align*}
j_{v\kappa_i}
 &= A_v*j_{\kappa_i}&&
\\
&=A_v*\left(j_{\pi_{i^s}}j_{t_{-\varpi_i^\lor}}\right)&&\text{(by Proposition \ref{prop:kappa})}
\\
&=
j_{t_{-\varpi_i^\lor}}A_v*j_{\pi_{i^s}}
&& \text{(by Lemma \ref{lem:j of anti-dom W invariant})}\\
&=
j_{t_{-\varpi_i^\lor}}
j_{v\pi_{i^s}}=
j_{t_{-\varpi_i^\lor}}
j_{\pi_{i^s}\bar{v}}\\
&=
\left(j_{\pi_{i^s}}^{-1}j_{\kappa_i}\right)
\left(
j_{\pi_{i^s}}
u_{i^s}(j_{\bar{v}})\right)&& \text{(by Proposition \ref{prop:kappa} and \eqref{eq:pi acts on j})}\\
&=
j_{\kappa_i}
u_{i^s}(j_{\bar{v}}).
\end{align*}
\end{proof}

\begin{ex} In type $A_2^{(1)}$, all the nodes are special. Let $v=s_2s_1$ and $i=1$. Then $\bar{v}=s_1s_0.$
We have $t_{-\om_1^\vee} = \pi_1^{-1} \kappa_1$ with $\kappa_1=s_2s_0.$
Then $j_{v\kappa_1}=j_{s_2s_1s_2s_0} = j_{s_2s_0}\pi_1(j_{s_1s_0})$.
\end{ex}

\begin{ex} In type $C_2^{(1)}$ let $v=s_2$ and $i=2$, which is special. Then $t_{-\varpi_2} = \pi_2^{-1} s_0 s_1 s_0$. We have
$\bar{v}= s_0$ and
\begin{align*}
  j_{(s_2)(s_0 s_1 s_0)} = j_{s_0s_1s_0}\pi_2(j_{s_0}).
\end{align*}
For $v=s_0$ and $i=1$, which is not special, we have $\kappa_1=s_1s_2s_1s_0,$
and 
\begin{align*}
j_{s_0s_1s_2s_1s_0} = j_{s_1s_2s_1s_0}j_{s_0} .
\end{align*}
\end{ex}
\begin{ex}
In type $B_3^{(1)}$, for $i=3$ we have $i^s=1$ and $t_{-\varpi_3}=\pi_1 \kappa_3$, where $$\kappa_3=\rho_3\rho_3'\rho_3, \quad
\rho_3=s_3s_2s_0,\quad\rho_3'=\pi_1(\rho_3)=s_3s_2s_1.$$
Take $v=s_1$. 
Then we have
$$
j_{v\kappa_3}=\pi_1(j_{\overline{v}})
j_{\kappa_3}.
$$
\end{ex}








\section{Notation for Type $C_n^{(1)}$ setting}
\label{sec:C}

\subsection{Lie algebra $\fsp_{2n}(\C)$ of type $C_n$}

Let $P$ be the weight lattice of $\fg=\fsp_{2n}(\C).$
We choose a $\Z$-basis of $P$ so that the fundamental weight are 
\begin{equation}
\varpi_i=a_1+\cdots+a_i\quad\text{for $1\le i\le n$}.
\end{equation}
We have $P=\bigoplus_{i=1}^n \Z a_i$ and the simple roots are given by 
\begin{equation}
\alpha_i=a_i-a_{i+1}\quad
\text{for $1\le i\le n-1$},\quad
\alpha_n=2a_n.
\end{equation}
Set $I=\{1,2,\ldots,n\}$ as the Dynkin node set.

Let $\{\eps_i\mid 1\le i\le n \}\subset\Hom_\Z(P,\Z)$ be the dual basis of $\{a_i\mid 1\le i\le n\}$.
The simple coroots are \begin{equation}
\alpha_i^\vee=\eps_i-\eps_{i+1} \quad\text{for $1\le i\le n-1$,}\quad \alpha_n^\vee=\eps_n.
\end{equation}
We have $Q^\vee = \bigoplus_{i=1}^n \Z\eps_i.$
The fundamental coweights are \begin{equation}   \om_i^\vee=\eps_1+\dotsm +\eps_i\quad \text{ for $1\le i\le n-1$},\quad
\om_n^\vee = \frac{1}{2} \sum_{k=1}^n \eps_i.
\end{equation}
The coroot lattice $Q^\vee$ is an index $2$ sublattice of $P^\lor=\bigoplus_{i=1}^n \Z\varpi_i^\lor$.
The highest root is $\theta=2a_1$ and the corresponding coroot is $\theta^\lor=\eps_1.$
For $1\le i\le n-1$, the simple reflection $s_i$ acts on $P$ by exchanging $a_i$ and $a_{i+1}$ and fixing other $a_j$. $s_n$ acts on $P$ by negating $a_n$ and fixing other $a_j$.

We have $I_\af^s=\{0,n\}$ and $\Sigma=\{\pi_0,\pi_n\}\cong \Z/2\Z.$
The element $u_n$ is given by 
\begin{align}\label{E:u n}
u_n =  s_n(s_{n-1}s_n) \dotsm(s_{2}\dotsm s_{n-1}s_n)(s_1\dotsm s_{n-1}s_n).
\end{align}

Since the action of translation elements on $S$ is trivial, the actions of $\pi_n$ and $u_n$ coincide on $S$; in particular, for $1\le i \le n,$ 
\begin{align}\label{E:un on P}
    \pi_n(a_i) = u_n(a_i) = -a_{n+1-i}.
\end{align}

\subsection{Affine Weyl group of type $C_n^{(1)}$}\label{ssec:Cnil-Hecke}
The Dynkin diagram  of type $C_n^{(1)}$is
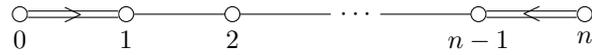
\vspace{0.3cm}
\begin{center}
\begin{tikzpicture}[node distance=1.5cm and 1.2cm, auto]
    \node[draw, circle, inner sep=2pt, label=below:{$0$}] (0) {};
    \node[draw, circle, inner sep=2pt, label=below:{$1$}, right=of 0] (1) {};
    \node[draw, circle, inner sep=2pt, label=below:{$2$}, right=of 1] 
    (2) {};
    \node[right=of 2] (dots) {\(\cdots\)};
    \node[draw, circle, inner sep=2pt, label=below:{$n-1$}, right=of dots] (n-1) {};
    \node[draw, circle, inner sep=2pt, label=below:{$n$}, right=of n-1] (n) {};
    \draw[double distance=2pt] (0) -- (1) node[midway, yshift=-1.5ex] {\(>\)};
    \draw (1) -- (2);
    
    \draw (2) -- (dots);
    \draw (dots) -- (n-1);
    
    \draw[double distance=2pt] (n) -- (n-1) node[midway, yshift=1.5ex] {\(<\)};
\end{tikzpicture}
\end{center}
Let $W_\af=\langle s_0,s_1,\ldots,s_n\rangle $ be the corresponding affine Weyl group.
The set of affine Dynkin nodes is $I_\af=\{0,1,\ldots,n\}$.
As an affine Dynkin diagram  automorphism, $\pi_n$ reverses the Dynkin diagram, that is, $\pi_n(i)=n-i$ for $i\in I_\af$.


\subsection{Translation elements 
in terms of $j$-basis elements}

The localization value $\xi^v(w)\in S$ of the $v$-th affine flags $T$-equivariant Schubert class at the $w$-th fixed point, is defined by the following expression in $\Ax$ \cite{KK}:
\begin{align}\label{E:xi def}
  w = \sum_{v\le w} \xi^v(w) A_v.
\end{align}
One may compute $\xi^v(w)$ as follows. We only need it for $v,w\in W_\af$.
Let $w=s_{i_1} s_{i_2} \dotsm s_{i_\ell}$ be a reduced decomposition of $w$.
Let $b_1b_2\dotsm b_\ell$ be a binary vector such that the ordered sequence $\{ i_k \mid b_k=1\}$ is a reduced word for $v$. 
Summing over all such binary vectors, we obtain a left-to-right ordered composition of operators acting on $1$, where each reflection acts to its right and the negative simple roots act by left multiplication: \cite{AJS},\cite{Bi}
\begin{align}\label{E:AJSBilley}
  \xi^v(w) = \sum_{\substack{(b_1,\dotsc,b_\ell)\in \{0,1\}^\ell \\ A_{i_1}^{b_1}A_{i_2}^{b_2}\dotsm A_{i_\ell}^{b_\ell} = A_v}} \left(\prod_{i=1}^\ell
  \begin{cases}
    (-\alpha_{i_k}) s_{i_k}& \text{if $i_k=1$} \\
    s_{i_k}&\text{otherwise}
  \end{cases}\right) \cdot 1\qquad\text{for $v,w\in W_\af$.}
\end{align}



\begin{lem}\label{L:loc value}
\begin{align}
  \xi^{\rho_k}(t_{\eps_i}) = \xi^{\rho_k}(\rho_i) = (\!(a_i|a)\!)^k\qquad\text{for $1\le k\le i\le n$.}
\end{align}
\end{lem}
\begin{proof}
We have 
$$t_{\eps_i} = s_{i-1}\dotsm s_1 s_0 s_1 \dotsm s_{n-1} s_n s_{n-1}\dotsm s_{i+1} s_i
\qquad\text{
for $1\le i\le n$.}$$ 
The elements of $\Wxgr$ below $t_{\eps_i}$ are $1$ and $\rho_k$ for $1\le k\le i$.
Computing $\xi^{\rho_k}(t_{\eps_i})$ the only vector $b$ that works, has $0$'s for reflections down to $s_k$, then $1$'s for $s_{k-1}$ through $s_0$, and then $0$'s for the rest.
Therefore,
\begin{align*}
\xi^{\rho_k}(t_{\eps_i}) &= \xi^{\rho_k}(\rho_i) = \left( s_{i-1}\dotsm s_k (-\alpha_{k-1} s_{k-1})
\dotsm(-\alpha_1 s_1) (-\alpha_0 s_0) \right) \cdot 1\\
&=(a_i-a_{k-1}) (a_i-a_{k-2})\dotsm (a_i-a_1) 2 a_i. 
\end{align*}
\end{proof}

For $w\in \Wxgr$, denote by $t^w\in\Wx$ the unique translation element such that $t^w W = wW$.
\begin{prop}\label{P:translation to j} For $w\in\Wxgr$
\begin{align}\label{E:translation to j}
t^w= \sum_{\substack{v\in\Wxgr \\ v\le w}} \xi^v(t^w) \,j_v.
\end{align}
\end{prop}
\begin{rem} See \cite[(4.12)]{LLS:coproduct} for the $K$-theoretic analogue.
\end{rem}

Combining Proposition \ref{P:translation to j} and Lemma \ref{L:loc value} we have
    \begin{equation}\label{E:t eps} 
    t_{\eps_i}=1+\sum_{k=1}^i (\!(a_i|a)\!)^k j_{\rho_k}\qquad\text{for $1\le i\le n$.}
    \end{equation}


\section{
Realization of $H_*^T(\Gr_{\Sp_{2n}(\C)})$ by Symmetric Functions
}\label{sec:Sym}

{
Recall $S$ and $\hLa_S$ from \S \ref{ssec:affine dual double Schur P}.
The goal of this section is to realize the homology ring $H_*^T(\Gr_{C_n})$ as an $S$-subalgebra $\hGn$ of $\hLa_S$ in such a way that the Schubert basis is identified with an explicitly constructed family of symmetric functions.
\subsection{Action of $\nilA$ on $\hat{\La}_S$}

The Weyl group $W$ acts on $S=\C[\mathfrak{h}].$
We extend this action of $W$ to an action on $\hLa_S$ by letting each $w\in W$ act only on the coefficients.

\begin{prop}\label{P:Waf on symmetric functions}
The action of $W$ on $\hLa_S$ extends to an action of the affine Weyl group $W_\af$ on $\hLa_S$ by defining the action of translation elements $t_{\eps_i}\;(1\le i\le n)$ via 
$$
t_{\eps_i}f(y|a)=
\Omega(a_i|y)f(y|a)
\quad 
\text{for $f(y|a)\in \hLa_S$}.
$$
In particular, the simple reflection $s_0$ acts as   \begin{equation}
    (s_0 f)(y|a)
    =\Omega(a_1|y)f(y|-a_1,a_2,\ldots,a_n).\label{eq:WafGa}
\end{equation}  Together with the left $S$-module structure, this action of $W_\af$ endows   $\hLa_S$ with the structure of a left $\A_\af$-module.
\end{prop}
\begin{proof} 
It is straightforward to verify that, as linear operators on $\hLa_S$, one has  $wt_{\eps_i}w^{-1}=t_{w(\eps_i)}$ for $w\in W$ and $1\le i\le n$. 
Thus, we obtain an action of $W_\af$ on $\hLa_S.$
Equation \eqref{eq:WafGa} follows from the relation  $s_0=t_{\theta^\lor}s_\theta$. 
Moreover, the $W_\af$ action and the left $S$-module structure satisfy the defining relation of the twisted group ring of $W_\af$. 
Finally, just as in the infinite rank case, we can define an action of $A_i$ ($i\in I_\af$) on $\hLa_S.$ 
Hence, the Proposition holds.
\end{proof}
\begin{rem}\label{rem:Aaf commute with y}
Note that the action of $\A_\af$ on $\hLa_S$ commutes with the multiplication by elements of $\hLa.$
\end{rem}

\subsection{Dual affine Grassmannian double Schur $P$-function}
\begin{defn}
For $w\in W_\af^0$, the \emph{dual affine Grassmannian double Schur $P$-function} associated with $w$ is by definition  
\begin{equation}
\hP^{(n)}_w(y|a)=A_w (1)\in \hLa_S.
\end{equation}
\end{defn}

Let $\hGn$ denote the $S$-span of $\hP^{(n)}_w(y|a)\;(w\in W_\af^0).$
 
\begin{thm}\label{thm:gamma}
$\hGn$ is an $S$-subalgebra of $\hLa_S$, and there is an  
isomorphism $$\gamma: H_*^T(\Gr_{\Sp_{2n}})\rightarrow \hGn$$ of $S$-algebras and left $\A_\af$-modules that sends the Schubert class $\sigma_w$ to the dual affine Grassmannian double Schur $P$-function $\hP^{(n)}_w(y|a)$ for $w\in \Wafgr$.
\end{thm}

\subsection{Proof of Theorem \ref{thm:gamma}}

Via the isomorphism Theorem \ref{thm:Pet and homology} due to Peterson, we identify the equivariant homology ring $H_*^T(\Gr_{\Sp_{2n}})$ with $\Pet$, and hence we do not distinguish between them.

We use the following fact. 
\begin{lem}\label{lem:Jac} 
The elements $\Omega(a_i|y)\;(1\le i\le n)$ are algebraically independent over the fraction field of $S.$
\end{lem}
\begin{proof} 
Since $\frac{\partial\Omega(a_i|y)}{\partial y_j}=\Omega(a_i|y)\cdot \frac{2a_i}{1-a_i^2y_j^2}$, we have
\[
\begin{aligned}
\det\left(\frac{\partial\Omega(a_i|y)}{\partial y_j}\right)_{1\le i,j\le n}
&=\prod_{i=1}^n 2a_i\Omega(a_i|y)
\times
\det\left(
{(1-a_i^2y_j^2)}^{-1}
\right)_{1\le i,j\le n}
.
\end{aligned}
\]
This is nonzero for generic $y$.
Hence, the Lemma follows by the Jacobian criterion.
\end{proof}


Let $\gamma: \Pet\to \hGn$ be the $S$-module homomorphism defined by 
\begin{align}
\label{E:j map}
  \gamma(b) = b\cdot 1\qquad\text{for $b\in \Pet$.}    
\end{align}

\begin{lem}\label{lem:gamma A}
The map $\gamma$ is a left $\A_\af$-module map.  
\end{lem}
\begin{proof} It is known that $\A_\af$ acts on both $\Pet$ and $\hGn$. Working over $\mathrm{Frac}(S)$, it suffices to verify the equivariance only for $w\in W_\af$ and $s\in S$.
We have
\begin{align*}
  \gamma(w*b) &= \gamma(wbw^{-1}) = wbw^{-1} \cdot 1 = w b\cdot 1 = w \gamma(b) \\
  \gamma(s*b )&= \gamma(s b) = s b \cdot 1 = s \gamma (b).
\end{align*}
Hence, the action of $\A_\af$ is equivariant with respect to $\gamma$.
\end{proof}

From this lemma, we immediately have the following.
\begin{lem} 
\label{L:j to P}
Let $w\in \Wafgr$, and $1\le i\le n.$ We have
\begin{align}
\gamma(j_w)&=
P^{(n)}_w(y|a)\label{eq:gamma j},\\
\gamma(t_{\eps_i}) &=\Omega(a_i|y).
\label{E:jt}
\end{align}
\end{lem}
\begin{proof}

By Lemma \ref{lem:gamma A}, we have
\begin{align*}
\gamma(j_w)&=\gamma(A_w*1)=A_w \gamma(1)
=A_w(1)=P^{(n)}_w(y|a),\\
\gamma(t_{\eps_i})
&=\gamma(t_{\eps_i}*1)
=t_{\eps_i}\gamma(1)
=\Omega(a_i|y).
\end{align*}
\end{proof}

%


\begin{prop} 
The map $\gamma$ is a homomorphism of $S$-algebras.
\end{prop}
\begin{proof}
Since \( \gamma \) is \(S\)-linear by definition, it remains to show that it preserves the product. 
This follows easily by working over $\mathrm{Frac}(S)$: $\Pet$ is spanned by translation elements, and the action of translation elements is multiplication by corresponding Cauchy kernels, and the map on the translation elements is a group homomorphism.
\end{proof}

By Lemma \ref{L:j to P}, the map $\gamma$ is surjective, as its image contains the set of $S$-module generators $\hP_w^{(n)}(y | a)$ for $\hGn$. 
In order to complete the proof of Theorem \ref{thm:gamma}, it remains only to prove the following.


\begin{prop}\label{P:j injective}
The map $\gamma$ is injective.
\end{prop}
\begin{proof}
Localizing, define
\[
\gamma^\reg := 1_{S^\reg}\otimes_S \gamma \colon \Pet^\reg \longrightarrow \hGn^\reg.
\]
Note that  \(\Pet^\reg \cong S^\reg[Q^\lor]\).
By \eqref{E:jt} the map \(\gamma^\reg\) sends Laurent monomials in the \(t_{\eps_i}\) to Laurent monomials in the \(\Omega(a_i|y)\). By Lemma~\ref{lem:Jac}, these Laurent monomials are linearly independent over \(S^\reg\), so \(\gamma^\reg\) is injective.

Finally, as \(\hGn\) is an integral domain, it naturally embeds into \(\hGn^\reg\). Hence, the restriction
\[
\gamma = \gamma^\reg|_{\Pet}\colon \Pet \to \hGn
\]
is injective as required.
\end{proof}
}




\subsection{Realization of the Extended Peterson algebra}
\label{ssec:extended_sym}

The adjoint group of type \( C_n \) is given by
\[
G_\mathrm{ad} = \Sp_{2n}(\mathbb{C})/\{\pm \mathrm{id}\}.
\]
Let \( T_\mathrm{ad} \) denote the image of a maximal torus \( T \subset \Sp_{2n}(\mathbb{C}) \) under the natural projection, so that \( T_\mathrm{ad} = T/\{\pm \mathrm{id}\} \).  
\begin{rem}
The equivariant cohomology ring of a point with respect to \( T_\mathrm{ad} \) is then identified with the symmetric algebra
\[
H_{T_\mathrm{ad}}^*(\mathrm{pt};\mathbb{Z}) \cong \mathrm{Sym}_\mathbb{Z} Q,
\]
where \( Q \) is the root lattice of type \( C_n \), viewed as a sublattice of the weight lattice \( P \).  
Over the complex numbers, the cohomology rings \( H_T^*(\mathrm{pt}) \) and \( H_{T_{\mathrm{ad}}}^*(\mathrm{pt}) \) can be naturally identified, since they are both polynomial rings on the same Lie algebra; thus, the distinction is often unnecessary in practice.  
\end{rem}

We discuss the extended Peterson algebra $\Petx$ and its realization in terms of symmetric functions.
For this purpose, we will extend the $\A_\af$-action on $\hLa_S$ to  $\tilde{\A}_\af$.
 
The crucial step is to define the action of the translation element $t_{\varpi_n^\lor}$ in $\Wx$. 
This action is given by multiplication by 
 $$\eta(y|a):=
\sqrt{\prod_{i=1}^n\Omega(a_i|y)},
 $$
 which is a well-defined formal power series in $\hLa_S$ with a constant term of $1.$ Here, the square root is taken in the ring of formal power series and is uniquely determined by the requirement that it also has constant term $1$.

Recall \eqref{E:un on P}.

\begin{lem} The affine Weyl group $W_\af$ acts on $\eta(y|a)$
as follows:
 \begin{eqnarray}
s_i(\eta(y|a))&=&
\begin{cases}
\eta(y|a) & (0\le i\le n-1)\\
\Omega(a_n|y)^{-1}\eta(y|a)
&(i=n) 
\end{cases},
\label{eq:si eta}\\
u_n(\eta(y|a)) 
&=&\eta(y|u_n(a))
=\eta(y|a)^{-1}. 
\label{eq:un eta}
\end{eqnarray}
\end{lem}
\begin{proof}
We verify 
\eqref{eq:si eta} for $i=0$.
While the sign choice in the square root remains ambiguous, we can justify the following formal computation
\begin{eqnarray*}
s_0(\eta(y|a))&=&
\Omega(a_1|y)s_\theta\left(\sqrt{\prod_{i=1}^n\Omega(a_i|y)}\right)\\
&=&
\Omega(a_1|y)\sqrt{
\Omega(a_1|y)^{-1}
\prod_{i=2}^n\Omega(a_i|y)}\\
&=&\sqrt{\Omega(a_1|y)^2
\Omega(a_1|y)^{-1}
\prod_{i=2}^n\Omega(a_i|y)}\\
&=&\eta(y|a).
\end{eqnarray*}
In fact, we can verify
$s_0(\eta(y|a)^2)=\eta(y|a)^2$ and see the sign is correct because clearly  $s_0(\eta(y|a))$ has constant term  $1$. The proofs of the remaining cases follow similarly. 
\end{proof}
\begin{prop}
There is a left $\Ax$-module structure on $\hLa_S$ such that the action of the generator $\pi_n$ of $\Sigma$ is given by 
   \begin{equation}
   \pi_n f(y|a)=\eta(y|a)f(y|u_n(a))
 \quad
 \text{for $f(y|a)\in \hLa_S,$}  
 \label{eq:pi n action on hat La}
   \end{equation}
 in particular we have \begin{equation}
\pi_n(\eta(y|a))=1.
\label{eq:pi n eta}
 \end{equation}
\end{prop}
\begin{proof}
It is straightforward to verify that the $W_\af$-module $\hLa_S$ can be extended to a $\Wx$-module via \eqref{eq:pi n action on hat La}. Together with its $S$-module structure, $\hLa_S$ becomes an $\Ax$-module. 
\eqref{eq:pi n eta} is a consequence of \eqref{eq:un eta}.
\end{proof}

We extend the definition of dual affine Grassmannian  Schur \( P \)-functions
to $w=\pi_n v\in \Wxgr$ for 
$v\in \Wafgr$ so that
\begin{equation}
\hP_w^{(n)}(y|a)
=\pi_n \left(\hP_v^{(n)}(y|a)\right).
\end{equation}
In particular, we have
\begin{equation}
\hP_{\pi_n}^{(n)}(y|a)
=\eta(y|a).
\label{eq:P hat pi n}
\end{equation}

Let $\hGn^e$ be the $S$-span of $\hP_w^{(n)}(y|a)\;(w\in \Wxgr).$
\begin{prop}
$H_*^T(\Gr_{G^{\mathrm{ad}}})$ is isomorphic to
$\hGn^e$ as an $S$-algebra and left $\tilde{\A}_\af$-module.
\end{prop}
\begin{proof}
We extend the map $\gamma: \Pet\to \hLa_S$ to a map $\gamma^e:\Petx\to \hLa_S$ by setting 
\begin{equation}
j^e(b) = b\cdot 1 \qquad\text{for $b\in \Petx$.}
\end{equation}
It is straightforward to see that \(\gamma^e\) is an injective homomorphism of \(S\)-algebras and \(\Petx\)-modules, in the same way as in the non-extended case. 
Since by definition \(\hGn^e\) is the image of \(\gamma^e\), the result follows.
\end{proof}

\subsection{Factorization formula}\label{ssec:factorization}
Recall the definition of the element $\kappa_i\in \Wafgr$, given in \S \ref{ssec:factor}.
For type $C$, they are given by
\begin{align}
    \kappa_i
 = w_{(2n-i+1)}^i \qquad
    \text{for $1\le i\le n-1$},
    \quad
\kappa_n=w_{(n,n-1,\ldots,1)}.
\end{align}

\begin{prop}\label{prop:-varpi n} 
We have
$$\hP_{t_{-\varpi_n^\lor}}^{(n)}(y|a)=\eta(y|a)^{-1}\cdot\hP_{\kappa_n}^{(n)}(y|a).$$
Moreover, the function $\hP_{t_{-\varpi_n^\lor}}^{(n)}(y|a)$ is $W$-invariant.
\end{prop}
\begin{proof}
Proposition \ref{prop:kappa} and Lemma \ref{lem:j of anti-dom W invariant}.
\end{proof}

\begin{thm}\label{thm:hat P factor}
Let $v\in W_\af$. 
Assume $v\kappa_i\in W_\af^0$ and $\ell(v\kappa_i)=\ell(v)+\ell(\kappa_i).$ If $i=n$, then $\overline{v}:=\pi_n^{-1} v \pi_n\in \Wafgr,$ and otherwise $v\in \Wafgr.$ 
Moreover, we have 
    \begin{equation}
\hP_{v\kappa_i}^{(n)}(y|a)
    =\begin{cases}
        \hP_{\kappa_i}^{(n)}(y|a)\hP_{v}^{(n)}(y|a) & \text{for $1\le i\le n-1$,}\\\hP_{\kappa_i}^{(n)}(y|a)
        \hP_{\overline{v}}^{(n)}(y|u_n(a)) & \text{for $i=n$.}
    \end{cases}
    \end{equation}
\end{thm}

\begin{proof}This follows from Theorem \ref{thm:factor kappa i}. 
\end{proof}

The functions $\hP^{(n)}_{t_{-\varpi_i^\lor}}(y|a)$ play the role of $\tau$-\emph{functions} of the $G^\vee$-Toda lattice.
The detailed discussion on the Toda lattice and its connection to $QH^*_T(G/B)$ will be given in the forthcoming paper. 

\subsection{Small classes}
\label{ssec:small}
Recall the bijection $\Wafgr\to \Ps_C^{(n)}$ denoted $w\mapsto \la_w$ from \S \ref{SS:relation to infinite rank}.
We have $\ell(w)=|\la_w|.$
Note that the inclusion order $\subset $ for $\mathscr{P}_C^{(n)}$ is not consistent with the Bruhat order of $W_\af^0.$

Recall $a^{(n)}$ from \eqref{E:an}. We have the following result.
\begin{thm}\label{thm:small}
If $w\in W_\af^0$ satisfies $\ell(w)\le 2n$ then $\la_w\in\SP$ and 
\begin{equation}
    \hP_{w}^{(n)}(y|a)=\hP_{\la_w}(y|a^{(n)}).
\end{equation}
\end{thm}

We regard \( W_{\mathrm{af}}^0 \) as a poset with respect to the left weak order \( \le_L \).
This order is generated by the covering relations: specifically, \( w \lessdot v \) if and only if \( v = s_i w \) for some \( i \in I_{\mathrm{af}} \) and \( \ell(v) = \ell(w) + 1 \).

The set \( \SP \) is equipped with the transitive action of $W_\infty$; this was specified in \S \ref{ssec:Cinf}.
Let \( I_{\mathrm{af}} = \{ 0,1,\ldots,n \} \) be the affine Dynkin node set.

Define a map \( \pi: I_\infty \to I_\af \) by
\begin{equation}
\pi(i+2kn) =
\begin{cases}
i & \text{if } 1 \leq i \leq n, \\
2n - i & \text{if } n+1 \leq i \leq 2n,
\end{cases}
\end{equation}
for $k\ge0$. 
This induces an action of the affine Weyl group $W(C_n^{(1)})$ on $\SP$ by
\begin{align}
    s_i \cdot \la = \prod_{j\in \pi^{-1}(i)} s_j \cdot \la
    \qquad\text{for $i\in I_\af$.}
\end{align}

The orbit $W(C_n^{(1)})\cdot \emptyset$ is the set \( \mathscr{C}_{2n}^{\mathrm{sym}} \) of symmetric (self-transpose) \( 2n \)-cores, partitions $\la$ such that $\la^t=\la$ and no hook of $\la$ has size $n$.

\begin{ex} 
Let $n=3$. In the pictures below, the diagonal $\pm j$ for $j\ge0$, is labeled by $i\in I_\af$ where $\pi(j)=i$. 
So the diagonals labeled $i$ are in the set $\pm i + 6\Z$.
We have
\[
s_2 \cdot \begin{ytableau} 0 & 1 & 2 &3\\
1 & 0 & 1\\
2 & 1 \\
3
\end{ytableau} = \begin{ytableau} 0 & 1 & 2 &3 & *(green)2\\
1 & 0 & 1 & *(green) 2\\
2 & 1 \\
3 & *(green) 2\\
*(green) 2
\end{ytableau},\quad
s_3 \cdot \begin{ytableau} 0 & 1 & 2 &3 & 2\\
1 & 0 & 1 & 2 \\
2 & 1 \\
3 & 2\\
 2 
\end{ytableau}=\begin{ytableau} 0 & 1 & 2 &3 & 2\\
1 & 0 & 1 & 2 &*(green) 3\\
2 & 1 \\
3 & 2\\
 2 &*(green) 3
\end{ytableau}
\]
\[
s_3 \cdot \begin{ytableau} 0 & 1 & 2 &3\\
1 & 0 & 1\\
2 & 1 \\
3
\end{ytableau} = \begin{ytableau} 0 & 1 & 2 &*(lightgray)\\
1 & 0 & 1 \\
2 & 1\\
*(lightgray)
\end{ytableau},\quad
s_1 \cdot \begin{ytableau} 0 & 1 & 2 \\
1 \\
2 
\end{ytableau} = \begin{ytableau} 0 & 1 & 2 \\
1 \\
2 
\end{ytableau}
\]
\end{ex}

A fundamental result we use is the following.

\begin{prop}[Hanusa and Jones {\cite{HJ}}]\label{prop:HJ}
There is a poset isomorphism
\begin{equation}
(W_\af^0, \le_L) \cong 
(\mathscr{C}_{2n}^{\mathrm{sym}}, \subset).
\end{equation}
Moreover, this isomorphism is compatible with the actions of $W_\af$ on both sides.
\end{prop}
Let $\mathfrak{c}: (W_\af^0,\le_L)\rightarrow (\mathscr{C}_{2n}^{sym},\subset)$ denote the isomorphism.

We shall prove the following.
\begin{prop}\label{prop:poset}
There is an isomorphism of posets
\begin{equation}
(W_\af^0,\le_L)_{\le 2n}
\cong(\SP,\subset)_{\le 2n}.
\end{equation}    
Moreover, the isomorphism is label compatible via $\pi$.
\end{prop}

The key fact we need to prove Proposition \ref{prop:poset} is the following.
\begin{lem}\label{lem:key comb}
Let \(\lambda\) be a strict partition with \(|\lambda| < 2n\).
Then, the $C_n^{(1)}$- residues of addable boxes lying on or above the diagonal are mutually distinct.   
\end{lem}
\begin{proof}

If \(\lambda_1 \leq n\), then the \(C_n^{(1)}\)-residue of each addable box coincides with its residue in \(\SP\).
Therefore, there are no repetitions among the residues of the addable boxes.

Suppose \(\lambda_1 > n\).
Write \(\lambda_1 = 2n - i\) for some \(1 \leq i \leq n-1\).
Since \(|\lambda| < 2n\), it follows that \(\lambda_k < i\) for all \(k \geq 2\).
The addable box in the first row has \(C_n^{(1)}\)-residue \(i\).
Meanwhile, any box in a lower row with the same \(C_n^{(1)}\)-residue \(i\) would need to be farther to the right than column \(i\), but such boxes are not addable because \(\lambda_k < i\) for all \(k \geq 2\).
Thus, there is no repetition between the \(C_n^{(1)}\)-residue of the addable box in the first row and the residues of addable boxes below.

Moreover, among the addable boxes below the first row, no repetition of \(C_n^{(1)}\)-residues occurs, because they correspond to distinct residues as in the case for \(\SP\) with smaller entries.

Therefore, the residues of all addable boxes are mutually distinct.
\end{proof}

\begin{proof}[Proof of Proposition \ref{prop:poset}] 
We claim that $(\mathscr{P}_C^{(n)})_{\le 2n}=\SP_{\le 2n}.$ 
In fact, for each $i$ such that $n<i\le 2n$, $i$ does not appears more than once in $\la$, because $|\la|\le 2n. $ Hence the claim holds. 
Since we have a bijection $W_\af^0\cong\mathscr{P}_C^{(n)}$, it yields a bijection \begin{equation}
    (W_\af^0)_{\le 2n}\cong \SP_{\le 2n}.\label{eq:main bij}
\end{equation}
We can describe the bijection \eqref{eq:main bij} as follows.
For any partition \(\lambda\), let \(\lambda^+\) denote the strict partition formed by the boxes lying on or above the main diagonal.
By construction of the map $\mathfrak{c}$, it is straightforward to see that the bijection \eqref{eq:main bij} is given by $W_\af^0\ni w\mapsto \mathfrak{c}(w)^+\in \SP.$

We need to prove that this is an isomorphism of posets.
Let $w,v\in W_\af^0$ such that $w\lessdot v$, with 
$v=s_i w$, and $\ell(v)\le 2n.$
Let $\la=\mathfrak{c}(w), \mu=\mathfrak{c}(v).$
The symmetric $2n$ core $\mu$ is obtained from $\la$ by adding all $i$-addable boxes to $\la$. 
From Lemma \ref{lem:key comb}, there is exactly one $i$-addable box of $\la$ lying on or above the diagonal. 
This means that the strict partition $\mu^+$ is obtained by adding the corresponding single box to $\la^+$. 
Hence $|\mu^+|=|\la^+|+1.$ 

On the other hand, suppose that \(\alpha\) and \(\beta\) are strict partitions such that 
\(\alpha \subset \beta\), \(|\alpha|,|\beta| \leq 2n\), and \(\beta\) is obtained from \(\alpha\) by adding a box of \(C_n^{(1)}\)-residue \(i\).
By Lemma \ref{lem:key comb}, there is exactly one such box.
Let \(\lambda\) and \(\mu\) be the corresponding elements of \(\mathscr{C}_{2n}^{\mathrm{sym}}\) associated with \(\alpha\) and \(\beta\), respectively, so that \(\lambda^+ = \alpha\) and \(\mu^+ = \beta\).
Then \(\mu\) is obtained from \(\lambda\) by the action of \(s_i\): either by adding a single box on the main diagonal or by adding a pair of boxes in symmetric positions off the diagonal, both of \(C_n^{(1)}\)-residue \(i\), depending on the case. Then the corresponding elements $w, v\in W_\af^0$ satisfy $w\lessdot v$ by Proposition \ref{prop:HJ}.
\end{proof}



We record the following consequence.  
\begin{cor}
If $w\in W_\af^0$ satisfies $\ell(w)\le 2n$, then 
$\ell(w)=|\mathfrak{c}(w)^+|.$    
\end{cor}

\begin{proof}[Proof of Theorem \ref{thm:small}]
Suppose $w,v\in W_\af^0$ satisfies $w\lessdot v$ and $v=s_i w$ with $i\in I_\af$. If $\ell(v)\le 2n$, there exists unique $\hat{i}\in \{0,1,\ldots,2n-1\}$ and unique $\hat{i}$-addable box of $\mathfrak{c}(w)^+$ such that $\mathfrak{c}(v)^+$ is obtained from $\mathfrak{c}(w)^+$ by adding the box in $\SP$.

 It is straightforward to see for all $f(y|b)\in \hat{\Lambda}_{\Z[a]}$, we have \begin{equation}
\pi(A_{\hat{i}}f(y|b))=A_{i}f(y|b^{(n)}),\label{eq:TT}
 \end{equation}
 where $\pi: \Z[a]\rightarrow S$ is the ring homomorphism defined by sending $a_i$ to $a_i^{(n)}$.

Suppose $w,v\in W_\af^0$ satisfies $w\lessdot v$ and $v=s_i w$ such that $\ell(v)\le 2n.$
We use induction on $\ell(v)$. 
By induction, we have
$$
\hP_w^{(n)}(y|a)=
\hP_{\la_w}(y|a^{(n)}).
$$
By this equation together with \eqref{eq:TT}, we obtain
\begin{align*}
\hP_{\la_v}(y|a^{(n)})&=\pi(\hP_{\la_v}(y|a))=\pi(A_{\hat{i}}\; \hP_{\la_w}(y|a))\\
&=A_{i}\; \hP_{\la_w}(y|a^{(n)})
=A_i\;\hP_w^{(n)}(y|a)
=\hP_v^{(n)}(y|a).
\end{align*}
This completes the proof.
\end{proof}

{
\section{Centralizer family for $\mathrm{SO}_{2n+1}(\mathbb{C})$ }\label{sec:cent}
This section focuses on the centralizer family $\mathscr{Z}_{G^\vee}$ for $G^\vee=\mathrm{SO}_{2n+1}(\mathbb{C})$.

\subsection{Notation for Lie algebra $\fg^\vee=\fso_{2n+1}$}
\label{sec:fgv}
Denote by $M_{r,s}(\C)$ the set of $r\times s$ matrices with entries in $\C$ and $M_r(\C)=M_{r,r}(\C)$.
Let $E_{k,l}\in M_{r,s}(\C)$ be the matrix that is zero in all positions except for a $1$ at position $(k,l)$. Let
$J_{B_n}=\sum_{i=1}^{2n+1}(-1)^{i+1}E_{i,2n+2-i}$. For example, 
\begin{align*}
	J _{B_2}= \begin{pmatrix}
		0&0&0&0&1\\
		0&0&0&-1&0\\	
		0&0&1&0&0\\	    
		0&-1&0&0&0\\	    
		1&0&0&0&0
	\end{pmatrix}
\end{align*}

\begin{rem} \label{R:JB} This choice of $J_{B_n}$ has the property that the resulting centralizer family $\cZGv$ is a subvariety of the centralizer family for type $SL_{2n+1,\C}$ studied in \cite{LS:double Kostant}. See the proof of Proposition \ref{prop:Zeq}.
\end{rem}

Let $\mathfrak{g}^\vee=\fso_{2n+1}(\C)$ be the set of matrices $X\in M_{2n+1}(\C)$ such that
\begin{align}\label{E:so def}
	{}^t\!X J_{B_n} + J_{B_n} X &= 0
\end{align}
where ${}^t\!X$ is the transpose of $X$.

\begin{rem} \label{R:so}
Let $X=(c_{ij})\in M_{2n+1}(\C)$. Then 
$X\in \fso_{2n+1}$ if and only if
\begin{align}\label{E:so condition}
c_{2n+2-j,2n+2-i}=(-1)^{i+j-1}
c_{i,j}.
\end{align}
In particular, the anti-diagonal entries are all zero.
\end{rem}

The following elements form a set of Chevalley generators for $\mathfrak{g}^\vee$:
\begin{align*}
f_i^\vee&=E_{i+1,i}+E_{2n+2-i,2n+1-i}&\quad&\text{for $1\le i\le n$} \\
e_i^\vee&= {}^t(f_i^\vee) &&\text{for $1\le i\le n-1$} \\
e_n^\vee &= 2\, {}^t(f_n^\vee) \\
h_i^\vee&=E_{i,i}-E_{i+1,i+1}
+E_{2n+1-i,2n+1-i}-E_{2n+2-i,2n+2-i}&&\text{for $1\le i\le n-1$} \\
h_n^\vee&=2(E_{n,n}-E_{n+2,n+2}).
\end{align*}



There is a perfect pairing $\fhv\times\fh$ such that $\langle\alpha_i^\lor,\alpha_j\rangle=a_{ij}$, where $(a_{ij})_{i,j\in I}$ is the Cartan matrix of $\fg$. Thus, we have a canonical linear isomorphism $\fh \cong (\fhv)^*$.
There is an embedding \begin{equation}
    \fh\hookrightarrow (\mathfrak{g}^\lor)^*\label{eq:h in dual gv}
\end{equation}
such that $h\in \fh$ maps to the linear function on $\fg^\lor$, that vanishes on the root subspaces of $\fg^\lor$, and on $\fhv$ is given by the pairing with $\fh.$
Under the embedding \eqref{eq:h in dual gv}, for $1\le i\le n$, we have $E_{ii}-E_{2n+1-i,2n+1-i}\in\fh\mapsto \eps_{ii}\in (\fgv)^*.$

\subsection{The group $SO_{2n+1}(\C)$}
The group $\Gv$ dual to $G=\Sp_{2n}(\C)$ is $\mathrm{SO}_{2n+1}(\C)$ which we realize as the set of all $g\in \mathrm{GL}_{2n+1}(\C)$ such that
\begin{align}\label{E:O condition}
\trsp{g}J_{B_n} g &= J_{B_n}\qquad\text{and} \\
\label{E:det 1 condition}
\det(g)&=1.
\end{align}
Its Lie algebra is given by the above construction of $\fso_{2n+1}$.
The Borel subgroup $\Bv$ of $\Gv$ consists of the upper triangular matrices in $\Gv$.

\subsection{Centralizer family $\cZGv$}

We consider $\fh=(\fhv)^*$ as the subspace of $(\fgv)^*$ consisting of functionals on $\fg^\vee$ that vanish on all root subspaces of $\fg^\vee$. 
Define 
\[
E=\sum_{i=1}^n(f_i^\vee)^*\in (\fgv)^*,
\]
where $(f_i^\vee)^*$ denotes the linear function on $\fgv$ which takes value $1$ on $f_i^\vee$, and $0$ on all other root spaces of $\fgv$.  

$\Gv$ acts on $(\fgv)^*$ by the coadjoint action 
$$
(g\cdot \varphi)(x)=\varphi(\mathrm{Ad}(g^{-1})\cdot x)
\quad\text{for $g\in \Gv$, $\varphi \in(\fgv)^*$, and 
$x\in \fgv$.}
$$
The centralizer family associated with $G^\vee$ is, by definition
\[
\mathcal{Z}_{\Gv}=\{(b,h)\in \Bv\times\fh\;|\;b\cdot(-E+h)=-E+h \}.
\]
We consider $\mathcal{Z}_{\Gv}$ as a variety over $\fh$ via the second projection. 

Rather than working with functionals, we shall work with matrices instead. Consider the pairing $\fgv \times \fgv\to \C$ given by  $(A,B)\mapsto\langle A|B\rangle:= \frac{1}{2}\trace(AB)$. Explicitly for $A=(a_{ij})$ and $B=(b_{ij})\in \fgv$,
 \begin{equation}
   \langle A|B\rangle=\sum_{1\le i,j\le n,\;i+j\le 2n+1}a_{ij}b_{ji}.  
 \end{equation}
 We have
 \begin{equation}
    \langle A|\Ad(g)B\rangle 
    = \langle \Ad(g^{-1})A|B\rangle\quad 
    \mbox{for}\quad
    A,B\in \fgv,\quad g\in \Gv.
 \end{equation}
The pairing $\langle\cdot|\cdot\rangle$ induces an isomorphism
$(\fgv)^*\cong \fgv$ such that $(f_i^\vee)^*\mapsto E_{i,i+1}+E_{2n+1-i,2n+2-i}$, $E=\sum_{i=1}^n(f_i^\vee)^*\mapsto \sum_{i=1}^{2n-1}E_{i,i+1}\in \fgv$.
For $h\in \fh$, the element $-E+h\in(\fgv)^*$ corresponds to the following matrix
\begin{align}\label{E:-E+h}
L_0(h):=
\left[\begin{array}{ccccccccc}
a_{1}(h) & -1 & 0 & 0 & 0 & 0 & 0 &0 & 0
\\
0  & a_{2}(h) & -1 & 0 & 0 & 0 & 0 &0 & 0
\\
 0 & 0 & \ddots& \ddots & 0 & 0 & 0 &0 & 0
\\
 0 & 0 & 0 & a_n(h) & -1 & 0 & 0 &0 & 0
\\
 0 & 0 & 0 &0 & 0 & -1 & 0 &0 & 0
\\
 0 & 0 & 0 & 0 & 0 & -a_{n}(h) & -1 &0 & 0
\\
 0 & 0 & 0 & 0 & 0 & 0 & \ddots&\ddots & 0\\
 0 & 0 & 0 & 0&0 & 0 & 0 & -a_{2}(h)&-1 \\0& 0 &0&0&
   0 & 0 & 0 & 0 & -a_{1}(h) 
\end{array}\right]
\end{align}
in $\fgv$ via the identification $(\fgv)^*\cong \fgv$ given by the pairing $\langle\cdot |\cdot\rangle$ defined above. 

Thus, for $(b,h)\in B^\vee\times \fh$, we have $(b,h)\in \cZGv$ if and only if $\langle bL_0(h)b^{-1}|X\rangle=  \langle L_0(h)| X\rangle$ for all $X\in\fgv,$ which is equivalent to the matrix equation
\begin{equation}
 bL_0(h)=L_0(h)b. \label{eq:bL=Lb}  
\end{equation}

\subsection{$W$-action on $\cZGv$}
Let $U^\vee_{-}$ be the subgroup of $G^\vee$ generated by $\exp(uf_i^\vee)$ for $i\in I,\,u\in \C$.
Recall that $\{\al_i\;|\;i\in I\}$ is the set of simple roots for $\fg=\mathfrak{sp}_{2n}(\C)$, and the Cartan subalgebra $\fh$ of $\fg$ is considered to be a subspace of $(\fgv)^*.$

\begin{prop} [\cite{Pet}]\label{P:W action}
There is a unique family of morphisms $u_w^\vee: \fh\to U^\vee_{-}$ for $w\in W$, such that for all $h\in\fh$
\begin{align}
\label{E:u s}
u_{s_i}^\vee(h) &= \exp(-\pair{\alpha_i}{h} f_i^{\vee})&\qquad&\text{for $i\in I$} \\
\label{E:uvw}
    u_{vw}^\vee(h) &= u_v^\vee(w(h)) u_w^\vee(h)&&\text{for all $v,w\in W$.}
\end{align}
Moreover
\begin{align}
\label{E:uw}
    u_w^\vee(h) \cdot (-E+h) &= -E+ w(h)\qquad\text{for all $w\in W$ and $h\in\fh$}
\end{align}
where $\cdot$ denotes the coadjoint action of $G^\vee$ on $(\fgv)^*$.
\end{prop}
\begin{proof} 
For $w=s_i$, one may show by direct computation that \eqref{E:uw} holds using \eqref{E:u s}. It is immediate that if \eqref{E:uw} holds for $v$ and $w$, then it holds for $vw$ using \eqref{E:uvw}. 
Thus, the family $\{u_w^\vee\mid w\in W\}$ is unique if it exists. The value of $u_w^\vee$ seems to depend on the iterated factorizations of $w$ into shorter factors, but it follows from \cite[Theorem 2.6]{Ko} that $u_w^\vee$ is independent of these choices.
\end{proof}

Note that the construction works for any complex semisimple Lie 
algebra $\fgv.$
\begin{ex} 
In this example we use $G^\vee=SO_5(\C)$ and its Lie algebra $\fgv$ with the Chevalley generators $\{e_i^\vee,h_i^\vee,f_i^\vee\;|\;i=1,2\}$ given in \S~\ref{sec:fgv}.
We have
\begin{align*}
 -E+h &= \begin{pmatrix}
     a_1 &-1&0&0&\!0\\
     0&a_2&-1&0&\!0\\
     0&0&0&-1&\!0 \\
     0&0&0&-a_2&\!-1\\
     0&0&0&0&\!-a_1
 \end{pmatrix},\;
u^\vee_{s_1}(h)
=\begin{pmatrix}
 1&0&0&0&0 \\
 \!\!-a_1\!+\!a_2\!\!\!&1&0&0&0\\
 0&0&1&0&0\\
 0&0&0&1&0\\
 0&0&0&\!\!\!-a_1\!+\!a_2\!\!\!&1
\end{pmatrix},
\\
u^\vee_{s_2}(h)&=
\begin{pmatrix}
 1&0&0&0&0 \\
 0&1&0&0&0\\
 0&-2a_2&1&0&0\\
0&2a_2^2&-2a_2&1&0\\
 0&0&0&0&1
\end{pmatrix}
\end{align*}
where $\alpha_1=a_1-a_2,\alpha_2=2a_2$ are the simple roots of 
$\fg=\mathfrak{sp}_{4}(\C).$
Here we denote $a_i(h)$ simply by $a_i.$
Since 
$\theta=2a_1, s_\theta=s_1s_2s_1,
$, 
we have
\begin{equation*}
u_{s_\theta}^\vee(h)=
\begin{pmatrix}
  1 & 0 & 0 & 0 & 0 \\
  -2a_1 & 1 & 0 & 0 & 0 \\
  2a_1(a_1-a_2) & -2a_1 & 1 & 0 & 0 \\
  -2a_1^2(a_1-a_2) & 2a_1^2 & -2a_1 & 1& 0 \\
  2a_1^2(a_1-a_2)(a_1+a_2) & -2a_1^2(a_1+a_2) & 2a_1(a_1+a_2) & -2a_1 & 1
\end{pmatrix}.
\end{equation*}
\end{ex}
\bigskip

Define the action of $w\in W$ on $\Gv\times \fh$ by
\begin{align*}
  w(g,h) = (u^\vee_w(h)\cdot g, w(h)),\quad (g,h)\in \Gv\times \fh,
\end{align*}
where $\cdot$ is the adjoint action of $G^\vee$ on itself.

\begin{lem}[\cite{Pet}] \label{lem: W stabilizes contralizer}
The action of $W$ on $\Gv\times \fh$ described above stabilizes $\cZGv$.
\end{lem}
\begin{proof} Let $w\in W$ and $(b,h)\in \cZGv$. Then we compute:
\begin{align*}
(u^\vee_w(h) b u_w^\vee(h)^{-1}) \cdot (-E+w (h)) &= 
u^\vee_w(h) b \cdot (-E+h) \\
&= u^\vee_w(h) \cdot (-E+h)\\
&= -E+w(h).
\end{align*}
This shows that $w(b, h) \in \cZGv$, as required.
\end{proof}

Let us compute the $W$-action on $\C[\cZGv]$ given by Lemma \ref{lem: W stabilizes contralizer}. 
Let $Z$ be the upper triangular matrix with entries $z_{ij}$. Then the matrix  $s_k(Z)$ with entries $s_k(z_{ij})$ is given by 
\begin{equation}\label{eq: W action on Z}
s_k(Z)=
\exp(-\langle h,\alpha_k\rangle f_k^\vee )
Z
\exp(\langle h,\alpha_k\rangle f_k^\vee )
=
    \exp(-\mathrm{ad}(\langle h,\alpha_k\rangle f_k^\vee )(Z).
    \end{equation}
\begin{ex}\label{ex: Z acts on C[cent]} 
Let $\fgv=\mathfrak{so}_5(\C).$
$s_1(z_{ij})$ is given as follows:
The first two rows
$$
    \left(\begin{array}{ccccc}
        z_{11} & z_{12} & z_{13} &z_{14} &z_{15}
        \\
        0 & z_{22} & z_{23}&
        z_{24}& z_{25}
    \end{array}
    \right)
    $$
    are sent to 
    $$\left(
    \begin{array}{ccccc}
 z_{11}+\alpha_1 z_{12} & z_{12} & z_{13} &z_{14}+\alpha_1  z_{15} &z_{15}
        \\
        0 & z_{11} & z_{23}-\alpha_1 z_{13} &
        z_{24}+\alpha_1 
       (z_{25}-z_{14}) 
       -\alpha_1 ^2z_{15}& z_{25}-\alpha_1 z_{15}
    \end{array}\right).
    $$
    Other coordinate functions are determined by 
    the anti-diagonal symmetry. 
    Similarly, 
    $s_2(z_{ij})$ is given as 
    $$
    \left(
    \begin{array}{ccccc}
        z_{11} & z_{12}+\alpha_2 z_{13}+{2}^{-1}
        \alpha_2^2 z_{14}
        & z_{13}+\alpha_2 z_{14}&
        z_{14}& z_{15}\\
        0 & z_{22}+\alpha_2 z_{23}
        +{2}^{-1}
        \alpha_2^2z_{24}& z_{23}+\alpha_2 z_{24}&
        z_{24} & z_{25}
    \end{array}\right).
    $$
\end{ex}

\begin{prop} \label{P:W on diagonal matrix entries}
Under the action of $W$ on $\C[\cZGv]$, each element $w\in W$ acts on the diagonal entries $z_{11}^{\pm 1},\ldots,z_{nn}^{\pm 1}$ as a signed permutation.
\end{prop}
\begin{proof}
Straightforward. 
\end{proof}
\subsection{$\nilA$-action on $\C[\cZGv]$}
\label{SS:nil-Hecke action on centralizer family}
The $W$-action on $\C[\cZGv]$ extends to a $W_\af$-action by setting
\begin{align}\label{E:translation acting on matrix entries}
t_{\eps_i}f:=z_{ii}f\quad 
\text{for $f\in \C[\cZGv]$,} 
\end{align}
or equivalently, by letting the generator $s_0$ of $W_\af$ act by
$$
s_0 f:=z_{11}s_\theta(f)\quad 
\text{for $f\in \C[\cZGv].$} 
$$
We will show that this further extends to an $\A_\af$-action.
We now state the following important fact:
\begin{thm}\label{thm:flat}
   $\C[\cZGv]$ is flat as an $S$-module. 
\end{thm}
The proof is given in Section \ref{sec:flat}.

\begin{lem}\label{lem: div alpha}
Let $\alpha\in \Phi$ and $f\in \C[\Bv\times \fh]$. 
Then $f-s_\alpha f$ is divisible by $\alpha.$
\end{lem}
\begin{proof}
It is clear from \eqref{eq: W action on Z} that $f-s_i(f)$ is divisible by $\alpha_i.$
    Let $w\in W, i\in I$ such that $w(\alpha_i)=\alpha.$
    We have $w^{-1}f-s_i w^{-1}f=\alpha_i g$ for some $g\in \C[\Bv\times \fh].$
    Then $f-s_\alpha f=f-ws_i w^{-1}f=w(\alpha_i) w(g).$
\end{proof}

\begin{lem} \label{L:div alpha0}
For $f\in \C[\Bv\times\fh]$, $f-s_0f$ is divisible by $\alpha_0=-\theta$.
\end{lem}
\begin{proof} We have 
\begin{align*}
    (1 - s_0) f &= f - t_{\theta^\vee} f + t_{\theta^\vee} f - t_{\theta^\vee} s_\theta f \\
    &= (1-z_\theta) f + t_{\theta^\vee} (1 - s_\theta) f
\end{align*}
where $z_{\theta}=z_{11}\in \C[\cZGv]$ is by definition the element by which $t_{\theta^\vee}$ acts.

By Lemma \ref{lem: div alpha} $(1-s_\theta) f$ is divisible by $\theta$.
By \eqref{eq:EqtypeA} we have $z_{n,n}-z_{n+1,n+1} = (b_n-b_{n+1})z_{n,n+1} = a_n z_{n,n+1}$
so that $z_{n,n}-1  = a_n z_{n,n+1}$. Acting by any Weyl group element $w$ such that $w(a_n)=a_1$ (that is, sending $\alpha_n$ to $\theta$) we have $z_{1,1} - 1 = a_1 w(z_{n,n+1})$. Therefore $1-z_{11}$ is divisible by $\theta$ in $\C[\cZGv]$.
It follows that $(1-s_0)f$ is divisible by $\theta$ as required.
\end{proof}

Let $f\in \C[\Bv\times \fh].$ 
By Lemmas \ref{lem: div alpha} and \ref{L:div alpha0}, for $i\in I_\af$, we can define the linear endomorphism of $\C[\Bv\times \fh]$ by
\begin{equation}
   A_i: f\mapsto \frac{f-s_i f}{\alpha_i}\quad 
\text{for $f\in    \C[\Bv\times \fh]$}.
\end{equation}
\begin{cor}[\cite{Pet}]\label{cor:Aaf acts on O}
For $i\in I_\af$ the operator $A_i$ preserves the defining ideal $\II$ of  $\cZGv$ and is therefore well-defined on $\C[\cZGv]$. In particular the operators $A_i$ define an $\A_\af$-module structure on $\C[\cZGv]$.
\end{cor}
\begin{proof}
Let $f\in \II$. By Lemma \ref{lem: W stabilizes contralizer} we have $s_i(f)\in \II$. 
Therefore we have $\alpha_iA_i f=f-s_i(f)\in \II$.
Because $\C[\cZGv]$ is a flat $S$-module (Theorem \ref{thm:flat}), the linear endomorphism of $\C[\cZGv]$ given by the multiplication with $\alpha_i$ is injective.
Thus, we can deduce $A_i f\in \II$. 
\end{proof}
}

\subsection{The isomorphism $\beta:\C[\mathcal{Z}_{\Gv}] \rightarrow \hGn$}
The main goal of this section is to prove the following result.
\begin{thm}\label{thm:beta} Let 
\begin{equation}\label{eq:b seq}
    (b_1,\ldots,b_{2n+1})=(a_1,a_2,\ldots,a_n,0,-a_n,\ldots,-a_1).
\end{equation}
There is a $S$-algebra and $\A_\af$-module isomorphism 
$\beta:\C[\mathcal{Z}_{\Gv}]
\rightarrow
\hGn
$
defined by
\begin{equation}\label{E:beta matrix entry}
    \beta(z_{ij})=
(-1)^{j-i} \hq_{j-i}(y|b_{i},\ldots,b_j).
\end{equation}
In particular
\begin{align}\label{E:beta diag}
\beta(z_{ii})=\Omega(a_i|y)\qquad\text{for $1\le i\le n$.}
\end{align}
\end{thm}
\subsection{Construction of the map $\beta$}

Let $t=(t_1,t_3,\ldots,t_{2n-1})\in \mathbb{C}^n$ be arbitrary. 
We can easily see that $L_0^{2i-1}\in \fbv$ for $1\le i\le n.$
We define
\begin{equation}\label{eq:b=exp()}
g(t)=\exp\left(2\sum_{i=1}^nt_{2i-1} L_0^{2i-1}\right)\in B^\vee.
\end{equation}
Then we have \begin{equation}\label{eq:exp L0}
g(t)\,L_0=L_0\,g(t).\end{equation}
Therefore $(h,g(t))\in \cZGv.$
Thus, we have a holomorphic map 
\begin{equation}
   \fh\times \C^n\rightarrow \cZGv\quad\text{for $(h,t)\mapsto (h,g(t)).$} 
\end{equation}

\begin{lem}\label{lem:coordinate=hQ}
Then $(i,j)$ entry of $g(t)$ defined by \eqref{eq:b=exp()} is equal to
\begin{equation}
\label{eq:hatb}
(-1)^{j-i} \hat q_{j-i}(y|b_{i},\ldots,b_j)
\end{equation}
for all $(i,j)$ such that $1\le i\le j\le 2n+1.$ 
\end{lem}

\begin{proof}
Let
 $D=\mathrm{diag}(b_1,b_2,\dots,b_{2n+1})$ and 
$
E_\eps=
\begin{pmatrix}
    0&1\\
    &0&\ddots\\
    &&\ddots & 1\\
    \eps & & &0
\end{pmatrix}
$, where $\eps$ is an auxiliary parameter.
Define $L_\eps=D-E_\eps$ and $g_\eps(t)=\exp\left(2\sum_{i=1}^nt_{2i-1} L_\eps^{2i-1}\right)$.
By \eqref{eq:generating_q}, $g_\eps(t)$ expands as
\[
\begin{aligned}
    g_\eps(t)
    &=\sum_{k=0}^\infty {q}_k(t)L_\eps^k
    =\sum_{k=0}^\infty {q}_k(t)(D-E_\eps)^k.
\end{aligned}
\]
Let $D_{[i]}:=E_\eps^i DE_\eps^{-i}=\mathrm{diag}(b_{1+i},b_{2+i},\dots,b_{2n+1+i})$, where $b_{2n+1+k}=b_k$.
Let us expand the power of $D-E_\eps$ into signed sums given by sequences of $D$'s and $E_\eps$. 
For such a sequence, let $p$ be the number of $D$'s and $r$ the number of $E_\eps$'s, and for all $s$, let $i_s$ be the number of $E_\eps$'s to the left of the $s$-th $D$. Summing over all such sequences, we have
\[
g_\eps(t)=\sum_{r=0}^\infty 
(-1)^r\sum_{p=0}^\infty 
{q}_{r+p}(t)\sum_{0\leq i_1\leq i_2\leq \dots\leq i_p\leq r}
D_{[i_1]}D_{[i_2]}\dots D_{[i_p]}
\cdot E_\eps^r.
\]
Taking the limit $\eps\to 0$, we have $E_0^{2n+1}\to 0$ and 
\[
g(t)=g_0(t)=\sum_{r=0}^{2n} 
(-1)^r\sum_{p=0}^\infty 
{q}_{r+p}(t)
\sum_{0\leq i_1\leq i_2\leq \dots\leq i_p\leq r}
D_{[i_1]}D_{[i_2]}\dots D_{[i_p]}
\cdot E_0^r.
\]
Since the $(i,i)$ entry of 
$
\sum_{0\leq i_1\leq i_2\leq \dots\leq i_p\leq r}
D_{[i_1]}D_{[i_2]}\dots D_{[i_p]}
$
is 
\[
\sum_{0\leq i_1\leq i_2\leq \dots\leq i_p\leq r}b_{i+i_1}b_{i+i_2}\dots b_{i+i_p}
=h_p(b_{i},\dots,b_{i+r}),
\]
the $(i,j)$ entry of $g(t)$ is equal to
\[
(-1)^{j-i}\sum_{p=0}^\infty h_p(b_{i},\dots,b_{j}){q}_{j-i+p}(t).
\]
\end{proof}

\begin{ex}
    Let $n=2.$
    The matrix $g$ is equal to the following:
\begin{equation*}
\begin{pmatrix}
     \hq_0(y|a_1)\!\! & -\hq_1(y|a_1,a_2)\!\!& \hq_2(y|a_1,a_2)\!\!&
     -\hq_3(y|a_1,a_2,-a_2)
     \!\!&\hq_4(y|a_1,a_2,-a_2,-a_1)\\
      & \hq_0(y|a_2) &
      -\hq_1(y|a_2)
      & \hq_2(y|a_2,-a_2) & 
      -\hq_3 (y|-a_1,-a_2,a_2)\\
      & & 1 & -\hq_1(y|-a_2)&
      \hq_2(y|-a_1,-a_2)\\
      & & & \hq_0(y|-a_2)
      & -\hq_1(y|-a_1,-a_2)\\
      &&&&\hq_0(y|-a_1)
    \end{pmatrix}.
\end{equation*}
\end{ex}

\begin{prop}
There is an $S$-algebra homomorphism 
$\beta:\C[\cZGv]\rightarrow \hLa_S.$
\end{prop}
\begin{proof}  The $(i,j)$ entry of $g(t)$ is interpreted as an element of $\hLa_S$ via the correspondence given in Lemma \ref{lem:q and p}. 
By \eqref{eq:exp L0}, we have $(h,g(t))\in \cZGv$. 
Therefore, the family of formal power series $g_{ij}(t)$ for $1\le i\le j\le 2n+1$ satisfies the defining equations of $\cZGv$. 
Thus sending $z_{ij}$ to $\iota(g_{ij}(t))$ gives a well-defined map from $\C[\cZGv]$ to $\hLa_S.$
\end{proof}

\begin{prop} For each $(i,j)$, the function 
 $g_{ij}(t)$ belongs to $\hGn$.
\end{prop}
\begin{proof}
This is a consequence of an explicit formula in Lemma \ref{lem:hatQ in hatP}, which  expresses $g_{ij}(t)$ as an $S$-linear combination of $\hP_{\rho_i}^{(n)}(y|a)\;(1\le i\le 2n).$ 
\end{proof}

{
\subsection{Structure of the localization $\C[\cZGv]^\reg$}
Let 
\begin{align}
    \mathfrak{h}^\reg= \{h\in \mathfrak{h}\mid \alpha(h)\ne0\quad\text{for all roots $\alpha$}\}.
\end{align}
An element \( h \in \mathfrak{h} \) is regular (i.e., its centralizer in \( \mathfrak{g} \) is \( \mathfrak{h} \)) if and only if \( h \in \mathfrak{h}^{\mathrm{reg}} \) \cite[Appendix D]{FH}.
Let $U_{+}^\vee$ be the subgroup of $\Gv$ generated by $\exp(ue^\vee_i)$ for
$i \in  I,\,u\in \C$.

\begin{prop}
{\rm(\cite[Lemma 3.5.2]{Ko}, \cite{Pet})} \label{P:centralizer versus torus}
There is a morphism $M: \mathfrak{h}^\reg \to U_{+}^\vee$
such that 
\begin{align}\label{E:M def}
  M(h) \cdot (-E+h) = h.
\end{align}
In particular, the map 
\[
(b,h) \mapsto ((\Ad\,M(h))(b),h)
\]
defines an isomorphism 
over $\mathfrak{h}^\Delta$ 
from $\cZGv$ to the trivial family with fiber
$T^\vee$.
\end{prop}
\begin{proof} 
The element $-E+h$ given in \eqref{E:-E+h} for type $\SO_{2n+1}(\C)$, is obtained from the corresponding element of type $\SL_{2n+1}(\C)$ (say with general diagonal entries $(d_1,\dotsc,d_{2n+1})$) by restricting Cartan variables as in \eqref{E:b def}.
It follows that we may obtain the formula for $M(h)$ in type $\SO_{2n+1}(\C)$ by restricting the formula in type $\SL_{2n+1}(\C)$, where one may verify that the latter is given by
\begin{align*}
M^{\SL_{2n+1}}(h)_{ij} = \prod_{k=i+1}^j (d_k-d_i)^{-1}
\qquad\text{for $1\le i\le j \le 2n+1$.}
\end{align*}
This proves the existence of $M(h)$. The condition \eqref{E:M def} shows that for $h\in \mathfrak{h}^\Delta$, $(b,h)\in \cZGv$ if and only if $(\Ad\,M(h))(b) \cdot h = h$, that is, $(\Ad\,M(h))(b) \in T^\vee$. 
The last statement uses the fact that the centralizer of a regular element of the Cartan is $T^\vee$ (see, for example \cite[Appendix D]{FH}).
\end{proof}

\begin{ex} For $\Gv=\mathrm{SO}_5(\C),$ we have
\begin{equation}
M(h)= \begin{pmatrix}
1 & -\frac{1}{a_{1}-a_{2}} & \frac{1}{a_{1} \left(a_{1}-a_{2}\right)} & -\frac{1}{a_{1} \left(a_{1}-a_{2}\right) \left(a_{1}+a_{2}\right)} & \frac{1}{2 a_{1}^{2} \left(a_{1}-a_{2}\right) \left(a_{1}+a_{2}\right)} 
\\
 0 & 1 & -\frac{1}{a_{2}} & \frac{1}{2 a_{2}^{2}} & -\frac{1}{2 a_{2}^{2} \left(a_{1}+a_{2}\right)} 
\\
 0 & 0 & 1 & -\frac{1}{a_{2}} & \frac{1}{a_{1} a_{2}} 
\\
 0 & 0 & 0 & 1 & -\frac{1}{a_{1}-a_{2}} 
\\
 0 & 0 & 0 & 0 & 1 
\end{pmatrix}.
\end{equation}
\end{ex}

\begin{cor} \label{C:centralizer functions as group algebra} 
There is an $S^\reg$-algebra isomorphism $S^\reg[T^\vee]\to \C[\cZGv]^\reg$.
\end{cor}
\begin{proof} This follows from Proposition \ref{P:centralizer versus torus}.
\end{proof}
}

\subsection{Proof of Theorem \ref{thm:beta}}

\begin{prop}\label{prop:beta-Waf linear}
The map $\beta$ is an  $\A_\af$-module isomorphism.
\end{prop} 
\begin{proof}
There is the following commutative diagram
\[
\xymatrix{
\mathbb{C}[\cZGv] \ar[r]^{\beta} \ar@{^{(}->}[d] & \hGn \ar@{^{(}->}[d] \\
\mathbb{C}[\cZGv]^\reg \ar[r]^{\beta^\reg} & \hGn^\reg,
}
\]
where the vertical arrows are the natural inclusions into the localized rings. 
Note that the left vertical arrow is injective because $\C[\cZGv]$ is flat over $S$. 
Here, the localization 
\[
\beta^\reg := 1_{S^\reg}\otimes_S \beta \colon \mathbb{C}[\cZGv]^\reg \to \hGn^\reg
:=S^\reg\otimes_S \hGn
\]
is defined by tensoring with \(S^\reg\).
The map $\beta^\reg $ is an isomorphism of \(S^\reg\)-algebras, since both \(\mathbb{C}[\cZGv]^\reg\) and \(\hGn^\reg\) are isomorphic to the group ring \(S^\reg[Q^\lor]\). In particular, \(\beta^\reg\) is determined by
\[
\beta^\reg(z_{ii}) \;=\; \Omega(a_i|y), \quad (1\le i\le n).
\]
The injectivity of $\beta$ follows from the diagram.

There is a natural action of \(W\) on \(S^\reg[Q^\lor]\). For \(w\in W\), the action is given by
\[
c_\xi\, t_\xi \;\mapsto\; w(c_\xi)\, t_{w(\xi)}, \quad \text{for } c_\xi\in S^\reg \text{ and } \xi\in Q^\lor.
\]
This action extends to an action of \(W_\af\) by letting the translation elements \(t_\xi\) (with \(\xi\in Q^\lor\)) act by multiplication by \(t_\xi\). 
By Proposition \ref{P:W on diagonal matrix entries}, \eqref{E:translation acting on matrix entries}, and Proposition \ref{P:Waf on symmetric functions},
the actions of \(W_\af\) on \(\mathbb{C}[\cZGv]\) and on \(\hGn\) coincide with the restriction of the natural action of \(W_\af\) on \(S^\reg[Q^\lor]\). Hence, the map \(\beta\) is \(W_\af\)-linear.

Moreover, since both rings carry natural left \(S^\reg\)-module structures, it follows that \(\beta^\reg\) is an \(S^\reg[W_\af]\)-module homomorphism. Restricting this structure to the appropriate subrings yields an \(\A_\af\)-module map from \(\mathbb{C}[\cZGv]\) to \(\hGn\).

The map $\beta$ is surjective because for $w\in \Wafgr$, we have $
\beta(A_w 1)=A_w\beta(1)=A_w(1)=\hP_w^{(n)}(y|a).$
\end{proof}

\subsection{Pieri classes generates $H_*^T(\Gr_{C_n})$}


\begin{prop}\label{prop: hGSn}
The $S$-subalgebra of $\hat{\Gamma}_S$ generated by $\hq_{j-i}(y|b_i,\ldots,b_j)$ is equal to 
\begin{equation}
   S[\hP_1(y|a^{(n)}),\hP_2(y|a^{(n)}),\ldots,\hP_{2n}(y|a^{(n)})]. 
\end{equation} 
\end{prop}
\begin{proof}
It follows from Lemma \ref{lem:hatQ in hatP}.
\end{proof}

\begin{cor}
$H_*^T(\Gr_{C_n})$ is generated by $\sigma_{\rho_i}(y|a)\;(1\le i\le 2n)$ as an $S$-algebra. 
\end{cor}

\section{Flatness of $\C[\cZGv]$ over $S$}\label{sec:flat}

\subsection{Equations for centralizer family}

{
Let $\{z_{ij}\mid 1\le i\le j\le 2n+1\}$ be the matrix entry functions for the Borel subgroup $B_{GL}$ of $GL_{2n+1}(\C)$. The coordinate ring is given by $\C[B_{GL}] = \C[z_{ii}^{\pm1} \mid 1\le i\le 2n+1][z_{ij}\mid 1\le i<j\le 2n+1]$. The inclusion $B^\vee \subset B_{GL}$ induces a restriction map $\C[B_{GL}]\to \C[B^\vee]$. Let $I_{B^\vee}$ be the ideal of polynomial functions on $B_{GL}$ vanishing on $B^\vee$. Let $g=(z_{ij}\mid 1\le i\le j\le 2n+1)$ 
and let $R= J-\trsp{g} J g$.
We make the notation $\bar{i} = 2n+2-i$ for $1\le i\le 2n+1$.
Explicitly we have
\begin{align}\label{E:R entry}
    R_{ij} &= - \delta_{i,\jj} (-1)^i + \sum_{p=\jj}^i 
    (-1)^p
    z_{p,i} z_{\bar{p},j}
\end{align}
The ideal $I_{B^\vee}$ is generated by the equations
\begin{align}
\label{E:OB condition} 
R_{ij} &= 0 \qquad \text{for $1\le i\le j\le 2n+1$} \\
\label{E:center 1}
  z_{n+1,n+1} &=1.   
\end{align}
In particular we have the ``above antidiagonal" relations $R_{ij}$ with $i+j\le 2n+2$:
\begin{align}\label{E:automatic zero relations}
    R_{ij} &= 0&\qquad&\text{for $1\le i\le n$ and $ 1\le j < \ii$}\\
\label{E: antidiagonal relation}
  R_{i,\ii} &= (-1)^i(   z_{ii}z_{\ii,\ii} - 1)&&\text{for $1\le i\le n$.}
\end{align}

\begin{rem} 
It is known that $\dim B^\vee=n^2+n$. We have that $B^\vee\cong (\C^\times)^n \times \C^{n^2}$ and that $$\C[B^\vee]=\C[z_{ii}^{\pm 1}\mid 1\le i\le n][z_{ij}\mid 1\le i< j\le 2n+1, \,i+j\le 2n+1].$$
In particular, the other matrix entries depend algebraically on these.
\end{rem}
}

The variety $\mathcal{Z}_{\Gv}\subset B_{GL} \times \fh$ is defined by an ideal $\II$ of $$\C[B_{GL}\times \fh] = \C[a_i\mid 1\le i\le n][z_{ij}\mid 1\le i\le j\le 2n+1][z_{ii}^{-1}\mid 1\le i \le 2n+1].$$
$$
\C[\cZGv]=
\C[B_{GL}\times \fh] /I_Z.
$$

Recall from \S~\ref{ssec:Cnil-Hecke} that $a_i\in S$ is a coordinate function on $\fh$, so $S$ is identified with the coordinate ring $\C[\fh]$ of $\fh$, which is also 
$\C\otimes H_T^*(\mathrm{pt})$. 
Recall the notation \eqref{E:b def}. We have
\begin{align}\label{E:b antisymmetry}
  b_{\ii} = - b_{i}\qquad\text{for all $1\le i\le 2n+1$.}
\end{align} Note that, in particular, we have $b_{n+1}=0.$
The centralizer equation \eqref{eq:bL=Lb} can be written as
     \begin{align}
      (b_i-b_j)z_{ij}+z_{i,j-1}-z_{i+1,j} \qquad\text{for $1\le i<j\le 2n+1$.}
      \label{eq:EqtypeA}
      \end{align}
(cf. \cite[(13)]{LS:double Kostant}).
The following is clear from the definition of $\II$.
\begin{prop}\label{prop:Zeq}
The ideal $\II$ is generated by the equations \eqref{eq:EqtypeA} together with \eqref{E:OB condition}, \eqref{E:center 1}.
\end{prop}

 Writing \eqref{eq:EqtypeA} in terms of $a$ variables, we obtain the following equations.    \begin{align}
\label{E:BA}
    &(a_i-a_j)z_{ij}+z_{i,j-1}-z_{i+1,j}=0&&\text{$1\le i<j\le n$} \\
    &a_i z_{i,n+1}+z_{in}-z_{i+1,n+1}=0&&\text{$1\le i\le n$}  \\
    &(a_i+a_j)z_{i,\bar{j}} + z_{i,\bar{j}-1}-z_{i+1,\bar{j}}&&\text{$1\le i\le n$, $i\le j\le n$}
\label{eq:Zeq}
\end{align} 
The above are equation \eqref{eq:EqtypeA} for positions $(i,j$) with $i+j\le 2n+2$, the positions on or above the antidiagonal. 
There are also equations given from the above by replacing $z_{ij}\mapsto z_{\jj\ii},\;a_i\mapsto -a_i$; these $z_{ij}$ live on or below the antidiagonal.

\subsection{Fewer equations}
Many of the defining equations for $\cZGv$ are redundant. We show that only a few of the defining equations for $SO_{2n+1}$ are necessary in the presence of the equations \eqref{eq:EqtypeA}.

{   
\begin{lem}\label{lem:fewer-equations}
For $\nn\le k\le \overline{1}$ we have
\begin{equation}
R_{k-1,k}=b_kR_{kk}.
\label{eq:R=bR}
\end{equation}
\end{lem}
\begin{proof}

By \eqref{E:R entry} and \eqref{eq:EqtypeA}
\begin{eqnarray}
R_{k-1,k}&=& \sum_{p=\kk}^{k-1} (-1)^{p} z_{p,k-1} z_{\pp k} \nonumber\\
     &=& \sum_{p=\kk}^{k-1} (-1)^{p} z_{p+1,k}z_{\pp k}-\sum_{p=\kk}^{k-1} (-1)^{p}(b_p-b_k)z_{pk} z_{\pp k}\label{eq:Rzero}\\
     &=&-\sum_{p=\kk}^{k-1} (-1)^{p}(b_p-b_k)z_{p k} z_{\pp k},\label{eq:Rk-1k}\end{eqnarray}
where the first sum of \eqref{eq:Rzero} is zero because the pair of terms corresponding to $p=i$ and $p=\overline{i+1}$ cancels. 

We divide the sum \eqref{eq:Rk-1k} into the parts corresponding to the decomposition:
\begin{eqnarray*}
[\kk,k-1]=I_1\cup I_2\cup I_3,\; I_1:=\{\kk\},\;
I_2:= [\kk+1,n]\cup [\nn, k-1],\;
I_3:=\{n+1\}\end{eqnarray*}
where $[a,b]=\{i\in \Z\mid a\le i\le b]$ for integers $a,b$ such that $a\le b$.
Let $S_1,S_2,S_3$ be the sum corresponding to $I_1,I_2,I_3.$
By \eqref{E:b antisymmetry}, one sees that 
\begin{eqnarray}
S_1&=&
- (-1)^{\kk}(b_{\kk}-b_k)z_{\kk k} z_{kk}=(-1)^{\kk}2b_kz_{\kk k} z_{kk},\label{eq:S1}\\
S_3&=&(-1)^{n+1}b_kz_{n+1,k}^2.\label{eq:S3}
\end{eqnarray}
     Noting that $[\kk+1,n]$ and $[\nn,k-1]$ are bijective via $p\mapsto \pp$ and
     $b_{\pp}=-b_p$ we have
\begin{eqnarray}
S_2
&=&-\sum_{p=\kk+1}^{n} (-1)^{p}(b_p-b_k)z_{p k} z_{\pp k}\nonumber
     -\sum_{p=\nn}^{k-1} (-1)^{p}(b_p-b_k)z_{pk} z_{\pp k}\nonumber\\
     &=&-\sum_{p=\kk+1}^{n} (-1)^{p}\left( b_p+b_{\pp}\right)z_{pk} z_{\pp k}
     + \sum_{p=\kk+1}^{n} (-1)^{p}2b_kz_{pk} z_{\pp k}\nonumber\\
  &=&   \sum_{p=\kk+1}^{n} (-1)^{p}2b_kz_{pk} z_{\pp k}.
  \label{eq:S2}
     \end{eqnarray}
    We can write $R_{kk}$ as 
            \begin{eqnarray*}
      R_{kk} =
(-1)^{\kk}2  z_{\kk k} z_{kk}
      +2 \sum_{p=\kk+1}^n (-1)^{p} z_{p k} z_{\pp k}
      + (-1)^{n+1} z_{n+1,k}^2.
       \end{eqnarray*}
       Thus \eqref{eq:R=bR} is clear from 
       \eqref{eq:S1}, \eqref{eq:S3}, and \eqref{eq:S2}.
\end{proof}

\begin{lem}For $ i+2\le j$ and $i+j\ge 2n+2$, we have
\begin{equation}
\label{E:connect bottom relations}
    R_{ij} = - R_{i+1,j-1} + (b_{i+1}+ b_j) R_{i+1,j}.
    \end{equation}
\end{lem}
\begin{proof}
First, consider the case $i=\jj$.
Applying \eqref{E: antidiagonal relation} we obtain
\begin{eqnarray*}
&&\quad\, R_{i\ii} + R_{i+1,\ii-1}=R_{i\ii} + R_{i+1,\overline{i+1}} \\ 
&=& (-1)^i z_{ii} z_{\ii\ii} +(-1)^{i+1} z_{i+1,i+1}z_{\overline{i+1},\overline{i+1}} \\
    &=& (-1)^i\left((z_{i+1,i+1} - (b_i - b_{i+1}) z_{i,i+1}) z_{\ii\ii} - z_{i+1,i+1}(z_{\ii\ii} - (b_{\overline{i+1}}-b_{\ii})z_{\overline{i+1},\ii})\right) \\
    &=& (-1)^i\left((b_{i+1}-b_i) z_{i,i+1} z_{\ii\ii} +(b_i-b_{i+1}) z_{i+1,i+1} z_{\overline{i+1},\ii} \right)\\
    &=& (-1)^i(b_{i+1}-b_i) (z_{i,i+1} z_{\ii\ii} -  z_{i+1,i+1} z_{\overline{i+1},\ii}) \\
    &=& (b_{i+1}+b_{\ii}) R_{i+1,\ii}.
\end{eqnarray*}

Assume $i>\jj$. Note that all the entries of $R$ that are involved here are strictly below the antidiagonal, so in particular, there is no constant term. We have
\begin{eqnarray*}
&&R_{ij} + R_{i+1,j-1} \\
&=& 
 \sum_{p=\jj}^i  (-1)^p  z_{pi}  z_{\pp j} + \sum_{p=\overline{j-1}}^{i+1} (-1)^p z_{p,i+1} z_{\pp,j-1} \\
&=& \sum_{p=\jj}^{i}  (-1)^p  (z_{p+1,i+1}-(b_p-b_{i+1})z_{p,i+1})  z_{\pp j} \\
&
+&\sum_{p=\jj+1}^{i+1} (-1)^p z_{p,i+1} (z_{\pp+1,j}-(b_{\pp}-b_j)z_{\pp j})  \\
&=&\sum_{p=\jj}^{i}  (-1)^p  (b_{i+1}-b_p)z_{p,i+1}  z_{\pp j} 
+\sum_{p=\jj+1}^{i+1} (-1)^p (b_j-b_{\pp})z_{p,i+1}z_{\pp j}   \\
&=&
 (-1)^{\jj}  (b_{i+1}-b_{\jj})z_{\jj,i+1}  z_{j j}
+\sum_{p=\jj+1}^{i}  (-1)^p  (b_{i+1}-b_p)z_{p,i+1}  z_{\pp j} \\
&+&\sum_{p=\jj+1}^{i} (-1)^p (b_j-b_{\pp})z_{p,i+1}z_{\pp j}
+(-1)^{i+1} (b_j-b_{\overline{i+1}})z_{\overline{i+1},i+1}z_{b_{\overline{i+1}}, j}\\
&=&   (-1)^{\jj}  (b_{i+1}+b_j)z_{\jj,i+1}  z_{j j}+\sum_{p=\jj+1}^{i}  (-1)^p  (b_{i+1}+b_j)z_{p,i+1}  z_{\pp j}\\
&+& (-1)^{i+1} (b_{{i+1}}+b_j)z_{\overline{i+1},i+1}z_{b_{\overline{i+1}}, j}\\
&=&(b_{{i+1}}+b_j)\sum_{p=\jj}^{i+1}  (-1)^p  z_{p,i+1}  z_{\pp j}\\
&=&(b_{{i+1}}+b_j)R_{i+1,j},
\end{eqnarray*}
where in the third equality we used 
$$
\sum_{p=\jj}^{i}  (-1)^p  z_{p+1,i+1}  z_{\pp j} 
+\sum_{p=\jj+1}^{i+1} (-1)^p z_{p,i+1} z_{\pp+1,j}=0 $$
and the fifth equality we used $b_p+b_{\pp}=0.$ 
\end{proof}
}

\begin{prop} \label{P:fewer equations} 
The ideal $\II$ is generated by the equations \eqref{E:automatic zero relations}, \eqref{E:center 1}, and $R_{kk}$ for $n+2\le k\le 2n+1$.
\end{prop}
\begin{proof} 
We must show that the given equations imply that
$R=0$. 
By \eqref{E:center 1} we have $R_{n+1,n+1}=0$, pictured by the cyan box in Figure \ref{F:relations}, and we have \eqref{E:automatic zero relations}, which are pictured by the gray boxes in Figure \ref{F:relations}. The entries associated with $R_{kk}$ are colored yellow in Figure \ref{F:relations}.

\begin{figure}

\[
\begin{ytableau}
*(lightgray) 11&*(lightgray) 12&*(lightgray) 13&*(lightgray) 14&*(lightgray) 15&*(lightgray) 16&17\\
\none&*(lightgray) 22&*(lightgray) 23&*(lightgray) 24&*(lightgray) 25&26&27\\
\none&\none&*(lightgray) 33&*(lightgray) 34&35&36&37\\
\none&\none&\none&*(cyan) 44&45&46&47\\
\none&\none&\none&\none&*(yellow) 55&56&57\\
\none&\none&\none&\none&\none&*(yellow) 66&67\\
\none&\none&\none&\none&\none&\none&*(yellow) 77\\
\end{ytableau}
\]

\caption{Picture of $R$ for $n=3$}
\label{F:relations}
\end{figure}

It remains to establish the relations corresponding to white boxes in Figure \ref{F:relations}.
To this end, we use Lemma \ref{lem:fewer-equations} and Lemma \ref{P:fewer equations}.
This suffices: every white entry of $R$ vanishes supposing the entries below and below left vanish, and the latter holds by induction. 
\end{proof}

\subsection{More efficient presentation of $\C[\cZGv]$}
Let $\hz_1,\dotsc,\hz_{2n}$ be indeterminates.
We will provide a streamlined  presentation of the $S$-algebra $\C[\cZGv]$ as a quotient ring of $S[\hz_1,\ldots,\hz_{2n}
].$
The idea is to use the matrix entries $z_{1j}\;(2\le j\le 2n+1)$ as a set of generator of $\C[\cZGv]$ as an $S$-algebra. 

We first express all the entries $z_{ij}$ as an $S$-linear combination of the first row entries.
For $1\le i\le j\le 2n+1$ define $y_{ij}
\in \sum_{k=1}^{2n+1}Sz_{1k}$ by induction on $i$ and then on $j$, by 
\begin{align}
  {y}_{1j} &= z_{1j} &\qquad&\text{for $1\le j\le 2n+1$} \\
  {y}_{ij} &= y_{i-1,j-1} + (b_{i-1}-b_j) y_{i-1,j} &&\text{for $2\le i\le j \le 2n+1$.}
\end{align}

\begin{rem} \label{R:yz} The $y_{ij}$ are $z_{ij}$ expressed in first row $z$ variables using \eqref{eq:EqtypeA} and as such satisfy $y_{ij} \equiv z_{ij} \mod \II$.
\end{rem}

To describe the precise coefficients of the $y_{ij}$ we use plethystic notation. For any variable sets $A$ and $B$ and $r\in \Z_{\ge0}$ and $e_r$ the elementary symmetric polynomial and $h_r$ the complete homogeneous symmetric polynomial, let
\begin{align}\label{E:e pleth}
  e_r[A-B] = \sum_{i+j=r} (-1)^j e_i[A] h_j[B]
\end{align}
with $e_0[A-B]=1$ by definition.
If $A=\{a_1,\ldots,a_i\}$ and $B=\{b_1,\ldots,b_j\}$, we denote
$e_k(A-B)$ by $e_k[(a_1+\cdots+a_i)-(b_1+\cdots+b_j)].$

Let $S_+$ be the augmentation ideal of $S$, the ideal generated by $a_i$ for $1\le i\le n$.

\begin{lem}\label{L:y} For $1\le i\le j\le 2n+1$
\begin{align}\label{E:y to z1}
  y_{ij} = \sum_{k=0}^{i-1} e_k[(b_1+b_2+\dotsm+b_{i-1})-(b_{j-i+k+1}+\dotsm+b_{j-1}+b_j)] z_{1,k+j-i+1}.
\end{align}
In particular 
\begin{align}\label{E:leading term of y}
y_{ij} &\in z_{1,j-i+1} + \sum_{k>j-i+1} S_+ z_{1k} \\
\label{E:diagonal y}
    y_{ii} &= z_{11} + (b_1-b_i) z_{12} + (b_1-b_i)(b_2-b_i) z_{13} + \dotsm \\ 
\notag    &+ (b_1-b_i)(b_2-b_i)\dotsm(b_{i-1}-b_i) z_{1i}\qquad\text{for $1\le i\le 2n+1$}\\
\label{E:middle y}
  y_{n+1,n+1} &= 
 z_{11} + a_1 z_{12} + a_1a_2z_{13}+\dotsm + a_1 a_2\dotsm a_n z_{1,n+1}.
\end{align}
\end{lem}
\begin{proof} Equation \eqref{E:y to z1} easily holds for $i=1$.
Suppose \eqref{E:y to z1} holds for $i-1$. We take the coefficient of $z_{1,k+j-i+1}$ in the equation $y_{ij} = y_{i-1,j-1} + (b_{i-1}-b_j) y_{i-1,j}$. 
It suffices to show that
\begin{align}
\notag&e_k[(b_1+\dotsm+b_{i-1})-(b_{j-i+k+1}+\dotsm+b_{j-1}+b_j)] \\
\label{E:coef rec}&= e_k[(b_1+\dotsm+b_{i-2})-(b_{j-i+k+1}+\dotsm+b_{j-1})] \\
\notag&+(b_{i-1} - b_j) e_{k-1}[(b_1+\dotsm+b_{i-2})-(b_{j-i+k-1}+\dotsm+b_j)].
\end{align}
Letting $A=(b_1+\dotsm+b_{i-2})-(b_{j-i+k+1}+\dotsm+b_{j-1})$ we have
\begin{align*}
e_k[A+b_{i-1}-b_j] &= \sum_{p+q=k} e_p[A] e_q[b_{i-1}-b_j] \\
&= e_k[A] + \sum_{q=1}^k e_q[b_{i-1}-b_j] e_{k-q}[A] \\
&= e_k[A] + \sum_{q=1}^k (b_{i-1}-b_j) (-b_j)^{q-1} e_{k-q}[A] \\
&= e_k[A] + (b_{i-1}-b_j) e_{k-1}[A-b_j]
\end{align*}
as required.
\end{proof}

\begin{ex}
Let $n=2$. Some of $y_{ij}$ are given below
\begin{eqnarray*}
y_{22}&=&z_{11}+(a_1-a_2)z_{12},\\
y_{23}&=&z_{12}+a_1z_{13},\\
y_{33}&=&z_{11}+a_1z_{12}+a_1a_2z_{13},\\
y_{34}&=&
z_{12}+(a_1+2a_2)z_{13}
+2(a_1+a_2)a_2z_{14},\\
y_{44} &=& z_{11}
+(a_{1} +a_{2}) z_{12}
+2( a_{1}+ a_{2} )a_2z_{13}+
2 (a_{1} + a_{2} )a_{2}^{2}z_{14}.
\end{eqnarray*}

\end{ex}

Then we eliminate the diagonal entries $z_{ii}\; (1\le i\le n)$ as follows. 
We set
\begin{equation}
\hz_i=z_{1,i+1}\quad (1\le i\le 2n).
\end{equation}
We define $w_{ij}\in S[\hz]:=S[\hz_{1},\dotsc,\hz_{2n}]$ by
\begin{align}
\label{E:w11}
    w_{11} &= 1 - (y_{n+1,n+1} - z_{11}) \\
\label{E:w11 in z-hat}
    &= 1 - a_1 z_{12} -a_1 a_2 z_{13} -\dotsm - a_1\dotsm a_n z_{1,n+1}\\
\label{E:wii def}
  w_{ii} &= -z_{11} + w_{11} + y_{ii}&&\text{for $2\le i\le 2n+1$} \\
  \label{E:wij def}
  w_{ij}&=y_{ij}&&\text{for $1 \le i<j\le 2n+1$.}
\end{align}
$w_{11}\in S[\hz]$ by \eqref{E:w11 in z-hat}.
For $2\le i\le 2n+1$, $w_{ii}\in S[\hz]$ by \eqref{E:wii def} and \eqref{E:diagonal y}.
For $1\le i<j\le 2n+1$, $w_{ij}\in S[\hz]$ by \eqref{E:leading term of y}.

\begin{rem}\label{R:yii to wii}
Since $z_{11}$ occurs in $y_{ii}$ with coefficient $1$, $w_{ii}$ is obtained from $y_{ii}$ by substituting $w_{11}$ for $z_{11}$.
\end{rem}

By \eqref{E:middle y}, \eqref{E:w11}, and \eqref{E:wii def} with $i=n+1$ we have
\begin{align}\label{E:middle w}
    w_{n+1,n+1} = 1.
\end{align}

\begin{rem}\label{R:wij} The $w_{ij}$ also satisfy the centralizer equations
\eqref{eq:EqtypeA}. 
If the relation does not involve diagonal entries, then this follows from \eqref{E:wij def} and Remark \ref{R:yz}. For $2\le i\le 2n$
\begin{align*}
  w_{i,i} - w_{i+1,i+1} + (b_i - b_{i+1}) w_{i,i+1} = y_{ii} - y_{i+1,i+1} + (b_i - b_{i+1}) y_{i,i+1}.
\end{align*}
Also
\begin{align*}
  w_{11} - w_{22} + (b_1-b_2) w_{1,2} &= w_{11} - (-z_{11}+w_{11}+y_{22}) + (b_1-b_2)w_{1,2} \\
  &= y_{11} - y_{22} + (b_1-b_2) y_{1,2}
\end{align*}
Both are relations by Remark  \ref{R:yz}.
\end{rem}

\begin{lem} \label{L:W} 
\begin{align}
    w_{ii} &\in 1 + \sum_{j=1}^{2n} S_+ \hz_{j} &\qquad&\text{for $1\le i\le 2n+1$} \\ 
 \label{eq:wi<j}
    w_{ij}&\in \hz_{j-i} + \sum_{k>j-i} S_+ \hz_{k} &&\text{for $1\le i<j\le 2n+1$.}
\end{align}
\end{lem}
\begin{proof} 
For $i<j$, this follows immediately from \eqref{E:leading term of y}. 
For $i=j$ this follows from sending $z_{11}\mapsto w_{11}$ in the expression \eqref{E:y to z1} for $y_{ii}$; see Remark \ref{R:yii to wii}.
\end{proof}

\begin{ex}
Let $n=2$. Using Lemma \ref{L:y} and \eqref{E:w11} we compute some of $w_{ij}$
\begin{align*}
w_{11} &= 1-a_1 \hz_{1} - a_1a_2\hz_{2},\\
w_{22} &
= 1 - a_2\hz_{1} -a_1a_2\hz_{2}, \\
w_{44}
&= 1+a_{2}\hz_{1}+a_{2}(a_{1}+2 a_{2}) \hz_{2}+2 a_{2}^{2}(a_{1}+a_2) \hz_{3},\\
w_{35}
&= \hz_{2} + 2(a_1+a_2) \hz_{3} + 2a_1(a_1+a_2)\hz_{4}, \\
 w_{45}
&= \hz_{1} +2 (a_1+a_2) \hz_{2} + 2 (a_{1}+a_2)^{2} \hz_{3}
+2 a_{1}^{2}(a_{1}+a_{2})\hz_{4}.
\end{align*}
\end{ex}

Let $(-1)^{n+1} r_{2i}$ be the $(n+1-i,n+1-i)$-th entry of the matrix $J - \trsp{Z}JZ$ for $1\le i\le n$.
Let $W$ be the $(2n+1)\times(2n+1)$ matrix with $(i,j)$-th entry $w_{ij}\in S[\hz]$. 
For $1\le i\le n$ let $(-1)^{n+1} \hat{R}_{2i} \in S[\hz]$ be the $(n+1+i,n+1+i)$-th entry of the matrix $J-{}^t{W}JW$.
The index indicates the degree of the element: $\deg(z_{ij})=j-i,\;\deg(a_i)=-1$.
Explicitly, we have
\begin{equation}
\hat{R}_{2i} =    w_{n+1,n+i+1}^2 +2\sum_{j=1}^i (-1)^j 
w_{n+j+1,n+i+1}\, w_{n-j+1,n+i+1}.\label{E:R into w} 
\end{equation}

There is a surjective $S$-algebra homomorphism 
\begin{align*}
S[z_{ij}\mid 1\le i\le j\le 2n+1]&\overset{\theta}{\longrightarrow} S[\hz_1,\dotsc,\hz_{2n}] \\
z_{ij} &\mapsto w_{ij}\qquad\text{for all $1\le i\le j\le 2n+1$.}    
\end{align*}

\begin{ex} For $n=2$ we have
\begin{align*}
\hat{R}_2 &=
\hz_{1}^{2}
-2 \hz_{2} 
+ 2\left(a_{1} +a_{2}\right) \hz_{1} \hz_{2} 
+ a_1\left(a_{1} + 2  a_{2}\right) \hz_{2}^{2} 
+ 2a_2\left(a_{1} + a_{2}\right) \hz_{1} \hz_{3}  \\
&+ 2a_1a_2\left(a_{1} + a_{2}\right) \hz_{2} \hz_{3}  - 2\left(a_{1} +a_{2}\right) \hz_{3}, \\
\hat{R}_4 &= 
\hz_{2}^{2} -2 \hz_{1} \hz_{3}+2\hz_{4}-2 a_{1}\hz_{1} \hz_{4} -2 a_{1} a_{2}\hz_{2} \hz_{4}. 
\end{align*}
\end{ex}

Define 
\begin{align}\label{E:Quo}
\Quo = S[\hz_1,\dotsc,\hz_{2n}]/(\hat{R}_2,\hat{R}_4,\dotsc,\hat{R}_{2n}).
\end{align}

\begin{thm} \label{T:iso} 
The map $\theta$ induces an $S$-algebra isomorphism $\C[\cZGv]\to \Quo$.
\end{thm}
\begin{proof} 
The centralizer relations \eqref{eq:EqtypeA} are sent to zero by $\theta$ by construction; see Remark \ref{R:wij}. 
By \eqref{E:middle w}, the relation \eqref{E:center 1} is sent to $0$ by $\theta$. 
Let $1\le i\le n$ and $k=n+1+i$. 
Recall that $R_{kk}$ is the $(k,k)$-entry of the matrix  $J-\trsp{Z}JZ$. Therefore $\theta(R_{kk})$ is obtained by taking the $(k,k)$-entry of $J - \trsp{W}JW$. 
The result is $(-1)^{n+1}\hat{R}_{2i}$, which is a relation in $\Quo$. 
By Proposition \ref{prop:Zeq}, we have verified all relations and the map is well-defined.
Define an $S$-algebra homomorphism $\eta:\Quo\to \C[\cZGv]$ by 
$\eta(\hz_i)=z_{1,i+1}$ for $1\le i\le 2n$. 
It is well-defined because $\eta(\hat{R}_{2i}) = r_{2i}$, which is in the ideal $\II$.
$\theta$ is surjective since each $w_{ij}$ is in its image and Lemma \ref{L:W} implies that all the $\hat{z}$-variables are in the $S$-span of the $w_{ij}$.

Finally $\theta$ and $\eta$ are mutually inverse, being mutually inverse on the $S$-algebra generators $z_{1,i+1}$ and $\hz_i$ for $1\le i\le 2n$.
\end{proof}

Let $\Delta=\prod_{\alpha\in\Phi_+}\alpha$ and $S_{\Delta}:=S[\Delta^{-1}].$ 
We denote $S^\reg\otimes_S \Quo$, and $S^\reg\otimes_S \C[\cZGv]$ by $\Quo^\reg,$ and $\C[\cZGv]^\reg$. 
\begin{prop}\label{P:w generates A delta} 
$\Quo^\Delta$ is generated by $w_{kk}$ for $1\le k\le 2n+1$, $k\ne n+1$, as an $S^\Delta$-algebra.
\end{prop}
\begin{proof} 
Since $w_{n+1,n+1}=1$, it suffices to show that the elements $\{w_{kk}-w_{11} \mid 2 \le k\le 2n+1\}$ are expressible as $S$-linear combinations of $\hz_1,\dotsc,\hz_{2n}$ by a triangular matrix, each of whose diagonal coefficients divide a power of $\Delta$.

By \eqref{E:diagonal y} for $2\le k\le 2n+1$ we have
\begin{align}
\label{E:wk}
  w_{kk}-w_{11} &= \dotsm + \prod_{j=1}^{k-1} (b_j-b_k) \hz_{k-1}.
\end{align}
where $\dotsm$ in \eqref{E:wk} means an $S$-linear combination of $\hz_p$ for $p<k-1$.
This is a clearly triangular system. It remains to show that the leading coefficient $c_k=\prod_{j=1}^{k-1}(b_j-b_k)$ is nonzero and divides a power of $\Delta$ for all $k$.
When $2\le k\le n+1$ we have $c_k=\prod_{j=1}^{k-1}(a_j-a_k)$ which divides $\Delta$.
Suppose $n+2\le k\le 2n+1$. Write $k=2n+2-i$ for some $1\le i\le n$. 
We have $b_k=-a_i$ and 
\begin{align*}
c_k &= \prod_{j=1}^{2n+1-i} (b_j-b_k) \\
&= (a_1+a_i)(a_2+a_i)\dotsm (a_i+a_i)\dotsm(a_n+a_i)(0+a_i)(-a_n+a_i)\dotsm (-a_{i+1}+a_i)
\end{align*}
which divides $\Delta^2$; the square is needed because of the factors $a_i+a_i$ and $0+a_i$.
\end{proof}


\subsection{Flatness of $\C[\cZGv]$ over $S$}
The goal of this subsection is to prove 
Theorem \ref{thm:flat}

By Theorem \ref{T:iso} we work with $\Quo$ instead of $\C[\cZGv]$.
The set of maximal ideals of $S$ is naturally in bijection with $\fh$. If $\mathfrak{m}$ is a maximal ideal of $S$, we say it is regular if the associated element of $\fh$ is regular.

\begin{cor} 
 \label{cor:regular wiis generates}
 If $\mathfrak{m}\in \fh$ is regular, then
 $\Quo^{\mathfrak{m}}=(S/\mathfrak{m})\otimes_S \Quo$ is 
 generated by $\{w_{ii}\mid \text{$1\le i\le 2n+1$, $i\ne n+1$} \}$ as an $(S/\mathfrak{m})$-algebra.
 \end{cor}
 \begin{proof}
 By Proposition \ref{P:w generates A delta} if $\mathfrak{m}$ is regular, then we can express each  $\hz_i$ as a linear combination of $w_{jj}$ for $1\le j\le 2n$.  
 Hence, the lemma follows. 
 \end{proof}

\begin{lem}
The leading term of $\hat{R}_{2i}$ is ${\hz_i}^2$ with respect to the graded reverse lexicographic (grevlex) term order (on monomials in the $\hat{z}_i$, ignoring $a_i$ variables) such that $\hz_1>\hz_2>\ldots>\hz_{2n}$.
   \end{lem}
\begin{proof} 
For $i<j$, by \eqref{eq:wi<j} and
\eqref{E:R into w}, 
the leading term of $w_{n+1,n+i+1}^2$ is $\hz_i^2,$ while for each term in the sum corresponding to $1\le j<i$, the leading monomial is $\hz_{i-j}\hz_{i+j}.$
Finally we consider the term $w_{n+i+1,n+i+1}w_{n-i+1,n+i+1}$ corresponding to $j=i$ in the sum of \eqref{E:R into w}.
Since the leading term of the second factor $w_{n-i+1,n+i+1}$ is
$\hz_{2i}$, the leading monomial of $w_{n+i+1,n+i+1}w_{n-i+1,n+i+1}$ is smaller than $\hz_i^2.$
\end{proof}

\begin{prop}\label{prop: fiber ring} Let
$\mathfrak{m}$ be a maximal ideal of $S$.
\begin{enumerate}
\item The polynomials $\hat{R}_{2i}$ specialised at $\mathfrak{m}$ form a Gr\"obner basis of the polynomial ring $(S/\mathfrak{m})[\hz_1,\ldots,\hz_{2n}]$ with respect to the graded reverse lexicographic order such that $\hz_1>\cdots>\hz_{2n},$ with $\hz_i$ of degree $1.$
\item The following set of monomials form a basis of $\Quom$ over $S/\mathfrak{m}\cong \C$:
    \begin{equation}
      \hz_1^{i_1}\cdots\hz_n^{i_n}
      \hz_{n+1}^{i_{n+1}}
      \cdots
      \hz_{2n}^{i_{2n}},
      \label{eq:monomial basis}
    \end{equation}
    where $i_k$ are non-negative integers such that
$i_k\in\{0,1\}$ for $1\le k\le n$.
In particular, $\{\hz_{n+1},\hz_{n+2},\ldots,\hz_{2n}\}$ is algebraically independent over $S/\mathfrak{m}$.
\item 
The ring $\Quom$ is a global complete intersection over $\C$ of dimension $n$.
\end{enumerate}
\end{prop}
\begin{proof} 
(1) Since the leading monomials of the $\hat{R}_{2i}$ are pairwise relatively prime, the $\hat{R}_{2i}$ form a Gr\"obner basis by Buchberger's algorithm. See \cite[Lemma 4.3]{IMN2} for a similar argument.

(2) is a consequence of (1).

(3) 
Let $B$ denote the subring of $\Quom$ generated by $\hz_{n+1},\ldots,\hz_{2n}$.
From (2), we see that $B$ is isomorphic to the polynomial ring of $n$ variables, therefore $\dim(B)=n$.
Because $\Quo^{\mathfrak{m}}$ is a finitely generated $B$-module by (2), we have $\dim (\Quom)=\dim B=n$.
Hence $\Quom$ is a global complete intersection.
\end{proof}

\begin{lem}\label{lem: Quo is flat}
    $\Quo$ is flat over $S$.
\end{lem}
\begin{proof}
By Proposition \ref{prop: fiber ring} (3), we can deduce that $\Quo$ is a Cohen-Macaulay ring.  
Since $S$ is regular, the dimension statement of (3) in Proposition \ref{prop: fiber ring} implies $\Quo$ is flat over $S$ (\cite[Corollary to Theorem 23.1]{Mat}).
\end{proof}
\begin{rem}
    In \cite{Pet}, it was proved that $\C[\cZGv]$ is flat over $S$ for general semisimple linear algebraic groups, using the sheaf Euler characteristic on the Peterson variety, which is a certain closed subvariety of $\fh\times (\Gv/\Bv)$, closely related to the Toda lattice of $\Gv$. We also refer  \cite{Rietsch} for more properties of the Peterson variety.
\end{rem}

The following fact is a corollary of the proof of Lemma \ref{lem: Quo is flat}.
\begin{cor}\label{cor: regular sequence}
$\hat{R}_{2},\hat{R}_{4},\ldots,\hat{R}_{2n} $ form a regular sequence in $S[\hz_1,\hz_2,\ldots,\hz_{2n}].$
\end{cor}

\section{Construction of $\hP_\la(y|a)$ by the Free-Fermion formalism
}\label{sec:Fermion}
In this section, we construct the equivariantly deformed Schur $P$- and $Q$-functions within the framework of the free-fermion formalism (cf. \cite{DJKM}). This approach provides a unified algebraic realization that naturally incorporates both the homological and cohomological theories. As an application,  we prove Lemma \ref{lem:hatQ in hatP}.

\subsection{Neutral free Fermions}

Let $\cC$ be the $\C[a]$-algebra generated by the neutral-fermions $\{\phi_n\}_{n\in \Z}$ with
\[
[\phi_m,\phi_n]_+=2(-1)^m\delta_{m+n,0},
\]
where $[x,y]_+=xy+yx.$
Note that $\phi_0^2=1$ and $\phi_n^2=0$ ($n\neq 0$).
There exists a $\C[a]$-linear anti-involution $\phi\mapsto \phi^*$ of $\cC$ defined by
\[
\phi_i\mapsto 
\phi_i^*:=(-1)^i\phi_{-i}.
\]
We have
\begin{align}\label{E:anticommutator star}
[ \phi_m^*,\phi_n]_+ &= 2 \delta_{m,n} \\
\label{E:anticommutator duality}
[\phi^*,\psi^*]_+ &= [\phi,\psi]_+\qquad\text{for $\phi,\psi\in \cC$.}
\end{align}

The Fock space $\cF$ is the left $\mathcal{C}$-module defined by $\cF = \cC /\sum_{n=1}^\infty \cC \phi_{-n}$.
By convention, the element $1\in\cF$ is called the \emph{vacuum vector} and is denoted by $\ket{\mathrm{vac}}$.
The dual Fock space $\cF^\dagger$ is the right $\cC$-module defined by $\mathcal{F}^\dagger=\mathcal{C}/\sum_{n=1}^\infty \phi_n\mathcal{C}$.
The element $1\in \mathcal{F}^\dagger$ is denoted by $\bra{\mathrm{vac}}$.
The map $\phi\mapsto \phi^*$ induces a $\C[a]$-linear isomorphism $\cF\to\cFd$. We denote the image of $\ket{v}\in\cF$ by $\bra{v^*}\in\cFd$.
It satisfies $a\ket{v}\mapsto \bra{v^*}a^*$.
There is a unique $\C[a]$-linear isomorphism $M:=\cC / (\sum_{n=1}^\infty \cC \phi_{-n}\oplus \sum_{n=1}^\infty \phi_n\cC ) \cong \C[a]$ sending $\phi_0\mapsto 1$. Let $\pi:\cC\to M\cong \C[a]$ be the natural projection. 
Consider the map $\cC\otimes_{\C[a]} \cC \to \C[a]$ defined by $a\otimes b\mapsto \pi(ab)$. It induces a nondegenerate $\C[a]$-bilinear pairing $\cFd \otimes_{\C[a]} \cF \to \C[a]$ denoted $\bra{v} \otimes \ket{w}\mapsto \bra{v} w \rangle$.
We have $(\bra{v}a)\ket{w}=\bra{v}(a\ket{w})$ for all $v\in\cFd$, $w\in\cF$, and $a\in \cC$. 
This common value is denoted $\langle v | a | w\rangle$. 

There is a vector space decomposition $\cF = \cF_+ \oplus \cF_-$ where $\cF_+$  (resp. $\cF_-$) is the $\C[a]$-subspace of $\cF$ spanned by all elements of the form $\phi_{i_1}\dots\phi_{i_p}\ket{\mathrm{vac}}$ with $p$ even (resp. odd).
The dual Fock space has a similar decomposition $\mathcal{F}^\dagger=\mathcal{F}^\dagger_{+}\oplus \mathcal{F}^\dagger_{-}$.
For $\la\in\SP$, let $r'$ denote the least even number satisfying $r'\geq \ell(\lambda)$. 
Define the vectors $\bra{\la}\in\mathcal{F}^\dagger$ and $\ket{\la}\in\mathcal{F}$ by 
\begin{align}
  \bra{\la} &= \bra{\mathrm{vac}}\phi^\ast_{\lambda_{r'}}\dots \phi^\ast_{\la_2}\phi^\ast_{\la_1}, \\
  \ket{\la} &= \phi_{\la_1}\phi_{\la_2}\dotsm \phi_{\la_{r'}} \ket{\mathrm{vac}}.
\end{align}
Then we have
\begin{align} 
\label{E:dual basis element}
\bra{\la} &= (\ket{\la})^* \\
  \bra{\la} \mu\rangle &= \delta_{\la\mu} 2^{\ell(\la)}.
  \label{eq:pair la mu} 
\end{align}

Define the operators $\phi^{(a)}_j$, $\hat{\phi}^{(a)}_j$ $(j\in \Z)$ acting on the Fock and dual Fock spaces by
\begin{align}
\phi^{(a)}_j&=
\begin{cases}
\sum_{i=0}^{j-1}
(-1)^ie_i(a_1,a_2,\dots,a_{j-1})\phi_{j-i}
 & (j> 0),\\
\sum_{i=0}^\infty
h_i(a_{1},a_2,\dots,a_{-j})\phi_{j-i}
& (j\leq 0),
\end{cases}\label{eq:def_phi_j}\\
\hat{\phi}^{(a)}_j&=
\begin{cases}
\sum_{i=0}^\infty
h_i(a_1,a_2,\dots,a_{j})\phi_{j+i}
 & (j\geq 0),\\
\sum_{i=0}^{-j-1}
(-1)^i
e_i(a_{1},a_2,\dots,a_{-j-1})\phi_{j+i}
& (j< 0),
\end{cases}\label{eq:def_hatphi_j}
\end{align}
where $e_i(a_1,a_2,\dots,a_p)=h_i(a_1,a_2,\dots,a_p)=\delta_{i,0}$ if $p=0$.
In particular, we have
\[
\bra{\mathrm{vac}}\phi^{(a)}_{n}=\bra{\mathrm{vac}}\hat{\phi}^{(a)}_{n}=0,\qquad
(\phi^{(a)}_{n})^\ast\ket{\mathrm{vac}}
=
(\hat{\phi}^{(a)}_{n})^\ast\ket{\mathrm{vac}}=0
\]
for any $n>0$.
\begin{lem}\label{lemma:duality_of_fermions}
We have the anti-commutation relations
\begin{align}\label{eq:duality_B}
&[(\phi^{(a)}_m)^\ast,\hat{\phi}^{(a)}_n]_+=2\delta_{m,n}, \\
\label{eq:duality_Bprime}
&[\phi^{(a)}_m,(\hat{\phi}^{(a)}_n)^\ast]_+=2\delta_{m,n}.
\end{align}
\end{lem}
\begin{proof} By \eqref{E:anticommutator duality} it suffices to show \eqref{eq:duality_B}.
(i) When $m=0$, the equation $$[(\phi^{(a)}_m)^\ast,\hat{\phi}^{(a)}_n]_+=
[\phi_0,\hat{\phi}^{(a)}_n]_+=2\delta_{0,n}
$$ follows directly from \eqref{eq:def_hatphi_j}.
(ii) When $m>0$, we have
\[
\begin{aligned}
& [(\phi^{(a)}_m)^\ast,\hat{\phi}^{(a)}_n]_+\\
&=
 \begin{cases}
 0 & (n<0),\\
 \sum_{i=0}^{m-1}\sum_{j=0}^\infty (-1)^ie_i(a_1,\dots,a_{m-1})
 h_j(a_1,\dots,a_{n})
 [\phi^\ast_{m-i},\phi_{n+j}]_+   
 &(n\geq 0)
 \end{cases}
\end{aligned}
\]
by \eqref{eq:def_phi_j} and \eqref{eq:def_hatphi_j}.
Obviously, \eqref{eq:duality_B} holds for $n<0$.
If $n\geq 0$, then the function $[(\phi^{(a)}_m)^\ast,\hat{\phi}^{(a)}_n]_+$ equals
\[
\begin{aligned}
&2\sum_{i=0}^{m-n}(-1)^ie_i(a_1,a_2,\dots,a_{m-1}) h_{m-n-i}(a_1,a_2,\dots,a_{n})    \\
&=
2h_{m-n}(a_1,a_2,\dots,a_n;a_1,a_2,\dots,a_{m-1})\\
&=
\begin{cases}
0 & (m<n),\\
(-1)^{m-n}2e_{m-n}(a_{n+1},a_{n+2},\dots,a_{m-1}) & (m\geq n)
\end{cases}\\
&=2\delta_{m,n},
\end{aligned}
\]
which implies \eqref{eq:duality_B}.
(iii) When $m<0$, we similarly have
\[
\begin{aligned}
\relax
[(\phi^{(a)}_m)^\ast,\hat{\phi}^{(a)}_n]_+
&=
\begin{cases}
0 & (n\geq 0),\\
2h_{m-n}(a_1,a_2,\dots,a_{-m};a_1,a_2,\dots,a_{-n-1})
& (n<0),
\end{cases}   \\
&=2\delta_{m,n},
\end{aligned}
\]
which implies \eqref{eq:duality_B}.
\end{proof}

As a corollary, we have
\begin{equation}\label{eq:0_th_value}
\bra{\mathrm{vac}}(\phi_0^{(a)})^\ast \hat{\phi}_0^{(a)}\ket{\mathrm{vac}}=
\bra{\mathrm{vac}}\phi^2_0\ket{\mathrm{vac}}=1.
\end{equation}
and
\begin{equation}\label{eq:n_th_value}
\bra{\mathrm{vac}}(\phi_m^{(a)})^\ast \hat{\phi}_n^{(a)}\ket{\mathrm{vac}}=
2\delta_{m,n}-
\bra{\mathrm{vac}}
\hat{\phi}_n^{(a)}
(\phi_m^{(a)})^\ast 
\ket{\mathrm{vac}}
=2\delta_{m,n}
\end{equation}
for $m,n>0$.

\subsection{$Q_\lambda(x)$ as a vacuum expectation value}
Let $x=(x_1,x_2,\ldots)$ be a sequence of variables. Following \cite{DJKM}, we introduce the operator $H_m$ defined by $$H_m=\frac{1}{4}\sum_{i\in \Z}(-1)^i\phi_{-i-m}\phi_i$$ for each odd integer $m$. Each $H_m$ gives a well-defined endomorphism of $\mathcal{F}$, and satisfies the relations $[H_m,\phi_n]=\phi_{n-m}$.
We define the \emph{Hamiltonian} $H(x)$ by 
\[
H(x):=2\sum_{m\in 1 + 2 \Z_{\ge0}}\frac{p_m(x)}{m}H_m,
\]
where $p_m(x)=\sum_{i=1}^\infty x_i^m$ is the $m$-th powersum in $x.$
The operator $H(x)$ is well-defined as a linear map from $\mathcal{F}$ to $\Gamma(x)\otimes \mathcal{F},$ where $\Gamma(x)=\C[p_1(x),p_3(x),\ldots].$
By using
$H(x)\ket{\mathrm{vac}}=0$, one sees that $\bra{u}e^{H(x)}\ket{v}$ is a well-defined element of $\Gamma(x)$ for any $\bra{u}\in \mathcal{F}^\dagger$ and $\ket{u}\in \mathcal{F}$. 
In particular, for $\la\in \SP$, it follows from \cite[\S 6]{JM},\cite{Y} that 
\begin{equation}
Q_\la(x)=\bra{\mathrm{vac}}e^{H(x)}\ket{\la}.
\label{eq:Q}
\end{equation}
Historically, the function $Q_\lambda(x)$ appeared as the  $\tau$-function of the BKP-hierarchy \cite{DJKM}.
\subsection{$\hat{P}$-functions as vacuum expectation value}
We will deform the vector $\ket{\la}$ so that $\hat{P}_\lambda(y|a)$ and $Q_\lambda(x|a)$ are expressed in similar forms. 

Let $\lambda=(\la_1,\ldots,\la_{r})\in \SP$ where $r=\ell(\lambda)$ is the number of nonzero parts of $\lambda$.
Define $\ket{\lambda}_H,\ket{\lambda}_C\in \mathcal{F}_{+}$ and
$\bra{\lambda}_H,\bra{\lambda}_C\in \mathcal{F}^\dagger_{+}$ by
\[
\begin{gathered}
\ket{\lambda}_{H}=
\hat{\phi}_{\lambda_1}^{(a)}\dots \hat{\phi}_{\lambda_{r'}}^{(a)} \ket{\mathrm{vac}},
\quad
{}_C\bra{\lambda}=
\bra{\mathrm{vac}}(\phi^{(a)}_{\lambda_{r'}} )^\ast \dots (\phi^{(a)}_{\lambda_1} )^\ast,\\
\ket{\lambda}{}_C=
{\phi}_{\lambda_1}^{(a)}\dots {\phi}_{\lambda_{r'}}^{(a)} \ket{\mathrm{vac}}, 
\quad
{}_H\bra{\lambda}=
\bra{\mathrm{vac}}(\hat{\phi}^{(a)}_{\lambda_{r'}} )^\ast \dots (\hat{\phi}^{(a)}_{\lambda_1} )^\ast.
\end{gathered}
\]
The notation $H$ and $C$ means homology and cohomology, respectively.
Then, from Lemma \ref{lemma:duality_of_fermions}, we have
\begin{equation}\label{eq:duality_of_vectors}
{}_H\langle{\lambda}\ket{\mu}_C=2^{\ell(\lambda)}\delta_{\lambda\mu}
,\quad
{}_C\langle{\lambda}\ket{\mu}_H=2^{\ell(\lambda)}\delta_{\lambda\mu}.
\end{equation}
Note that 
$\ket{\mathrm{vac}}=\ket{\mathrm{vac}}_C=\ket{\mathrm{vac}}_H,\; \bra{\mathrm{vac}}={}_C\bra{\mathrm{vac}}={}_H\bra{\mathrm{vac}}.$


The following is a basic fact we use.
\begin{prop}\label{prop:isometry}
For $\ket{u},\ket{v}\in \mathcal{F}_+,$
define $f(x)=\bra{\mathrm{vac}}e^{H(x)}\ket{u},$ $g(y)=\bra{\mathrm{vac}}e^{H(y)}\ket{v}$. Then
\begin{equation}\langle f(x),g(y)\rangle=
\langle{u^*}|{v}\rangle
\end{equation}
where the pairing in the left-hand side is defined in \S \ref{SS:Schur Q}.
\end{prop}
\begin{proof}
We write 
$u=\sum_{\lambda}u_\lambda \ket{\la}$ and $v=\sum_{\lambda}v_\lambda \ket{\la}$ for $u_\la,v_\la\in \C[a]$. By \eqref{E:dual basis element} we have $u^* = \sum_\la u_\la \bra{\la}$.
Then by \eqref{eq:Q} and \eqref{eq:pair la mu},
\begin{align*}
 \langle f(x),g(y)
 \rangle
 &=
\sum_{\la,\mu}
 u_\la v_\mu\langle Q_\la(x),Q_\mu(y)\rangle=\sum_{\la,\mu}
 u_\la v_\mu
2^{\ell(\lambda)}\delta_{\la\mu}\\
 &=\sum_{\la,\mu}
 u_\la v_\mu
\langle\lambda|
\mu\rangle=\langle u^*|v\rangle.
\end{align*}
\end{proof}
\begin{thm}
For a strict partition $\lambda$, we have
\[
\hat{P}_\lambda(y|a)
=
2^{-\ell(\lambda)}\bra{\mathrm{vac}}e^{H(y)}\ket{\lambda}_H.
\]
\end{thm}
\begin{proof}
Let $\la\in \SP$, and we denote the function $\bra{\mathrm{vac}}e^{H(y)}\ket{\lambda}_H$ by $F_\la(y|a).$ We expand it as $F_\la(y|a)=\sum_{\mu\in \SP} C_{\lambda\mu}P_\mu(y).$
From \eqref{eq:Q} and $\langle Q_\lambda(x),P_\mu(y)\rangle=\delta_{\lambda,\mu}$, we have
\begin{align*}
C_{\lambda\mu}
&=
\langle Q_\mu(x),F_\la(y|a)\rangle\\
&=\langle \mu |\lambda\rangle_H\quad\text{(by Proposition \ref{prop:isometry})}\\
&=
\bra{\mathrm{vac}}\phi^\ast_{\mu_s}\dots \phi^\ast_{\mu_1}\hat{\phi}^{(a)}_{\lambda_1}\dots \hat{\phi}^{(a)}_{\lambda_{r'}}\ket{\mathrm{vac}},
\end{align*}
where $r,s$ are the least even integer satisfying $r\ge \ell(\la)$ and $s\ge \ell(\mu)$ respectively.
This vacuum expectation value is computed by Wick's theorem (see \cite{DJKM}).
When $r\neq s$, $C_{\lambda\mu}$ vanishes automatically.
When $r=s$, we have
\[C_{\lambda\mu}=\det(M),\] where $M=(
\bra{\mathrm{vac}}
\phi^\ast_{\mu_j}
\hat{\phi}^{(a)}_{\lambda_i}
\ket{\mathrm{vac}}
)_{1\leq i,j\leq r}
$.
Since $$\bra{\mathrm{vac}}
\phi^\ast_{\mu_i}
\hat{\phi}^{(a)}_{\lambda_j}
\ket{\mathrm{vac}}=
\begin{cases}
2h_{\mu_i-\lambda_j}(a_1,a_2,\dots,a_{\lambda_j}) & (\lambda_j>0)\\
\delta_{\mu_i,0} & (\lambda_j=0)
\end{cases},
$$
we have 
\[
C_{\lambda\mu}=\begin{cases}
    2^{\ell(\lambda)}\det\left(h_{\mu_i-\lambda_j}(a_1,a_2,\dots,a_{\lambda_j})\right)&\text{if $\la\subset \mu$ and $\ell(\la)=\ell(\mu)$},\\
    0&\text{otherwise.}
\end{cases}
\]
Hence, from Proposition \ref{prop: dual Ivanov},
we obtain the desired result.
\end{proof}

\begin{prop}\label{prop:dual C H} For $\la,\mu\in \SP$, we have
 $2^{-\ell(\lambda)}{}_C\langle \lambda|\mu \rangle_H =\delta_{\lambda\mu}$   
\end{prop}
\begin{proof}
By Wick's theorem, ${}_C\langle \lambda|\mu \rangle_H$ vanishes unless $\ell(\lambda)=\ell(\mu)$.
Assume $\ell(\lambda)=\ell(\mu)$.
Let $r'$ be the smallest integer satisfying $r'\geq \ell(\lambda)$.
Then, by Wick's theorem, we obtain
\[
\begin{aligned}
{}_C\langle \lambda|\mu \rangle_H&=
\bra{\mathrm{vac}}
(\phi^{(a)}_{\lambda_{r'}})^\ast
\dots
(\phi^{(a)}_{\lambda_{1}})^\ast
\hat\phi^{(a)}_{\mu_1}\dots
\hat\phi^{(a)}_{\mu_{r'}}\ket{\mathrm{vac}}\\
&=\det\left(\bra{\mathrm{vac}}
(\phi^{(a)}_{\lambda_{i}})^\ast\hat\phi^{(a)}_{\mu_j}
\ket{\mathrm{vac}}\right)_{1\leq i,j\leq r'}.
\end{aligned}
\]
Then we obtain ${}_C\langle \lambda|\mu \rangle_H=2^{\ell(\lambda)}\delta_{\lambda\mu}$ by \eqref{eq:0_th_value} and Lemma \ref{lemma:duality_of_fermions}.
\end{proof}

\begin{cor}[cf. \cite{ISY}]
For $\la\in \SP$, we have
\[
Q_\lambda(x|a)
=
\bra{\mathrm{vac}}e^{H(x)}\ket{\lambda}_C.
\]
\end{cor}
\begin{proof}
This follows from Proposition \ref{prop:dual C H} and Proposition \ref{prop:isometry}.
\end{proof}
\subsection{Expansion of $\hat{q}$-function in terms of $\hat{P}$-functions}\label{ssec:proof of hat q expansion}

\begin{lem}\label{lemma:Q_to_hatP}
For $i\geq 0$, we have 
\[
Q_i(y)=\delta_{i,0}+2\sum_{k=0}^\infty (-1)^{k}e_{k}(a_1,a_2,\dots,a_{i+k-1})\hat{P}_{i+k}(y|a).
\]
\end{lem}
\begin{proof}
When $i=0$, the claim is trivial.
Assume $i>0$.
From \eqref{eq:duality_of_vectors}, we have
\[
\mathrm{id}_{\mathcal{F}_+}=\sum_{\lambda:\mathrm{strict}} 2^{-\ell(\lambda)}\ket{\lambda}_H\cdot {}_C\bra{\lambda}.
\]
Then, we have
\begin{equation}\label{eq:Q_to_vac}
Q_i(y)=\bra{\mathrm{vac}}e^{H(y)}\phi_i\phi_0\ket{\mathrm{vac}}=
\sum_{\lambda}2^{-\ell(\lambda)}
\bra{\mathrm{vac}}e^{H(y)}
\ket{\lambda}_H\cdot 
{}_C\bra{\lambda}
\phi_i\phi_0\ket{\mathrm{vac}}.
\end{equation}
As $i>0$, 
$
{}_C\bra{\lambda}\phi_i\phi_0\ket{\mathrm{vac}}
$ vanishes unless $\ell(\lambda)=1$.
When $\ell(\lambda)=1$, we have $\lambda=(n)$ for some $n>0$, and
\[
\begin{aligned}
{}_C\bra{(n)}\phi_i\phi_0\ket{\mathrm{vac}}
&=
\bra{\mathrm{vac}}\phi_0\phi_i^\ast \ket{(n)}_C\\
&=
\bra{\mathrm{vac}}\phi_0\phi_i^\ast\phi^{(a)}_n\phi_0\ket{\mathrm{vac}}\\
&=
\sum_{j=0}^{n-1}
(-1)^je_j(a_1,a_2,\dots,a_{n-1})
\bra{\mathrm{vac}}\phi_0\phi_i^\ast\phi_{n-j}\phi_0\ket{\mathrm{vac}}\\    
&=(-1)^{n-i}2e_{n-i}(a_1,a_2,\dots,a_{n-1}).
\end{aligned}
\]
By substituting this equation into \eqref{eq:Q_to_vac}, we complete the proof.
\end{proof}

Recall the $\hat{q}$-functions defined in \eqref{eq:def_of_hatq}.
Lemma \ref{lem:hatQ in hatP} is a direct consequence of the following.

\begin{prop}\label{prop:hatq_to_hatP} 
Let $c_i$ be arbitrary parameters. 
Then, we have
\begin{align}
\label{eq:hat q in hat P}
\hat{q}_i(y|c_1,\dots,c_m)
=\delta_{i,0}
+2\sum_{k=0}^\infty h_{k}(c_1,c_2,\dots,c_m;a_1,a_2,\dots,a_{i+k-1})\cdot \hat{P}_{i+k}(y|a).
\end{align}
\end{prop}

\begin{proof}
By \eqref{eq:def_of_hatq} and Lemma \ref{lemma:Q_to_hatP}, we have
\begin{align*}
&\hat{q}_i(y|c_1,\dots,c_m)\\
&=
\sum_{k=0}^\infty h_k(c_1,c_2,\dots,c_m)
Q_{i+k}(y)\\
&=
\sum_{k= 0}^\infty h_k(c_1,c_2,\dots,c_m)
\left(\delta_{i+k,0}+2\sum_{l=0}^\infty (-1)^{l}e_{l}(a_1,a_2,\dots,a_{i+k+l-1})\hat{P}_{i+k+l}(y|a)\right)\\
&=\delta_{i,0}+2\sum_{l=0}^\infty\sum_{k=0}^\infty
(-1)^{l}h_k(c_1,c_2,\dots,c_m)e_{l}(a_1,a_2,\dots,a_{i+k+l-1})\hat{P}_{i+k+l}(y|a)\\
&=\delta_{i,0}+2\sum_{j=0}^\infty
h_{j}(c_1,c_2,\dots,c_m;a_1,a_2,\dots,a_{i+j-1})
\hat{P}_{i+j}(y|a).
\end{align*}
\end{proof}


\begin{thebibliography}{99}

 \bibitem{AJS}
 H.~H.~Andersen, J.~C.~Jantzen, and W.~Soergel.
 \emph{Representations of quantum groups at a pth root of unity and of semisimple groups in characteristic p: independence of p,} Ast\'erisque (1994), no.220, 321 pp.



 \bibitem{Bi} S.~Billey,
 \emph{Kostant polynomials and the cohomology ring for $G/B$,}
 Duke Math. J. 96 (1999), no. 1, 205–224.

\bibitem{CP} P.~E.~Chaput and N.~Perrin,
\emph{Affine symmetries in quantum cohomology: corrections and new results,}
Math. Res. Lett. 30(2):341-374


\bibitem{DJKM}
E.~Date, M.~Jimbo, M.~Kashiwara, T.~ 
Miwa, \emph{Operator Approach to the Kadomtsev-Petviashvili Equation–Transformation Groups for Soliton Equations III,} 
J. Phys. Soc. Japan, November 15, 1981, Vol. 50, No. 11: pp. 3806-3812


\bibitem{DJKM4}E.~Date, M.~ Jimbo, M.~Kashiwara, T.~ 
Miwa, 
\emph{Transformation groups for soliton equations: IV. New hierarchy of soliton equation of $KP$-type,} Physica 4
D, {\bf 4} (1982) 343--365. 

\bibitem{EE} 
H. Eriksson and K. Eriksson. Affine Weyl groups as infinite permutations. Elec. J. Combin.
5 (1998), \#R18.



\bibitem{FH} W.~Fulton and J.~Harris, \emph{Representation Theory --- A First Course}, 
Graduate Texts in Mathematics 129 (1991) Springer. 


\bibitem{Gi} V.~Ginzburg.
\emph{Perverse sheaves on a Loop group and Langlands’ duality,} preprint, 1995;
arXiv:alg-geom/9511007.


\bibitem{GZ}
J.~Guo and H.~Zou, \emph{Quantum cohomology of symplectic flag manifolds}, J. Phys. A: Math. Theor. 55 (2022) 275401

\bibitem{HJ} C. Hanusa and B. Jones,
\emph{Abacus models for parabolic quotients of affine Weyl groups,}
J. Alg. {\bf 361} (2012), 134–162.

\bibitem{Ike}
T.~Ikeda, \emph{Schubert classes in the equivariant cohomology of the Lagrangian Grassmannian,} Adv. Math. {\bf 215} (2007), 1--23.


\bibitem{IINY}
T. Ikeda, S. Iwao, S. Naito, K. Yamaguchi,
\emph{Relativistic Toda lattice and equivariant $K$-homology},
arXiv.2050.02941



\bibitem{IMN} T.~Ikeda, L.~Mihalcea, and H.~Naruse,
\emph{Double Schubert polynomials for the classical groups,}
Adv. Math. {\bf 226} (2011), no. 1, 840–886.

\bibitem{GPGQ} T.~Ikeda and H.~Naruse,
\emph{$K$-theoretic analogues of factorial Schur $P$- and $Q$-functions}, Adv. Math. {\bf 243} (2013), 22–66.

 \bibitem{IMN2}
 T.~Ikeda, L.~C.~Mihalcea, H.~Naruse, \emph{Factorial $P$- and $Q$-Schur functions represent equivariant quantum Schubert classes},
 Osaka J. Math. 53 (2016), no. 3, 591–619.




\bibitem{ISY}
T.~Ikeda, M.~Shimozono, and K.~Yamaguchi,
\emph{Equivariant $K$-homology of affine Grassmannian and $K$-theoretic double $k$-Schur functions}, 
arXiv:2408.10956

\bibitem{ISY2}
T.~Ikeda, M.~Shimozono, K.~Yamaguchi,
\emph{Affine double $gp$-functions and Equivariant $K$-homology of type C affine Grassmannian,}
in preparation.

 \bibitem{Ivanov-dim}
 V.~N.~Ivanov, 
 \emph{The dimension of skew shifted Young diagrams, and projective characters of the infinite symmetric group},
 Zap. Nauchn. Sem. S.-Peterburg. Otdel. Mat. Inst. Steklov. (POMI) 240 (1997), 115–135, 292–293.
 J. Math. Sci. (New York) 96 (1999), no. 5, 3517–3530.

\bibitem{Iv} V.~N.~Ivanov,
\emph{Interpolation analogues of Schur $Q$-functions,}(English summary)
Zap. Nauchn. Sem. S.-Peterburg. Otdel. Mat. Inst. Steklov. (POMI) 307 (2004), Teor. Predst. Din. Sist. Komb. i Algoritm. Metody. 10, 99–119, 281–282; reprinted in
J. Math. Sci. (N.Y.) {\bf 131} (2005), no. 2, 5495–5507
\bibitem{Iv:pri} V.~N.~Ivanov, private communication (2023).


\bibitem{JM} M.~Jimbo and T.~Miwa,
\emph{Solitons and infinite-dimensional Lie algebras,}
Publ. Res. Inst. Math. Sci. 19 (1983), no. 3, 943–1001.


\bibitem{Ko}
B.~Kostant.
\emph{The solution to a generalized Toda lattice and representation theory,}
Adv. in Math. {\bf 34} (1979), no. 3, 195–338.
\bibitem{KK} B.~Kostant and S.~ Kumar, \emph{The Nil Hecke Ring and Cohomology of $G/P$
for a Kac-Moody Group $G$},
Adv.  Math {\bf 62} (1986) 
187--237.



 \bibitem{Lam} T.~Lam,
\emph{Schubert polynomials for the affine Grassmannian},
 J. Amer. Math. Soc. {\bf 21} (2008), no. 1, 259–281.


\bibitem{LLMSSZ} T.~Lam, L.~Lapointe, J.~Morse, A.~Schilling, M.~Shimozono, and M.~Zabrocki.
\emph{$k$-Schur functions and affine Schubert calculus}, Fields Institute Monographs, 33. Springer, New York; Fields Institute for Research in Mathematical Sciences, Toronto, ON, 2014. viii+219 pp.

\bibitem{LLS:backstable} T. Lam, S.~J.~Lee, and M. Shimozono.
\emph{Back stable Schubert calculus}, 
Compos. Math. {\bf 157} (2021), no. 5, 883–962.

\bibitem{LLS:coproduct}
T.~Lam, S.~J.~Lee, and M.~Shimozono,
\emph{On the coproduct in affine Schubert calculus}, Facets of algebraic geometry. Vol. II, 115–146.
London Math. Soc. Lecture Note Ser., {\bf 473} 
Cambridge University Press, Cambridge, 2022.




 \bibitem{LSS:C} T.~Lam, A.~Schilling, and M.~Shimozono,
\emph{Schubert polynomials for the affine Grassmannian of the symplectic group},
Math. Z. {\bf 264} (2010), no. 4, 765–-811. 

\bibitem{LSS:K} T.~Lam, A.~Schilling, and M.~Shimozono,
\emph{$K$-theory Schubert calculus of the affine Grassmannian},
Comp. Math. {\bf 146} (2010), no. 4, 811--852. \bibitem{LS:Acta} T.~Lam and M.~ Shimozono, \emph{Quantum cohomology of $G/P$ and homology of affine Grassmannian},
Acta Math. {\bf 204} (2010), no. 1, 49--90.

\bibitem{LS:double Kostant} T.~Lam and M.~Shimozono, 
\emph{From double quantum Schubert polynomials to k-double Schur functions via the Toda lattice},
arXiv:1109.2193.

\bibitem{LS:pieri} T.~Lam and M.~ Shimozono,
\emph{Equivariant Pieri rule for the homology of the affine Grassmannian},
J. Algebraic Combin. {\bf 36} (2012), no. 4, 623–648. 

\bibitem{LS:double Schur}
T.~Lam and M.~Shimozono, \emph{$k$-Double Schur functions and equivariant (co)homology of the affine Grassmannian}, Math. Ann. (2013) {\bf 356}, 1379--1404. 

\bibitem{LaSc}  A.~Lascoux and M.-P.~Sch\"utzenberger, Polyn\^omes de Schubert. C. R. Math. Acad. Sci. Paris, S\'er. I Math., \textbf{294} (1982) no. 13, 447–450. 

\bibitem{LM}L.~Lapointe and J.~Morse, \emph{A $k$-Tableau Characterization of $k$-Schur Functions}, Adv. Math. {\bf 213} (2007) 183--204.

\bibitem{Mac91}
I. G. Macdonald,
Notes on Schubert Polynomials,
Publications du LaCIM, Laboratoire de combinatoire et d'informatique math\'ematique
(LaCIM),
Universit\'e du Qu\'ebec \`a Montr\'eal,
1991,{\bf 6}.

\bibitem{Mac} I.~G.~Macdonald,
\emph{Symmetric Functions and Hall-polynomials},
Oxford
Mathematical Monographs, Clarendon Press, Oxford, 1979, viii + 180 pp.

\bibitem{Mac:aff}
I.~G.~Macdonald, \emph{Affine Hecke Algebras and Orthogonal Polynomials}, Cambridge University Press, 2003.



\bibitem{Mat}
H.~Matsumura,
\emph{Commutative ring theory},
Translated from the Japanese by M. Reid. Second edition
Cambridge Stud. Adv. Math., 8
Cambridge University Press, Cambridge, 1989. xiv+320 pp.


\bibitem{Molev}A.~ I.~Molev, 
\emph{Comultiplication rules for the double Schur functions and Cauchy identities}, Electron. J. Combin. {\bf 16} (2009), no. 1, Research Paper 13, 44 pp.

\bibitem{NN:Uni}
M.~Nakagawa and H.~Naruse,
\emph{Universal Factorial Schur $P,Q$-functions and their duals}, arXiv.1812.03328v1.

\bibitem{NN23}M.~Nakagawa and H.~Naruse, \emph{The universal factorial Hall-Littlewood $P$- and $Q$-functions}, Fund. Math. {\bf 263} (2023), 133--166.

\bibitem{NN}
M.~Nakagawa and H.~Naruse, 
\emph{Generalized (co)homology of the loop spaces of classical groups and the universal factorial Schur $P$- and $Q$-functions}, (English summary) Schubert calculus—Osaka 2012, 337--417,
Adv. Stud. Pure Math., {\bf 71}, Math. Soc. Japan, Tokyo, 2016.

\bibitem{Pet} D. Peterson,
MIT lecture notes, 1997.

\bibitem{Pon} S.~Pon, \emph{Affine Stanley symmetric functions for classical types}, J. Alg. Combin. {\bf 36}(4), (2012). DOI 10.1007/s10801-012-0352-6

\bibitem{Rietsch}
K.~Rietsch, \emph{A mirror symmetric construction of $qH_T^*(G/P)_{(q)}$}, Adv. Math. {\bf 217} (2008) 2401--2442.

\bibitem{Y}
Y.~You, \emph{Polynomial solutions of BKP-hierarchy and projective representation of symmetric groups}, in { Infinite-dimensional Lie algebras} (Luminy-Marseille, 1988) Adv. Ser. Math. Phys. {\bf 7}, 1989, pp 449–464, World Science, Singapore, 1989.

\bibitem{YZ} Z.~Yun and X.~Zhu, 
\emph{Integral homology of loop groups via Langlands dual groups}, 
Represent. Theory {\bf 15} (2011), 347–369.

\end{thebibliography}
\end{document}